\documentclass[reqno]{amsart}     
\usepackage{amssymb}
\usepackage{amsmath,amssymb,amsfonts,amsthm,latexsym,amscd,epsfig,color}
\usepackage{mathrsfs}
 
\usepackage{a4wide} 
\usepackage[T1]{fontenc}
\usepackage[latin1]{inputenc}
\usepackage{pgf,pgflibraryshapes,tikz,yfonts}

\usepackage{anyfontsize}
\usepackage{t1enc}

\renewcommand{\(}{\left(}
\renewcommand{\)}{\right)}
\newtheorem{theo}{Theorem}
\newtheorem{prop}{Proposition}
\newtheorem{lemma}{Lemma}
\newtheorem{cor}{Corollary\!\!}

\newtheorem{ncor}{Corollary}

\theoremstyle{definition}

\newtheorem{df}{Definition}
\newtheorem{ex}{Example}

\theoremstyle{remark}

\newtheorem{rem}{Remark\!\!}
 
\newtheorem{nrem}{Remark}

\newcommand{\rcp}[1]{\frac{1}{#1}}

\newcommand{\bt}{\begin{theo}} 
\newcommand{\et}{\end{theo}} 
\newcommand{\bl}{\begin{lemma}} 
\newcommand{\el}{\end{lemma}} 
\newcommand{\bp}{\begin{prop}} 
\newcommand{\ep}{\end{prop}} 
\newcommand{\bdf}{\begin{df}} 
\newcommand{\edf}{\end{df}} 
\newcommand{\brem}{\begin{rem}} 
\newcommand{\erem}{\end{rem}} 
\newcommand{\bnrem}{\begin{nrem}} 
\newcommand{\enrem}{\end{nrem}} 
\newcommand{\bex}{\begin{ex}} 
\newcommand{\eex}{\end{ex}} 
\newcommand{\bcor}{\begin{cor}} 
\newcommand{\ecor}{\end{cor}} 
\newcommand{\bncor}{\begin{ncor}} 
\newcommand{\encor}{\end{ncor}} 
\newcommand{\bpf}{\begin{proof}} 
\newcommand{\epf}{\end{proof}}

\newcommand{\nti}{{n\to\infty}}

\begin{document}

\title{Counting  General Phylogenetic networks  } 
\author{Marefatollah Mansouri} 
\thanks{This research has been supported by bilateral Austrian-Taiwanese project FWF-MOST, grant
I~2309-N35.}

\address{Department of Discrete Mathematics and Geometry, Technische
Universit\"at Wien, Wiedner Hauptstra\ss e 8-10/104, A-1040 Wien, Austria.}
\email{marefatollah.mansouri@tuwien.ac.at}




\begin{abstract} 
We provide precise asymptotic estimates for the number of  general phylogenetic networks by using analytic combinatorial methods. Recently, this approach is studied by \textit{Fuchs, Gittenberger} and author himself (\textit{Australasian Journal of Combinatorics 73(2):385-423, 2019}), to count networks with few reticulation vertices for two subclasses: \textit{tree-child}  and \textit{normal} networks. We follow this line of research to show how to  obtain  results on the enumeration of  general phylogenetic networks.    

\end{abstract} 
 
\maketitle

\section{Introduction} 
A phylogenetic network is a generalization of a phylogenetic tree which can be
used to describe the evolutionary history of a set of species that is non-tree like. Phylogenetic trees are usually computed from molecular sequences that are commonly used to show evolutionary history. Phylogenetic trees provide a useful representation of many
evolutionary relationships, and have been well studied (see, for example \cite{huson2010phylogenetic, BoSe16a, Se16, SeSt06, Sc}). However, these trees are less suited to model mechanisms of \textit{reticulate evolution} \cite{10.2307/2412721}, such as hybridization, recombination, or
reassortment. Phylogenetic networks provide an alternative to phylogenetic trees when analyzing data sets whose evolution involves of reticulate events (for more details see, \cite{Bandelt1994,GEL2003, LinderRieseberg2004}).

In the literature, there are often two kinds of labeling for phylogenetic networks
:  \textit{leaf-labeled }and \textit{vertex-labeled}. In the latter case, all the vertices take different labels, and in the former case, leaves are labeled but non-leaf vertices are unlabeled. 
The importance of phylogenetic networks is because of precise  representations of reticulation event which is in particular, the formation of hybrid species. However,  there are particular principles in the procedure of evolution that cause some more restrictions on phylogenetic networks. Thus,
biologists have defined many subclasses of the class of phylogenetic networks.
Recently, studying of enumerative aspects of phylogenetic networks and related structures have become increasingly interested. We mentioned already the shape analysis of phylogenetic trees
\cite{Bo16, BoFl09, FoRo80, FoRo88} and the bounds for the counting sequences of some classes of
phylogenetic networks \cite{MSW15}. But other counting problems were studied in
\cite{AlAl16,CEJM13,DiRo17,LSS13,Ro07,Se17,SeSt06, unpublishedkeyA, unknown, MPM, tree-child, Scornavacca}. Though combinatorial counting problems are
often amenable to the rich tool box of analytic combinatorics \cite{AnaCombi}, generating
functions have been rarely used in phylogenetic networks enumeration problems.

There are a quite few research studies on  general phylogenetic networks. This paper is concerned with the counting of general phylogenetic networks, a basic and fundamental question which is of interest in mathematical biology \cite{MSW15}.
On the other
hand, the combinatorial view of phylogenetic networks as an interesting challenge has been addressed only in few papers. The goal of this study is to develop a more rigorous understanding of counting problem for general phylogenetic networks with a fixed number of reticulation events. 
Here we come back to the open problems of \cite{unpublishedkeyA}   which are left for general networks and   show that \textit{sparsened skeleton decomposition} is powerful tool for enumeration problems in general phylogenetic networks. The purpose of the current study is to show how the present method in \cite{unpublishedkeyA} for \textit{tree-child} networks can be extended for general networks.    
 Before stating our results in more detail, we recall some definitions and previous works. 
 A phylogenetic network is defined as a directed acyclic graph (DAG) which is
 connected and consists of the following vertices:\\
 $ (1) $ \textit{leaves}  which have out-degree 0 and in-degree 1;\\
 $ (2) $ \textit{tree vertices} have out-degree 2 and in-degree 1;\\
 $ (3) $ \textit{reticulation vertices} have out-degree 1 and in-degree 2;\\
 $ (4) $ and the \textit{root node} with out-degree 2 and in-degree 0.\\
 Also,  a phylogenetic network is called \textit{tree-child} network if for every non-leaf node at least one of its children is a tree node or a leaf. Equivalently,  every tree vertex must have at least one child
 which is not a reticulation vertex and every reticulation vertex is not directly followed by another reticulation vertex.\\
 Note that variations on the definition of rooted binary phylogenetic networks are around in the
 literature. In general phylogenetic networks, as defined before, multiple edges are not explicitly forbidden. Our goal is indeed to study the most general model of
   general  phylogenetic networks that could be counted if their number of reticulation vertices is fixed and provide a more detailed investigation regarding enumeration properties of general networks with or without multiple edges on their structures.

\begin{figure}[h]
	\includegraphics[scale=0.2]{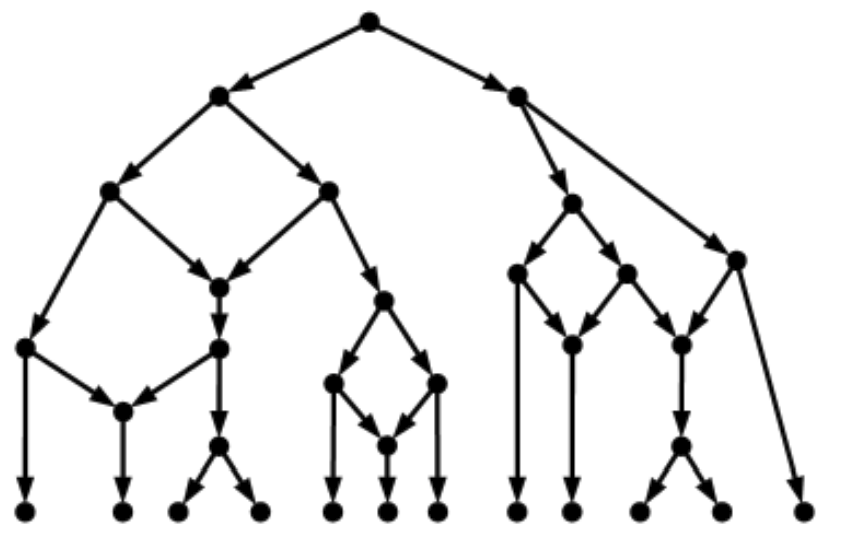}
	\caption{A phylogenetic network which is not a tree-child network.} \label{fig_skeleton}
\end{figure}
Now, denote by $ G_{k,n} $ resp. $ \tilde G_{k,\ell} $ the number of general networks  with $ k $ reticulation vertices in the vertex-labeled (leaf-labeled) case.  We 
  focus mainly on  proves the following results.
\bt\label{thm-G1}
For the number $G_{k,n}$ of vertex-labeled phylogenetic networks with $k\geq 1$ reticulation vertices, there is a positive constant
$d_k$ such that
\[
G_{k,n}\sim d_k\left(1-(-1)^n\right)\left(\frac{\sqrt{2}}{e}\right)^nn^{n+2k-1},\qquad (n\rightarrow\infty).
\]

In particular,

\begin{align*}
d_1=\frac{\sqrt{2}}{4};\qquad
d_2=\frac{\sqrt{2}}{32};\qquad
d_3=\frac{\sqrt{2}}{384}.
\end{align*}
\et
The result reveals that the first and second order asymptotics are the  same as  the once for vertex-labeled tree-child networks (see \cite{unpublishedkeyA}). In other words, we can show that for the general networks with fixed number of reticulation vertices, the additional networks that not satisfying the  tree-child conditions are asymptotically negligible as $\nti$. Also,  this approach help us to have following result.
\begin{theo}\label{thm-G2}
	For the numbers $\tilde{G}_{k,\ell}$  of leaf-labeled  general networks  with $k\geq 1$ reticulation vertices, we have
	\[
	\tilde{G}_{k,\ell}\sim2^{3k-1}d_k\left(\frac{2}{e}\right)^{\ell}\ell^{\ell+2k-1},\qquad(\ell\rightarrow\infty)
	\]
	where $d_k$ is as in Theorem \ref{thm-G1}.
\end{theo}

\brem
Note that this result only  holds for fixed $ k $ as $ n $ goes to infinity. The case when $ k $ approaches to infinity cannot be done in this way.  
\erem
The other objective of this paper is to study procedures that can be used to extract explicit formulas for the number of phylogenetic networks with fixed reticulation vertices.  So points of the presented argumentation can shed some light on open questions that are left in \cite{unknown} and \cite{tree-child}  about  an explicit formula for the count of phylogenetic networks.
\begin{table}[ht]
	\centering
		\begin{tabular}{c||c|c}
		  &  &  \\
		    &  The explicit formula  &  \\
		  &  &  \\
			\hline
			$ {G}_{1,\ell}^{\nshortparallel} $ & 
		$ \ell!\;  \Big(r_1(\ell)\displaystyle2^{-\ell}\binom{2\ell}{\ell}-2^{\ell}p_1(\ell)\Big) $
			& $r_1(\ell)= \ell, $\quad $ \text{and }$  $ p_1(\ell)=\frac{1}{2}.$\quad\quad\quad\quad\quad\quad\quad\quad\quad\quad\quad\quad\quad\quad\\
			\hline
			 ${G}_{2,\ell}^{\nshortparallel} $ & $ \ell!\;  \Big(r_2(\ell)\displaystyle2^{-\ell}\binom{2\ell+2}{\ell+1}-2^{\ell}p_2(\ell)\Big) $ & $r_2(\ell)= \frac{(\ell+1)(6\ell^4+19\ell^3+18\ell^2-7\ell)}{2(6\ell-3)(2\ell+1)},$ \quad $ \text{and }  $ $p_2(\ell)=\frac{2\ell^2+5\ell+3}{2}.$\\
			\hline
			 &  &$ r_3(\ell)=\quad\quad\quad\quad\quad\quad\quad\quad\quad\quad\quad\quad\quad\quad\quad\quad\quad\quad\quad\quad\quad$ \\
			$ {G}_{3,\ell}^{\nshortparallel} $ & $ \ell!\;  \Big(r_3(\ell)\displaystyle2^{-\ell}\binom{2\ell+4}{\ell+2}-2^{\ell}p_3(\ell)\Big) $&$ \frac{(\ell+1)(\ell+2)(280\ell^6+3072\ell^5+12834\ell^4+22386\ell^3+10949\ell^2-5211\ell-3990) }{840(2\ell+3)(2\ell+1)(2\ell-1)}. $    \\
			 &  & $p_3(\ell)=\frac{48\ell^4+415\ell^3+1326\ell^2+1799\ell+816}{768}.\quad\quad\quad\quad\quad\quad\quad\quad\quad$ \\
		\end{tabular}
		\smallskip
	\caption{The numbers of leaf-labeled general networks with $ \ell $ leaves and no multiple edges. \label{TablNumbers}}
\end{table}

\section{Generating functions and methods from Analytic combinatorics }
\label{sec:GF}
This section summarizes some  basic concepts on combinatorial classes and their generating functions that will be used in our work. 
Our presentation follows closely \cite{AnaCombi} (although with much less details), 
and the reader interested to know more on the topic is refered to \cite[mainly Chapters I.5, II.1, II.5, VI.3, VII.3, VII.4]{AnaCombi}.
\subsection{(Univariate) generating functions and counting}
Generally speaking, a \emph{combinatorial class}  is a collection $ \mathcal{C} $ of objects of a similar kind (e.g.
words, trees, graphs), endowed with a suitable notion of size or weight (which is a function $ f :  \mathcal{C}  \longrightarrow  \mathbb{N}  $) in a way that there are only finitely many objects of each size. We denote by $\mathcal{C}_n$ the set of objects of size $n$ in $\mathcal{C}$, and by $c_n$ the cardinality of $\mathcal{C}_n$. Specifically in this paper, 
we consider a family of general phylogenetic networks as a combinatorial class such that the size of a network is its number of vertices or leave.

 The arrangement of $ n $ atoms can be considered as objects of size $ n $  in $\mathcal{C}$ such that  each atom has size $ 1 $. 
 In our context, these atoms are the vertices (or leaves) of the networks.
In general, combinatorial objects may or may not be labelled, depending on whether the atoms constituting an object are distinguishable from one another (\emph{labelled} case) or not (\emph{unlabelled} case). Here, our networks will be labelled combinatorial objects.

To deal with a labelled combinatorial class $\mathcal{C}$, we introduce \emph{exponential}  \emph{generating function} $C(z) = \sum_{n\geq 0} c_n \frac{z^n}{n!}$, which is a formal power series in $z$ represents the entire counting sequence of $\mathcal{C}$. The neutral class $ \mathcal{E} $ is made of a single object of size $ 0 $, and its associated generating function is $ E(z) = 1 $. The atomic class $ \mathcal{Z} $ is made of a single
object of size $ 1 $, and its associated generating function is $ Z(z) = z $.

We now turn our attention to recursive \textit{specifications} of a combinatorial class. 
For instance, trees are best described recursively. Note that in the next  section we are going to  describe  decomposition of phylogenetic network that is based on tree structure which will then be translated into a functional equation involving their associated  exponential generating functions.
\begin{ex}A rooted plane tree is a root to which is attached
a (possibly empty) sequence of trees. In other words, the class $ \mathcal{T} $ of rooted plane trees is definable by the recursive equation:
 $ \mathcal{T} =\mathcal{Z}+\mathcal{Z} \times\mathcal{T} +\mathcal{Z} \times\mathcal{T}^2+ +\mathcal{Z} \times\mathcal{T}^3+\cdots =\mathcal{Z} \times SEQ(\mathcal{T}) $. 

Specifications describing a combinatorial class  that is iterative can be represented as a single term built on
$ \mathcal{E} $, $ \mathcal{Z} $ and the constructions $ +, ×, SEQ , CYC , MSET , PSET $. For instance, Cartesian products (with consistent relabellings in the case of labelled objects) correspond to products of series, and sequences (\emph{i.e.}, $m$-tuples of objects of a class $\mathcal{C}$, for any $m \geq 0$) to the quasi-inverse $\frac{1}{1-C(z)}$. This holds for exponential generating functions of labelled objects.
For labelled classes, the precise statement that we refer to is \cite[Theorem II.1]{AnaCombi}. 
Here we get that the corresponding generating function satisfies $T(z) = \frac{z}{1-T(z)}$.
The next step is to have access to the enumeration sequence $(t_n)$ of a class $\mathcal{T}$ from an equation satisfied by the generating function $T(z)$ of $\mathcal{T}$. 
To state it, we introduce the notation $n!\cdot [z^n] T(z)$ to denote the $n$-th coefficient of the series $T(z)$; 
that is to say, writing $T(z) = \sum_{n\geq 0} t_n \frac{z^n}{n!}$, we have $ t_n:= n!\cdot[z^n] T(z) $. From this point on, basic algebra  does the rest. First the original equation is equivalent to $ T-T^2-z=0 $. Solving this  quadratic equation gives
\begin{align*} 
T(z)&=
\frac{1}{2}(1-\sqrt{(1-4z)}) \\
&=z+ z^2 + 2 z^3 + 5 z^4 + 14 z^5 + 42 z^6 + 132 z^7 + 429 z^8+\cdots\\
&=\sum_{n\geq 1} \frac{1}{n}\binom{2n-2}{n-1}z^n, \\
\end{align*}
and  consequently, 
 $ t_n= \displaystyle n! \cdot \frac{1}{n}\binom{2n-2}{n-1}  $.
\end{ex}
The other possible way, especially in the case of tree-like objects, is to appeal to the \emph{transfer theorem}(see \cite{AnaCombi}, VI.1).  Before going ahead, first we illustrate   some concepts which help us to clarify the details.
A \textit{singularity} of an analytic function $ f (z) $ is
a point $ z_0 $ on the boundary of its region of analycity for which $ f(z) $ is not
analytically continuable. Singularities of a function analytic at $ 0 $, which lie
on the boundary of the disc of convergence, are called dominant singularities.
In this case, a dominant singularity is a singularity with smallest modulus. From
\textit{Pringsheim}'s theorem (\cite{AnaCombi}, Theorem IV.6) we know that  if $ f(z) $ is representable at the origin by a
series expansion that has non-negative coefficients and radius of convergence $ \rho $, then
the point $ z=\rho $ is a singularity of $f(z)$. The idea  behind the \textit{transfer  theorem} is that if $ A(z) $ and $ B(z) $ are two generating functions with the same positive real number $ \rho $ as dominant singularity; So when $ z{\rightarrow} \rho $,  we can write $ A(z){ \rightarrow} B(z) $. We obtain the asymptotic expansion of $ [z^n] A(z) $ by transferring the behaviour of $ A(z) $ around its dominant singularity from a simpler
function $ B(z) $, from which we know the analytic behaviour.
\bt[Transfer  Theorem]
If the generating function $A(z)$ admits an expansion of the form $A(z) \sim c\cdot (1-\frac{z}{\rho})^{-\alpha} $  as $ {\nti } $, around its (unique) dominant singularity $ \rho $,  then we have
$$ 
n! \cdot [z^n] A(z) \sim n! \cdot c\cdot \dfrac{n^{\alpha-1}}{\Gamma (\alpha)}\cdot
\rho^{-n},$$
as $ {\nti } $.
\et
\brem
Here  $ A(z)$ is analytic in the disk of radius $ \rho $ centered at the origin.
\erem
 Recall that $[z^n] A(z)$ is the coefficient of $z^n$ in $A(z)$, and so it is $\frac{c_n}{n!}$  (resp. $c_n$ ), when $A(z)$ is a exponential (resp.ordinary) generating function.
Note that the location of a dominant singularity will give the exponential growth
of the sequence, and the nature of this singularity the subexponential term.
 If $ A(z) $ has several dominant singularities coming from\textit{ pure periodicities} (for more details see \cite{AnaCombi}, IV.6.1  ), then the contributions from each of them must be combined.

These methods are fundamental results from complex analysis that  allow to set up generating function in its disk of convergence, but not always. In particular, the \emph{transfer theorem}   (Theorem VI.1 of~\cite{AnaCombi})  is one of the  suitable tools, 
which allows us to derive asymptotic estimates of the coefficients of generating functions.

\subsection{Additive parameters and bivariate generating function}

It is sometimes interesting to analyze the behaviour of other parameters than size. For example, interesting parameters for plane trees can be: height, number of leaves, path length, etc. These parameters are important for algorithm analysis as they correspond to the performance of algorithms that compute with or are modeled by plane trees. We now consider \emph{multivariate} generating functions, where additional variables ($x$, $y$, \dots) record the value of other parameters of our objects. One variable  is used to track the size of the structure (e.g. number of nodes in a plane tree) and the other is used to track the parameter of interest (e.g. height, number of leaves, path length). 
 
In our cases, we will consider one more such parameter, which are numbers of certain "unary nodes" occuring in our objects. 
Namely, denoting $c_{n,\ell}$ the number of objects of size $n$ in the combinatorial class $\mathcal{C}$ 
such that the parameter has value $\ell$, the multivariate exponential generating function we consider is $C(z,y) = \sum_{n,\ell} c_{n,\ell} y^{\ell} \tfrac{z^n}{n!}  $. 

For inistanc,on the previous example of rooted plane trees consider one additional parameter, 
which is the number of leaves nodes.
The coefficient of $z^ny^\ell$ in the generating function $T(z,y)$ is then the number of rooted plane trees with $n$ nodes and exactly $\ell$ leaves, divided by $n!$. 

The \textquotedblleft dictionnary\textquotedblright$\;$translating combinatorial specifications to equations satisfied by the generating function extends to multivariate series, and our  specification that shows any such tree is leaf or sequences $ (\geq 1) $ of trees that attached to  the root nod. This
gives $T(z,y) = zy + \dfrac{z T(z,y)}{1-T(z,y)}$. 
\section{Decomposing binary phylogenetic networks}

In order to count general phylogenetic networks, we will adjust the procedure of sparsened skeleton decomposition for general networks. This method  is well studied   for  tree-child networks, with $ k $ of reticulation vertices in \cite{unpublishedkeyA}.  We will  use the
decomposition to obtain a reduction which can be easily analyzed by means of generating functions.
Consider a general phylogenetic network having $k$
reticulation vertices. Then each such vertex has two incoming edges. If one edge is removed for each of the $ k $ reticulation vertices, then  
the remaining graph is again Motzkin tree (labeled and
nonplane). Depending on our choice for removed edges, this Motzkin tree has at most  $ 2k $ unary nodes.
Recall that for tree-child networks this method gives exactly $ 2k $ unary nodes.

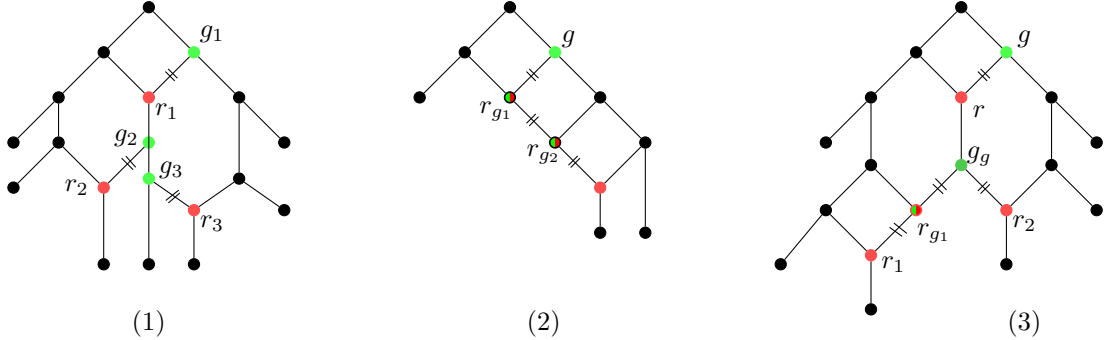
\begin{figure}[h]
	\begin{center}
		\begin{tikzpicture}[scale=0.6]

		\draw (-9cm,0cm) node[inner sep=1.5pt,circle,draw,fill] (1) {};
		\draw (-10cm,-1cm) node[inner sep=1.5pt,circle,draw,fill] (2) {};
		\draw (-8cm,-1cm) node[color=green!70,inner sep=1.5pt,circle,draw,fill] (3) {};
		\draw (-11cm,-2cm) node[inner sep=1.5pt,circle,draw,fill] (4) {};
		\draw (-9cm,-2cm) node[color=red!70,inner sep=1.5pt,circle,draw,fill] (5) {};
		\draw (-7cm,-2cm) node[inner sep=1.5pt,circle,draw,fill] (6) {};
		\draw (-12cm,-3cm) node[inner sep=1.5pt,circle,draw,fill] (7) {};
		\draw (-6cm,-3cm) node[inner sep=1.5pt,circle,draw,fill] (8) {};
		\draw (-9cm,-3.8cm)node[color=green!70,inner sep=1.5pt,circle,draw,fill] (10) {};
		\draw (-7cm,-3.8cm) node[inner sep=1.5pt,circle,draw,fill] (11) {};
		\draw (-12cm,-4cm) node[inner sep=1.5pt,circle,draw,fill] (12) {};
		\draw (-10cm,-4cm) node[color=red!70,inner sep=1.5pt,circle,draw,fill] (13) {};
		\draw (-8cm,-4.5cm) node[color=red!70,inner sep=1.5pt,circle,draw,fill] (14) {};
		\draw (-6cm,-4.5cm) node[inner sep=1.5pt,circle,draw,fill] (15) {};
		\draw (-10cm,-5.7cm) node[inner sep=1.5pt,circle,draw,fill] (16) {};
		\draw (-8cm,-5.7cm) node[inner sep=1.5pt,circle,draw,fill] (17) {};
		\draw (-11cm,-3cm) node[inner sep=1.5pt,circle,draw,fill] (101) {};
		\draw (-9cm,-3cm)node[color=green!70,inner sep=1.5pt,circle,draw,fill] (102) {};
		\draw (-9cm,-5.7cm) node[inner sep=1.5pt,circle,draw,fill] (103) {};
		\draw (-8.6cm,-2.3cm) node[inner sep=1.5pt,circle] (18) {$r_1$};
		\draw (-10.6cm,-4cm) node[inner sep=1.5pt,circle] (19) {$r_2$};
		\draw (-7.6cm,-4.8cm) node[inner sep=1.5pt,circle] (20) {$r_3$};
		\draw (-7.6cm,-.6cm) node[inner sep=1.5pt,circle] (21) {$g_1$};
		\draw (-9.5cm,-2.9cm) node[inner sep=1.5pt,circle] (22) {$g_2$};
		\draw (-8.5cm,-3.6cm) node[inner sep=1.5pt,circle] (23) {$g_3$};
		\draw (-8.73cm,-1.33cm) node[inner sep=1.5pt,circle] (21) {};
		\draw (-8.33cm,-1.73cm) node[inner sep=1.5pt,circle] (22) {};
		\draw (-8.67cm,-1.27cm) node[inner sep=1.5pt,circle] (23) {};
		\draw (-8.27cm,-1.67cm) node[inner sep=1.5pt,circle] (24) {};
		\draw (-9.70cm,-3.25cm) node[inner sep=1.5pt,circle] (25) {};
		\draw (-9.25cm,-3.70cm) node[inner sep=1.5pt,circle] (26) {};
		\draw (-9.60cm,-3.20cm) node[inner sep=1.5pt,circle] (27) {};
		\draw (-9.20cm,-3.65cm) node[inner sep=1.5pt,circle] (28) {};
		\draw (-8.23cm,-3.97cm) node[inner sep=1.5pt,circle] (29) {};
		\draw (-8.73cm,-4.37cm) node[inner sep=1.5pt,circle] (30) {};
		\draw (-8.21cm,-4.03cm) node[inner sep=1.5pt,circle] (31) {};
		\draw (-8.67cm,-4.43cm) node[inner sep=1.5pt,circle] (32) {};
		\draw (-9cm,-7cm) node[inner sep=1.5pt,circle] (44) {$(1)$};

		\draw (1)--(2);\draw (1)--(3);\draw (2)--(4);\draw (2)--(5);\draw (103)--(10);
		\draw (3)--(5);\draw (3)--(6);\draw (4)--(7);\draw (6)--(8);
		\draw (4)--(101);\draw (5)--(10);\draw (6)--(11);\draw (101)--(12);
		\draw (101)--(13);\draw (102)--(13);\draw (10)--(14);\draw (11)--(14);
		\draw (11)--(15);\draw (13)--(16);\draw (14)--(17);\draw (21)--(22);
		\draw (23)--(24);\draw (25)--(26);\draw (27)--(28);\draw (29)--(30);
		\draw (31)--(32);

		\draw (-1cm,0cm) node[inner sep=1.5pt,circle,draw,fill] (1) {};
		\draw (-2cm,-1cm) node[inner sep=1.5pt,circle,draw,fill] (2) {};
		\draw (0cm,-1cm) node[color=green!70,inner sep=1.5pt,circle,draw,fill] (3) {};
		\draw (-3cm,-2cm) node[inner sep=1.5pt,circle,draw,fill] (4) {};
		\draw (-1cm,-2cm) node[inner sep=1.5pt,circle,draw,fill] (5) {};
		
		\coordinate (center) at (-1,-2);
		
		\fill[green] (center) + (0, .1) arc (90:270:.1);
		
		\fill[red] (center) + (0, -.1) arc (270:450:.1);
		
		\draw (1cm,-2cm) node[inner sep=1.5pt,circle,draw,fill] (6) {};
		\draw (0cm,-3cm) node[inner sep=1.5pt,circle,draw,fill] (7) {};
		\coordinate (center) at (0,-3);
		\fill[green] (center) + (0, .1) arc (90:270:.1);
		\fill[red] (center) + (0, -.1) arc (270:450:.1);
		
		\draw (2cm,-3cm) node[inner sep=1.5pt,circle,draw,fill] (8) {};
		\draw (1cm,-4cm) node[color=red!70,inner sep=1.5pt,circle,draw,fill] (9) {};
		\draw (1cm,-5cm) node[inner sep=1.5pt,circle,draw,fill] (20) {};
		\draw (2cm,-5cm) node[inner sep=1.5pt,circle,draw,fill] (21) {};
		\draw (-1.3cm,-2.3cm) node[inner sep=1.5pt,circle] (10){$r_{g_1}$};
		\draw (-.25cm,-3.25cm) node[inner sep=1.5pt,circle] (11) {$r_{g_2}$};
		\draw (.3cm,-.7cm) node[inner sep=1.5pt,circle] (30){$g$};
		\draw (-.73cm,-1.33cm) node[inner sep=1.5pt,circle] (12) {};
		\draw (-.33cm,-1.73cm) node[inner sep=1.5pt,circle] (13) {};
		\draw (-.67cm,-1.27cm) node[inner sep=1.5pt,circle] (14) {};
		\draw (-.27cm,-1.67cm) node[inner sep=1.5pt,circle] (15) {};
		\draw (-.33cm,-2.27cm) node[inner sep=1.5pt,circle] (16) {};
		\draw (-.73cm,-2.67cm) node[inner sep=1.5pt,circle] (17) {};
		\draw (-.27cm,-2.33cm) node[inner sep=1.5pt,circle] (18) {};
		\draw (-.67cm,-2.73cm) node[inner sep=1.5pt,circle] (19) {};
        \draw (.20cm,-3.60cm) node[inner sep=1.5pt,circle] (22) {};
        \draw (.60cm,-3.20cm) node[inner sep=1.5pt,circle] (23) {};
        
        \draw (.25cm,-3.65cm) node[inner sep=1.5pt,circle] (24) {};
        \draw (.65cm,-3.25cm) node[inner sep=1.5pt,circle] (25) {};
        \draw (-.25cm,-7cm) node[inner sep=1.5pt,circle] (55) {$(2)$};
		\draw (1)--(2);\draw (1)--(3);\draw (2)--(4);\draw (2)--(5);
		\draw (3)--(5);\draw (3)--(6);\draw (5)--(7);\draw (6)--(8);
		\draw (7)--(9);\draw (6)--(7);\draw (12)--(13);\draw (14)--(15);
		\draw (16)--(17);\draw (18)--(19);\draw (8)--(9);\draw (9)--(20);\draw (8)--(21);\draw (22)--(23);\draw (24)--(25);
		\definecolor{red}{rgb}{2,0,0}
		\definecolor{OliveGreen}{rgb}{0,0.7,0}

		\draw (9cm,0cm) node[inner sep=1.5pt,circle,draw,fill] (1) {};
		\draw (8cm,-1cm) node[inner sep=1.5pt,circle,draw,fill] (2) {};
		\draw (10cm,-1cm) node[color=green!70,inner sep=1.5pt,circle,draw,fill] (3) {};
		\draw (7cm,-2cm) node[inner sep=1.5pt,circle,draw,fill] (4) {};
		\draw (9cm,-2cm)  node[color=red!70,inner sep=1.5pt,circle,draw,fill] (5) {};
		\draw (11cm,-2cm) node[inner sep=1.5pt,circle,draw,fill] (6) {};
		\draw (6cm,-3cm) node[inner sep=1.5pt,circle,draw,fill] (7) {};
		\draw (12cm,-3cm) node[inner sep=1.5pt,circle,draw,fill] (8) {};
		\draw (7cm,-3.5cm) node[inner sep=1.5pt,circle,draw,fill] (9) {};
		\draw (9cm,-3.5cm)  node[color=OliveGreen!70,inner sep=1.5pt,circle,draw,fill] (10) {};
		\draw (11cm,-3.5cm) node[inner sep=1.5pt,circle,draw,fill] (11) {};
		\draw (6cm,-4.5cm) node[inner sep=1.5pt,circle,draw,fill] (12) {};
		\draw (8cm,-4.5cm)  node[color=red!70,inner sep=1.5pt,circle,draw,fill] (13) {};
		\coordinate (center) at (8cm,-4.5cm);
		\fill[green] (center) + (0, .1) arc (90:270:.1);
		\fill[red] (center) + (0, -.1) arc (270:450:.1);

		\draw (10cm,-4.5cm)  node[color=red!70,inner sep=1.5pt,circle,draw,fill] (14) {};
		\draw (12cm,-4.5cm) node[inner sep=1.5pt,circle,draw,fill] (15) {};
		\draw (7cm,-5.5cm) node[color=red!70,inner sep=1.5pt,circle,draw,fill] (16) {};
		\draw (10cm,-5.7cm) node[inner sep=1.5pt,circle,draw,fill] (17) {};
		\draw (7cm,-6.7cm) node[inner sep=1.5pt,circle,draw,fill] (33) {};
		\draw (5cm,-5.7cm) node[inner sep=1.5pt,circle,draw,fill] (34) {};

		\draw (9.4cm,-2.3cm) node[inner sep=1.5pt,circle] (18) {$r$};
		\draw (10.4cm,-.7cm) node[inner sep=1.5pt,circle] (48) {$g$};
		\draw (9.4cm,-3.3cm) node[inner sep=1.5pt,circle] (58) {$g_g$};
		\draw (8.4cm,-5cm) node[inner sep=1.5pt,circle] (19) {$r_{g_1}$};
		\draw (10.4cm,-4.8cm) node[inner sep=1.5pt,circle] (20) {$r_2$};
		\draw (7.5cm,-5.7cm) node[inner sep=1.5pt,circle] (200) {$r_1$};
		\draw (9.27cm,-1.33cm) node[inner sep=1.5pt,circle] (21) {};
		\draw (9.67cm,-1.73cm) node[inner sep=1.5pt,circle] (22) {};
		\draw (9.33cm,-1.27cm) node[inner sep=1.5pt,circle] (23) {};
		\draw (9.73cm,-1.67cm) node[inner sep=1.5pt,circle] (24) {};
		\draw (8.27cm,-3.83cm) node[inner sep=1.5pt,circle] (25) {};
		\draw (8.70cm,-4.23cm) node[inner sep=1.5pt,circle] (26) {};
		\draw (8.33cm,-3.70cm) node[inner sep=1.5pt,circle] (27) {};
		\draw (8.77cm,-4.15cm) node[inner sep=1.5pt,circle] (28) {};
		\draw (9.67cm,-3.77cm) node[inner sep=1.5pt,circle] (29) {};
		\draw (9.27cm,-4.17cm) node[inner sep=1.5pt,circle] (30) {};
		\draw (9.73cm,-3.83cm) node[inner sep=1.5pt,circle] (31) {};
		\draw (9.33cm,-4.23cm) node[inner sep=1.5pt,circle] (32) {};
			\draw (7.80cm,-5.16cm) node[inner sep=1.5pt,circle] (291) {};
			\draw (7.30cm,-4.80cm) node[inner sep=1.5pt,circle] (301) {};
			\draw (7.4cm,-4.7cm) node[inner sep=1.5pt,circle] (311) {};
			\draw (7.9cm,-5.1cm) node[inner sep=1.5pt,circle] (321) {};
			\draw (291)--(301);\draw (311)--(321);
		\draw (1)--(2);\draw (1)--(3);\draw (2)--(4);\draw (2)--(5);
		\draw (3)--(5);\draw (3)--(6);\draw (4)--(7);\draw (6)--(8);
		\draw (4)--(9);\draw (5)--(10);\draw (6)--(11);\draw (9)--(12);
		\draw (9)--(13);\draw (10)--(13);\draw (10)--(14);\draw (11)--(14);\draw (16)--(33);\draw (16)--(12);\draw (12)--(34);
		\draw (11)--(15);\draw (13)--(16);\draw (14)--(17);\draw (21)--(22);
		\draw (23)--(24);\draw (25)--(26);\draw (27)--(28);\draw (29)--(30);
		\draw (31)--(32);
	\draw (10.4cm,-7cm) node[inner sep=1.5pt,circle] (20) {$(3)$};
		\end{tikzpicture}	
	\end{center}\caption{Three  general phylogenetic networks, where  $ (1) $ after adding marked edges there is bijection between green and red (reticulation) vertices in the Motzkin skeleton. General networks depicted in  $ (2)$  and $(3) $ which the red-green and double-green  vertices appear after adding the marked edges.  Edges are directed downwards.} \label{GD}
\end{figure}

Now consider the  following procedure: start with a Motzkin tree $T$ with not more that  $2k$ unary nodes and $n$ vertices in total.
\begin{itemize}
	\item Add directed edges such that each edge connects two unary nodes and any two edges do not have a vertex in common. Color the started vertices of  the added directed edges green and their  end vertices red. Note that if Motzkin tree has exactly $ 2k $ unary nodes, the coloring procedure imposes that  there will be equal $ k $ colored green  and red  nodes (see Figure~\ref{GD}, (1)).
	\item Consider two unary vertices and joint them by using  sequence $\mathcal{S}$ ($ 2\leq   \mid\mathcal{S}\mid  \leq k  $) of fixed number of directed edges  in the following way.
	One of edge in the $\mathcal{S}$ connects the first unary vertex of the Motzin tree to a leaf which we call $r_{g}$. Then  connect $r_{g}$ to another leaf and continue this process for disjoint leaves until use all directed edges but one. Now connect the last leaf to a second unary vertex by the remaining edge. As similar before color first unary vertex green and consider red color for second ones, then mark (color) all leaves on the path (leaf) \textit{red-green} ( Figure~\ref{GD}, (2)). Note that for a general network with $ k $ reticulation vertices, the number of directed edges in the considered sequence, cannot be exceeded of $k$, because each marked red-green  vertex is reticulation vertex. 	
	\item  Consider a leaf  $ g_g $  of the Motzkin tree.   As similar before connect $ g_g  $ to the two distinct unary vertices  by using two sequences $ \mathcal{S} \text{and}  \mathcal{S^\prime} $ of outgoing directed  edges ($ 2\leq   \mid\mathcal{S}\cup\mathcal{S^\prime} \mid  \leq k  $). Mark $ g_g $ as  \textit{double-green} vertex and then  color   targeted unary  vertices red. Also, consider  red-green color for  all the leaves on the paths of $ g_g $ to the unary  vertices.  ( see Figure~\ref{GD} (3)).  
\end{itemize}
Note that in the above procedure the resulting graph must be
a general phylogenetic network $ \mathcal{G} $. We say then that $ T $ (keeping the colors from the above generation
of $ \mathcal{G} $, but not the edges) is a \textit{colored Motzkin skeleton} (or simply Motzkin skeleton) of $ \mathcal{G} $. 
Now, consider two sets, but not necessarily disjoint,  of  colored vertices obtained of above procedure. The member of first set is all colored  vertices with  outgoing edges and then assume all colored vertices with ingoing edges in the second set. Call them \textit{pointer} and \textit{target} sets respectively. In this way,  red-green vertices are considered in both pointer and target sets. It is not hard to see that  the size of  target set  is correspondent with number of  reticulation nodes on a general phylogenetic network.
Note that in this procedure  any general network with no multiple edges  and $ n $ vertices is generated and each of them exactly $ 2^k $ times, so in this case
every network $ \mathcal{G} $ with $  k $ reticulation vertices has exactly $ 2^k $ different Motzkin skeletons. 
However, note that  as opposed to defined subclasses of phylogenetic networks like  tree-child  networks, here we assume   multiple edges (reticulation vertex with one parent) are allowed to be in general networks. 
So for a reticulation vertex  with just one parent, any arbitrary choice  and removing of multiple edges, causes the same Motzkin skeleton. It means the described procedure generates  a network with $ k $ reticulation vertices and $ r $ multiple edges exactly $ 2^{k-r} $ times. 
  In the first step our aim is  to set up an exponential generating functions for general networks with no multiple edges and then get  the correspondent exponential generating function for other networks with at least one multiple edge.

\begin{figure}[h]
	\begin{center}
		\begin{tikzpicture}[scale=0.6]
		\draw (4cm,0cm) node[inner sep=1.5pt,circle,draw,fill] (1) {};
		\draw (3cm,-1cm) node[color=green!70,inner sep=1.5pt,circle,draw,fill] (2) {};
		\draw (5cm,-1cm) node[color=green!70,inner sep=1.5pt,circle,draw,fill] (3) {};
		\draw (2cm,-2cm)node[color=green!70,inner sep=1.5pt,circle,draw,fill] (4) {};
		
		\draw (1)--(2);\draw (1)--(3);\draw (2)--(4);

		\draw (10cm,0cm) node[inner sep=1.5pt,circle,draw,fill] (1) {};
		\draw (9cm,-1cm) node[inner sep=1.5pt,circle,draw,fill] (2) {};
		\coordinate (center) at (9,-1);
		\fill[green] (center) + (0, .1) arc (90:270:.1);
		\fill[red] (center) + (0, -.1) arc (270:450:.1);
		\draw (11cm,-1cm) node[color=green!70,inner sep=1.5pt,circle,draw,fill] (3) {};
		\draw (12cm,-2cm) node[inner sep=1.5pt,circle,draw,fill] (4) {};
		\coordinate (center) at (12,-2);
		\fill[green] (center) + (0, .1) arc (90:270:.1);
		\fill[red] (center) + (0, -.1) arc (270:450:.1);
		
		\draw (1)--(2);\draw (1)--(3);\draw (4)--(3);

		\definecolor{OliveGreen}{rgb}{0,0.7,0}
		\draw (19cm,0cm) node[inner sep=1.5pt,circle,draw,fill] (30) {};
		\draw (18cm,-1cm) node[inner sep=1.5pt,circle,draw,fill] (31) {};
		\draw (19cm,-2cm) node[color=OliveGreen!70,inner sep=1.5pt,circle,draw,fill] (32) {};
		\draw (17cm,-2cm) node[inner sep=1.5pt,circle,draw,fill] (33) {};
		\draw (20cm,-1cm) node[color=green!70,inner sep=1.5pt,circle,draw,fill] (34) {};
		\coordinate (center) at (17,-2);
		\fill[green] (center) + (0, .1) arc (90:270:.1);
		\fill[red] (center) + (0, -.1) arc (270:450:.1);
		\draw (30)--(31);\draw (30)--(34);\draw (31)--(33);\draw (31)--(32);
		\end{tikzpicture}
	\end{center}
	\caption{Corresponding sparsened skeletons of Figure \ref{GD} networks.}
	\label{sparsened}
\end{figure}
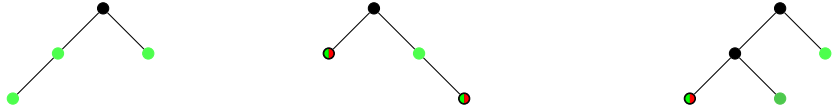

In order to set up exponential generating functions for the number of  general phylogenetic networks, we will construct them as
follows: For a given network $\mathcal{G}$ 
fix one of its possible Motzkin tree skeletons, that shows us how the pointer set vertices are distributed within $\mathcal{G}$ (for instance consider networks in Figure \ref{GD} without marked edges ). Now look for  sparsened skeleton of $\mathcal{G}$ which contains all pointer set vertices and contract all paths between any two vertices
which are either pointer  vertices or an ancestor of  them to one edge.  Note that this ancestor may
be pointer vertices itself (also see, \cite{unpublishedkeyA}).
In order to construct general networks with $k$ target  vertices (reticulations), we consider a sparsened skeleton
having no more than $ k $  pointer  vertices. Then we replace all edges by paths that are made of red vertices or binary vertices with a Motzkin tree (whose unary vertices are all colored red) as second child and add a
path of the same type on top of the root of the sparsened skeleton. Moreover, we attach a Motzkin
tree (again with all unary vertices colored red) only  to those leaves of the sparsened skeleton such that are just colored green (not red-green or double-green). Note that red-green and double-green  nodes lie on  leaves of sparsened skeleton. Do all of the above in such a way that the new structure has $ k $ target vertices (red and red-green) altogether. What we obtain so far is a Motzkin skeleton of a phylogenetic network. Finally, add
edges connecting the pointer vertices to the target ones in such a way that the general  phylogenetic networks  condition  is respected.
As an  advantage, a similar  procedure can be used to set up generating functions for other kinds of phylogenetic networks, with fixed number of reticulation vertices, such as  ``stack-free''  and ``galled'' networks that are  defined in \cite{Semple2018, huson2010phylogenetic}. \\
Let us set up the exponential generating function for the Motzkin trees which appear in the above
construction. After all, the unary  vertices in those trees will be
the red nodes of our network. 

Denote by $M_{\ell,n}$ the number of all vertex-labeled Motzkin trees $n$ vertices and $\ell$
unary vertices. Furthermore, let $\mathcal M$
be the set of all these Motzkin trees. The exponential generating function associated to $\mathcal
M$ is 
\[
M(z,y)=\sum_{n\ge 1}\sum_{\ell\ge 0} M_{\ell,n} y^\ell \frac{z^n}{n!}.
\]
Furthermore, let $M(z,y)$ denotes the generating function associated to all Motzkin trees
in $\mathcal M$ whose root is a unary node or a binary node, so we have
\[
M(z,y)=z+zyM(z,y)+\frac{z}{2}M^2(z,y).
\] 
and thus 
\begin{align} \label{fctM}
M(z,y)=\dfrac{1-zy-\sqrt{1+(y^2-2)z^2-2zy}}{z}. 
\end{align} 
The first few coefficients can be seen from 
\[
z+yz^2+(y^2+\frac{1}{2})z^3+(y^3+\frac{3}{2}y)z^4+\cdots.
\]
\section{Counting Vertex And Leaf-Labeled General Phylogenetic Network}
In this section,  Our main goal  is to present a precise asymptotic result for the number of  general phylogenetic networks with a fixed number  $  k  $ of reticulation vertices. To clear up the methods  we  start with simple cases, determine the asymptotic number of general phylogenetic networks with up to $ 3 $ reticulation nodes. After that, we will show how this approach helps us to present  explicit formulas for the exact number of vetrex and leaf-labeled of them.
Finally, we will focus on the general case and show how previous results  lead us to prove  theorems \ref{thm-G1} and \ref{thm-G2}, for general phylogenetic networks with $ k $ reticulation vertices. 
As a warm-up consider general phylogenetic network with only one reticulation node we use the procedure to obtain
\eqref{fctM} and the (sparsened) skeleton, as described in the previous section: Consider a general network with no multiple edges,  delete 
one of the two incoming edges of the reticulation node which then gives a unary-binary tree
with exactly two unary nodes which are colored green and red  (we will consider general networks with multiple edges separately).
Conversely, we can start with the general tree or even the sparsened skeleton and then construct the network from this. For more explicitly, Let $ G^{\nshortparallel}_i(z) $ resp. $ G^{\shortparallel}_i(z)  $  denote exponential generating functions for networks with no multiple edges (with multiple edges) and $ i $ reticulation vertices.   
\bp
The exponential generating function for  vertex-labelled   general phylogenetic networks with
one reticulation node is
\begin{align} \label{fctG1}
G_1(z)= G_1^{\nshortparallel}(z)+ G_1^{\shortparallel}(z)
=z\frac{\tilde{a}_1(z^2)-\tilde{b}_1(z^2)\sqrt{1-2z^2}\,}{(1-2z^2)^\frac{3}{2}},
\end{align}
where,
\begin{align} \label{polyG1}
\tilde{a}_1(z)=\tilde{b}_1(z)= 1-z.
\end{align}
\ep
\newcommand{\ket}[1]{$\left|#1\right\rangle$}
\newcommand{\El}[1]{\small $\ell_{#1}$}
\newcommand{\Ke}[1]{$k_{#1}$}
\newcommand{\Rl}[1]{$r_{#1}$}
\begin{figure}[h]
	\begin{center}
		\begin{tikzpicture}[
		scale=0.7,
		level/.style={thick},
		virtual/.style={thick,densely dashed},
		trans/.style={thick,<->,shorten >=2pt,shorten <=2pt,>=stealth},
		classical/.style={thin,double,<->,shorten >=4pt,shorten <=4pt,>=stealth}
		]

		\draw (0cm,0cm) node[inner sep=1.5pt,circle,draw,fill] (1) {};
		\draw (0.7cm,0.7cm) node[inner sep=1.5pt,circle,draw,fill] (2) {};
		\draw (1.4cm,1.4cm) node[inner sep=1.5pt,circle,draw,fill] (3) {};
		\draw (2.4cm,2.4cm) node[inner sep=1.5pt,circle,draw,fill] (4) {};
		\draw (3.1cm,3.1cm) node[inner sep=1.5pt,circle,draw,fill] (5) {};
		\draw (1.6cm,1.6cm) node[inner sep=1.5pt,circle] (11) {};
		\draw (2.2cm,2.2cm) node[inner sep=1.5pt,circle] (12) {};
		\draw (3.6cm,3.3cm) node[inner sep=1.5pt,circle] (14) {$ x $};
		\draw (-0.3cm,0.3cm) node[inner sep=1.5pt,circle] (13) {$g$};
		\draw (0cm,-1.8cm) node[regular polygon,regular polygon sides=3,draw,inner sep=0.1cm] (6) {};
		\draw (.8cm,-1cm) node[regular polygon,regular polygon sides=3,draw,inner sep=0.1cm] (7) {};
		\draw (1.5cm,0cm) node[regular polygon,regular polygon sides=3,draw,inner sep=0.1cm] (8) {};
		\draw (2.4cm,.8cm) node[regular polygon,regular polygon sides=3,draw,inner sep=0.1cm] (9) {};
		\draw (3.2cm,1.7cm) node[regular polygon,regular polygon sides=3,draw,inner sep=0.1cm] (10) {};
		
		\draw (1)--(2);\draw (2)--(3);\draw (4)--(5);\draw (1)--(6);
		\draw (2)--(7);\draw (3)--(8);\draw (4)--(9);\draw (5)--(10);
		\draw[loosely dotted,line width=0.7pt] (11)--(12);
		\draw[trans] (0cm,2em) -- (3.0cm,10em) node[midway,left] {\El{}};
		\draw (2cm,-2cm)  node[inner sep=1.5pt,circle] (1) {(a)};
		
		\draw (10cm,0cm) node[inner sep=1.5pt,circle,draw,fill] (1) {};
		\draw (10cm,-1cm) node[inner sep=1.5pt,circle,draw,fill] (20) {};
		\draw (10.7cm,0.7cm) node[inner sep=1.5pt,circle,draw,fill] (2) {};
		\draw (11.4cm,1.4cm) node[inner sep=1.5pt,circle,draw,fill] (3) {};
		\draw (12.4cm,2.4cm) node[inner sep=1.5pt,circle,draw,fill] (4) {};
		\draw (13.1cm,3.1cm) node[inner sep=1.5pt,circle,draw,fill] (5) {};
		\draw (11.6cm,1.6cm) node[inner sep=1.5pt,circle] (11) {};
		\draw (12.2cm,2.2cm) node[inner sep=1.5pt,circle] (12) {};
		\draw (13.6cm,3.3cm) node[inner sep=1.5pt,circle] (14) {$ x $};
		\draw (9.7cm,0.3cm) node[inner sep=1.5pt,circle] (13) {$g$};
		\draw (10cm,-2cm) node[regular polygon,regular polygon sides=3,draw,inner sep=0.1cm] (6) {};
		\draw (10.8cm,-.7cm) node[regular polygon,regular polygon sides=3,draw,inner sep=0.1cm] (7) {};
		\draw (11.5cm,0.3cm) node[regular polygon,regular polygon sides=3,draw,inner sep=0.1cm] (8) {};
		\draw (12.4cm,1.1cm) node[regular polygon,regular polygon sides=3,draw,inner sep=0.1cm] (9) {};
		\draw (13.2cm,2cm) node[regular polygon,regular polygon sides=3,draw,inner sep=0.1cm] (10) {};
		
		\draw (1)--(2);\draw (2)--(3);\draw (4)--(5);\draw (20)--(6);
		\draw (2)--(7);\draw (3)--(8);\draw (4)--(9);\draw (5)--(10);
		\draw[loosely dotted,line width=0.7pt] (11)--(12);
		\draw[trans] (10cm,2em) -- (13.0cm,10em) node[midway,left] {\El{}};
		\draw  [red,->] (1) to[out=45,in=45] (20);
		\draw  [blue,->] (1) to[out=145,in=145] (20);	
		\draw (12.6cm,-2cm)  node[inner sep=1.5pt,circle] {(b)};
		\end{tikzpicture}
	\end{center}
		\caption{(a) The structure of Motzkin skeletons of networks with one
			reticulation vertex. It originates from a sparsened skeleton which consists of only one green vertex. It has one
			green vertex, denoted by $g$, and one red vertex which is hidden within the forest made of the
			triangles in the picture, which are attached to $g$ and all the vertices on the path of length $\ell$.
			Note that the position of the red vertex in this forest is restricted by the general condition. (b) There is  multiple edge when the root of the subtree which is attached to $ g $ is the red node. } \label{fig_skeleton_1}
\end{figure}
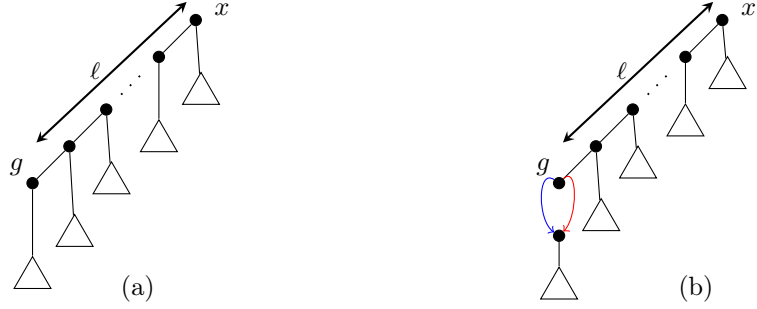
\bpf
We start with the general Motzkin tree as depicted in Figure~\ref{fig_skeleton_1} $ (a) $ and add an edge starting from $g$ and ending at a red vertex. 
Note that for all phylogenetic network, this edge is not allowed
point to a node on the path from $g$ to the root (since the phylogenetic network must be a DAG). Thus, when starting from the sparsened skeleton, \emph{i.e.}, the single green vertex
$g$, then we must add a sequence of trees on top of $g$ which consist of a root (these vertices
make the path from $g$ to the root of the network) to which either a leaf or a binary node with
two trees is attached. The red vertex must be contained in the forest made by this sequence or
the tree attached to $g$. Note that the second expression refers to the depicted structure $ (b) $ which is for general networks with a multiple  edge. In terms of generating functions altogether gives
\[
G_1(z) =\frac{1}{2}\frac{\partial}{\partial y} \frac{z\tilde M(z,y)}{1-z M(z,y)}\Big\vert_{y=0}+ \frac{z^2 M(z,0)}{1-z M(z,0)},
\]
where,
\begin{align}\label{tree-without-root}
\tilde{M}(z,y)&=M(z,y)-zyM(z,y)=(1-zy)M(z,y).
\end{align}
The factor $1/2$  makes up for the fact that  each network in case $ (a) $ is counted by the above procedure exactly twice.  

\epf

From this result we can now easily obtain the asymptotic number of networks. 

\bp
Let $G_{1,n}$ denote the number of vertex-labelled general phylogenetic network with $n$ vertices
and one reticulation vertex. If $n$ is even then $G_{1,n}$ is zero, otherwise 
\[
G_{1,n}=n![z^n]G_1(z)=\(\frac{\sqrt2}{e}\)^n n^{n+1}\(\frac{\sqrt{2}}{2}-\frac{\sqrt{\pi}}{2}\cdot\frac{1}{\sqrt{n}}
+{\mathcal O}\left(\frac{1}{n}\right)\),
\]
as $\nti$.
\ep
\bpf
The function \eqref{fctG1} has two dominant singularities, namely at $\pm1/\sqrt{2}$, with
singular expansions
\[
G_1(z)\stackrel{z\rightarrow 1/\sqrt{2}}{\sim}\frac{1}{8(1-\sqrt{2}z)^{3/2}},\qquad
G_1(z)\stackrel{z\rightarrow -1/\sqrt{2}}{\sim}-\frac{1}{8(1+\sqrt{2}{z})^{3/2}}.
\]
Applying a transfer lemma  
for these two singularities completes the proof.
\epf
\subsubsection{Exact value of vertex-labeled general phylogenetic networks with one reticulation vertex}
First, set $n=2m+1$. Then, from \eqref{fctG1} we obtain 
\[
[z^n]G_1(z)=[z^m]\bar G_{1}(z)
\]
with 
\[
\bar G_{1}(z)= \frac{\tilde{a}_1(z)-\tilde{b}_1(z)\sqrt{1-2z}}{(1-2z)^{3/2}},
\]
where $\tilde{a}_1(z)$ and $\tilde{b}_1(z)$ are as in \eqref{polyG1}.
So we have
\[
[z^m]\bar G_{1}(z)=[z^m]\dfrac{\tilde{a}_1(z)}{(1-2z)^\frac{3}{2}}-[z^m]\dfrac{\tilde{b}_1(z)}{(1-2z)}.
\]
after some  computation we have 
\[
\displaystyle [z^m]\bar G_{1}(z)=2^{m}\big((m+1)\frac{\displaystyle\binom{2m}{m}}{4^m}-\frac{1}{2}\big).
\]
By replacing $ m=(n-1)/2 $ this implies
\begin{align}\label{ExactG1}
G_{1,n}= n!2^{(n-3)/2}\big((n+1)\frac{\displaystyle \binom{n-1}{(n-1)/2}}{2^{n-1}}-1\big).
\end{align}

\subsubsection{Counting Leaf-Labeled General Phylogenetic Network}\label{LtoV}
Let $ G_{n,k} $ (resp.$\tilde G_{\ell,k} $) denote the number of vertex-labeled (leaf-labeled) general phylogenetic networks with $ n $ vertices ($ \ell $ leaves) and $ k $ reticulation nodes.
It is well studied in \cite{unpublishedkeyA},  that for all subclasses of general networks containing  only networks  in which any two vertices have different sets of descendant,  we have the following equation
\begin{align}\label{GVT}
\displaystyle G_{k,n}=\binom{n}{\ell}(n-\ell)!\;\tilde{G}_{k,\ell}.
\end{align}
\begin{center}
	\begin{figure}[h]
		\begin{minipage}{1\textwidth}
			\begin{center}
				{\includegraphics[width=.8\textwidth]{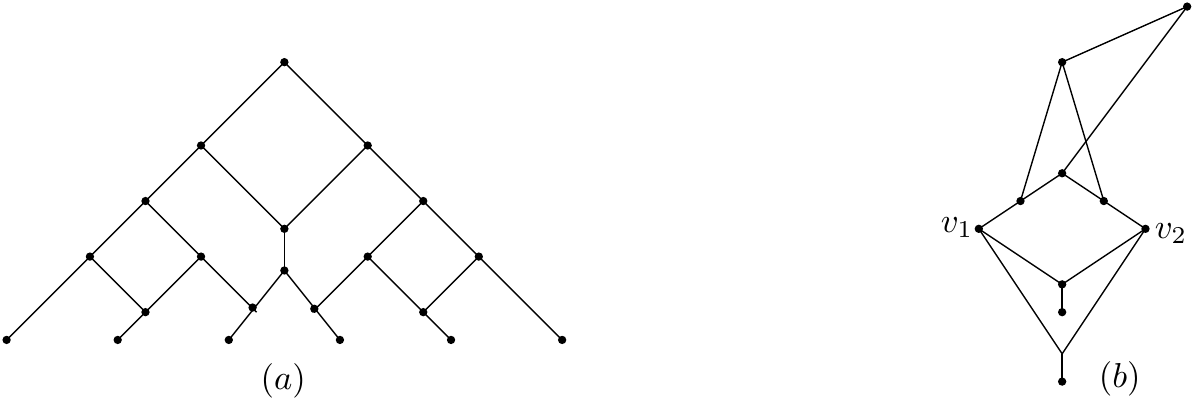}}
			\end{center}
		\end{minipage}
		\caption{Two general phylogenetic networks,  where in (a) the descendent set for any two vertices are different, and (b) is a general network which  vertices $ v_1 $ and $ v_2 $ have a set of same descendent.}
				\label{Twogen}
	\end{figure}
\end{center}

To see this, first recall from  \cite{MSW15} that for any phylogenetic network with $\ell$ leaves, $k$ reticulation vertices and $n$ vertices, we have
$\displaystyle \ell+k=\frac{n+1}{2} $ (Recall that $n$ is always odd). Now  all vertex-labeled general networks with $n$ vertices and $ k $ reticulation vertices can be constructed as follows: start with a (fixed) leaf-labeled general network with $\ell$ leaves and $k$ reticulation vertices. Then, choose $\ell$ labels from the set $n$ labels and re-label the leaves of the fixed network such that the order is preserved. Finally, label the remaining $n-\ell$ vertices by any permutation of the set of remaining $n-\ell$ labels. By the above structure, in this way every vertex-labeled general network is obtained exactly once. \\  
But for classes of networks where not all networks have the above mentioned property  it is difficult to obtain a simple connection between the vertex-labeled and  leaf-labeled phylogenetic networks. For that we have to  cope with symmetry in some of generated general networks. 
Here, we will present complete details to show how to deal with symmetry for general networks with up to $ 3 $ fixed reticulation vertices.  However, it will later be shown that as $ n $ goes infinity (resp.$ \ell $),  the family of general networks that need to deal with symmetry are asymptotically negligible and thus one again expects  $
\displaystyle \tilde{G}_{k,\ell}\sim\frac{\ell!}{n!}G_{k,n}
$, be a good  approximation   for all  leaf-labeled general networks with fixed number of reticulation vertices when $ n $ goes to infinity.  
As a warm up, we are going  to take exact formula for leaf-labeled general phylogenetic networks with one reticulation vertex.
By the above points  we get,
\[
\tilde{G}_{1,\ell}=\frac{\ell!}{n!}G_{1,n}.
\]
After seting $ n=2\ell+1 $ in (\ref{ExactG1}) we have
\begin{align}\label{ExactGL1}
\tilde{G}_{1,l}= \ell!\;2^{\ell}\big((\ell+1)\frac{\displaystyle\binom{2\ell}{\ell}}{4^\ell}-\frac{1}{2}\big).
\end{align}
\subsubsection{Relationship to Tree-child networks}\label{rel-uni}

In \cite{unpublishedkeyA}, the authors counted tree-child networks which are vertex-labeled or leaf-labeled.
On the other hand, general phylogenetic networks with exactly one reticulation vertex and no multiple edges  are tree-child networks. It means that $ G_1^{\nshortparallel}(z) $ is exactly corresponding to generating function for vertex-labeled  tree-child networks with one reticulation vertex. 
This translates into
\[
T_1(z)=G_1^{\nshortparallel}(z) =\frac{1}{2}\frac{\partial}{\partial y} \frac{z\tilde M(z,y)}{1-z M(z,y)}\Big\vert_{y=0}.
\]
Solving this equation gives
\begin{equation*}
T_1(z)=G_1^{\nshortparallel}(z) =\frac{z^3(1-\sqrt{1-2z^2})}{2(1-2z^2)^{3/2}}
\end{equation*}
as it must be. In the same way as before the mentioned approaches immediately implies that

\begin{align} \label{ExactT1}
T_{1,n}=G_{1,n}^{\nshortparallel} =\displaystyle 
n!\;2^{(n-3)/2}\big((n-1)\frac{\displaystyle \binom{n-1}{(n-1)/2}}{2^{n-1}}-1\big), 
\end{align}
and for leaf-labeled 
\begin{align}\label{T1L}
\displaystyle \tilde{T}_{1,\ell}={G}_{1,\ell}^{\nshortparallel} =\ell!\; 2^{\ell}\; \big(\ell\dfrac{\displaystyle\binom{2\ell}{\ell}}{4^{\ell}}-\frac{1}{2}\big).
\end{align} 
This approach  for leaf labeled case can be saw in \cite{unknown} with different methods. 
 \subsection{General Phylogenetic Network With Two Reticulation Vertices}
\begin{center}
 	\begin{figure}[h]
 		\begin{minipage}{1\textwidth}
 			\begin{center}
 				{\includegraphics[width=1\textwidth]{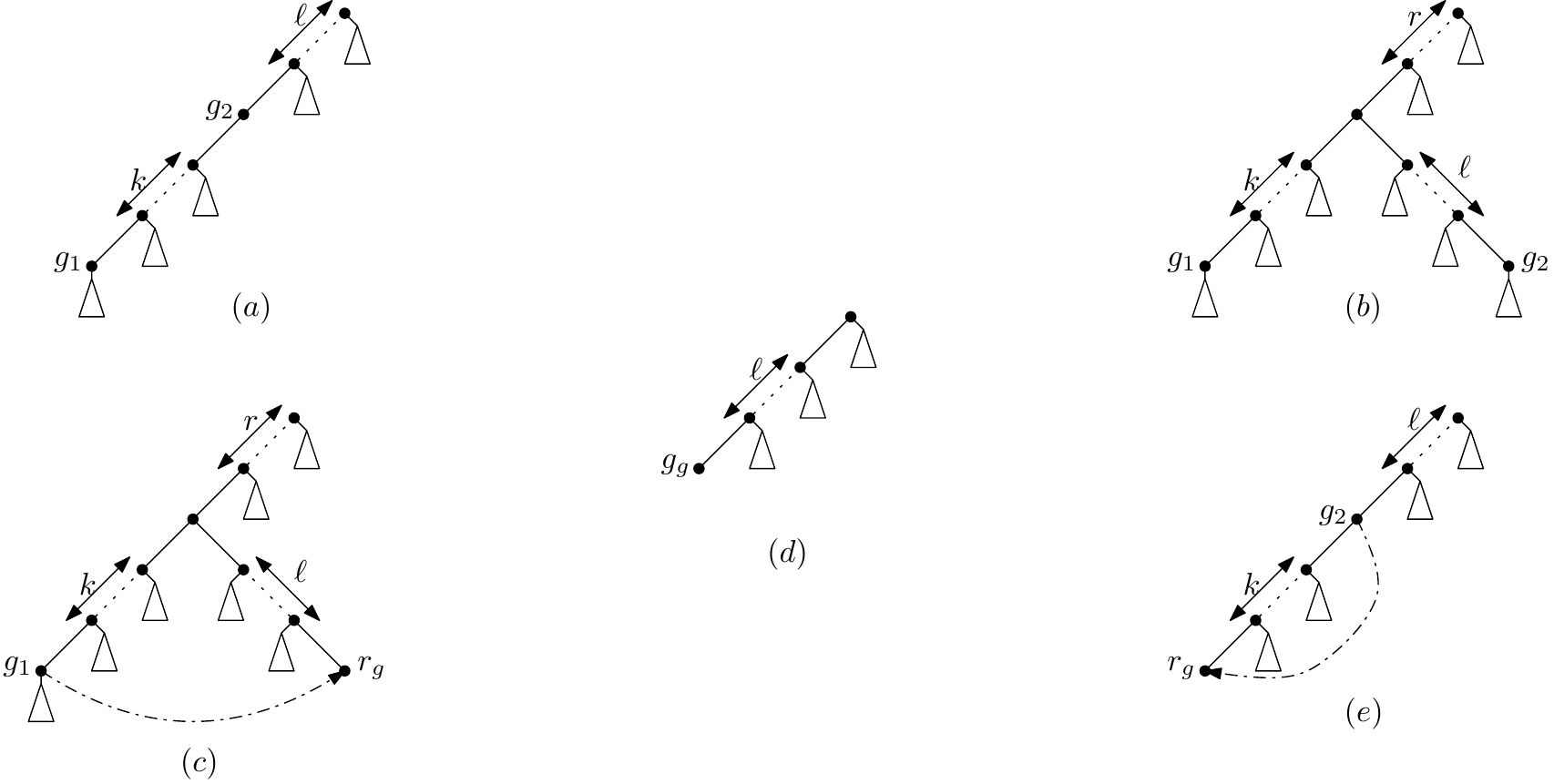}}
 			\end{center}
 		\end{minipage}
 		\caption{The possible structures of the Motzkin skeletons of general phylogenetic networks with $ 2 $ reticulation nodes such that added edges not allowed to make multiple edges.}
 		  		\label{Gen2}
 	\end{figure}
 \end{center}

\begin{center}
  	\begin{figure}[h]
  		\begin{minipage}{1\textwidth}
  			\begin{center}
  				{\includegraphics[width=1\textwidth]{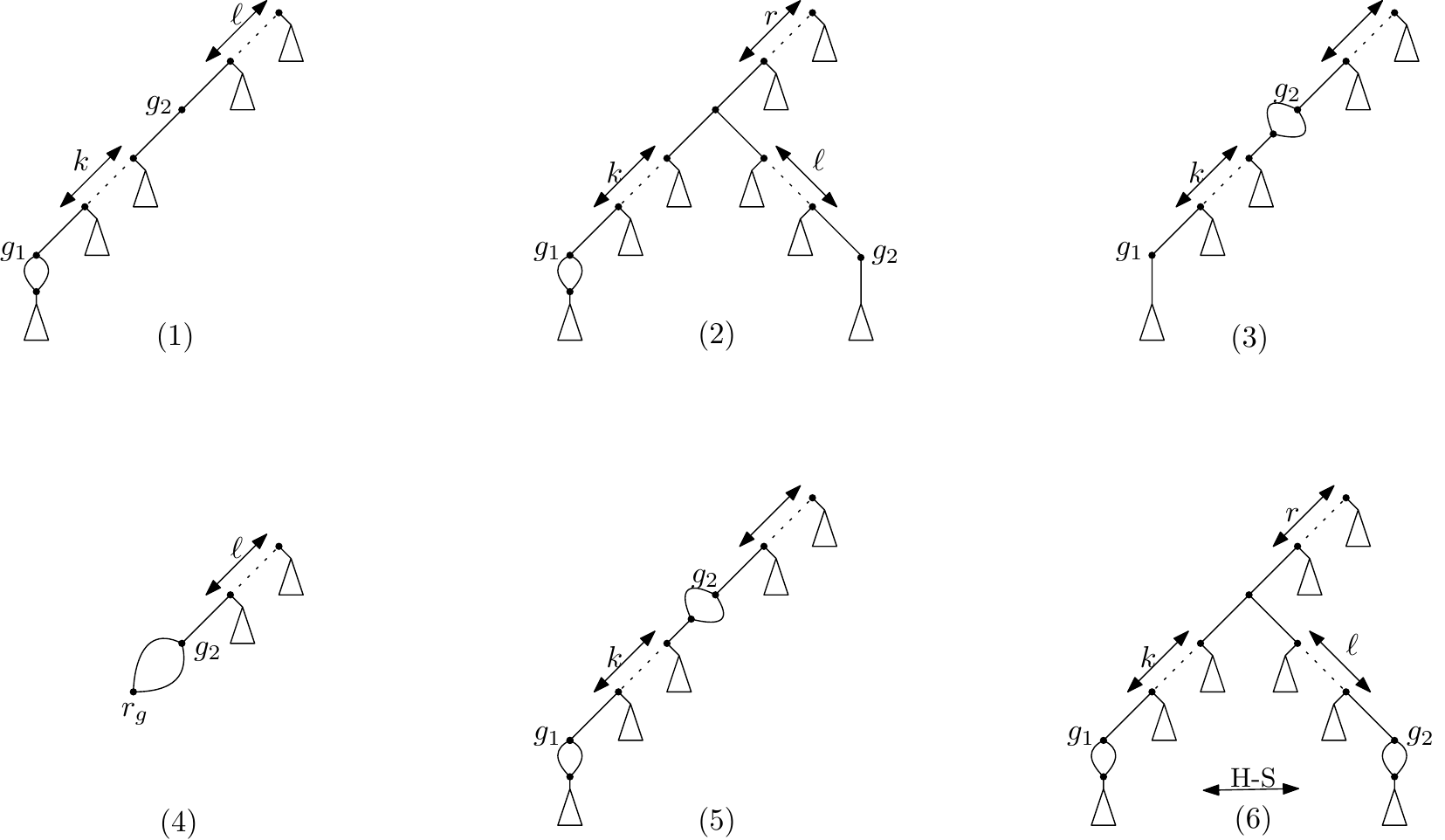}}
  			\end{center}
  		\end{minipage}
  		\caption{The possible structures of the Motzkin skeletons of  phylogenetic networks with $ 2 $ reticulation nodes and  all possible considered of fixed  multiple edges contributions.}
  		\label{PG2}
  	\end{figure}
  \end{center}

 Now we expand the methods for  general phylogenetic networks with  $ 2 $ reticulation nodes.
 For this case, we use some variables $y_1,y_2,y_{r_g},y_{g_g}$ to express the possible pointing of the  pointer set vertices of the Motzkin skeletons.
 Furthermore, we have now more complicated paths (and attached trees) which replace the edges of the sparsened skeleton and thus we first set up the generating function corresponding to theses paths. To govern the situation where an edge from one of the two pointer set  vertices must not point to a certain  vertices on the paths in order to avoid  multiple edges in the first step, we distinguish three types of unary vertices, which are the red vertices of our construction. We will  define is a class $\mathcal P$ of paths which serve as the essential building blocks for Motzkin skeletons. In this class the rules for pointing to particular red vertices differ, depending on whether (i) the red vertex lies on the path itself,(ii) it is one of the vertices of one of the attached subtrees (iii) the red vertex is the first vertex of the path. We will mark the red vertices of type (i) with the variable $y$, those of type
 (ii) with $\tilde y$ and the vertex of type (iii) with $\hat{y} $ which  is introduced in order to manage  structures analysis multiple edges in phylogenetic networks.\\
 \begin{figure}[h]
 	\begin{center}
 		\begin{tikzpicture}[
 		scale=0.6,
 		level/.style={thick},
 		virtual/.style={thick,densely dashed},
 		trans/.style={thick,<->,shorten >=2pt,shorten <=2pt,>=stealth},
 		classical/.style={thin,double,<->,shorten >=4pt,shorten <=4pt,>=stealth}
 		]
 		
 		\draw (0cm,0.8cm) node[inner sep=1.5pt,circle] (2) {};
 		\draw (-0.3cm,0.5cm) node[inner sep=1.5pt,circle] (1) {${\mathcal P}$};
 		\draw (1cm,1.8cm) node[inner sep=1.5pt,circle,draw,fill] (3) {};
 		\draw (1.7cm,1.1cm) node[regular polygon,regular polygon sides=3,draw,inner sep=0.1cm] (4) {};
 		
 		\draw (2)--(3);\draw (3)--(4);
 		
 		\draw (6.5cm,0.2cm) node[inner sep=1.5pt,circle] (2) {};
 		\draw (6.2cm,-0.1cm) node[inner sep=1.5pt,circle] (1) {${\mathcal P}$};
 		\draw (7.5cm,1.6cm) node[inner sep=1.5pt,circle,draw,fill,red!75] (5) {};
 		
 		\draw (2)--(5);
 		
 		\draw (4.1cm,1.2cm) node[inner sep=1.5pt,circle] (1) {$+$};
 		\draw (-2.2cm,1.2cm) node[inner sep=1.5pt,circle] (1) {$+$};
 		\draw (-4.2cm,1.2cm) node[inner sep=1.5pt,circle] (1) {$\{\varepsilon\}$};
 		\draw (-6.2cm,1.2cm) node[inner sep=1.5pt,circle] (1) {$=$};
 		\draw (-8.2cm,1.2cm) node[inner sep=1.5pt,circle] (1) {${\mathcal P}$};
 		\end{tikzpicture}
 	\end{center}	
 	\caption{The specification of the class $\mathcal P$.}	
 	\label{P}
 \end{figure}
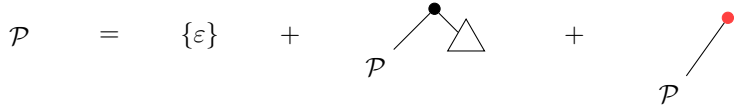 
 To simplify the explanation, let us use the following conventions where the $\varepsilon$  denotes the empty tree.   Each path is a sequence of vertices with trees attached. Note that each red vertex may  belong to different categories  (if it is first vertex of path marked with $ \hat{y} $, otherwise with $ y $). In our analysis the variables $y$,  $\tilde y$ and $ \hat{y} $ will be replaced by a
  sum of variables $y_i$ for $ i \in \{1,2, r_g,g_g\} $,  where the presence of a particular $y_i$ indicates that the corresponding $g_i$ is allowed to point,
  its absence that pointing is forbidden. In particular, $y$
  represent the permission to point to vertices of the path (except its first vertex) and
  $\tilde y$ describes the permission to point to  vertices of attached trees and $ \hat{y} $ allows pointing to the first vertex of the path.
  We specify $\mathcal P$ as
  \begin{align*}\label{specP}
  \mathcal P= \frac{1-zy}{1-z(y+M(z, \tilde{y}))}+\frac{z\hat y}{1-z(y+M(z, \tilde{y}))}.
  \end{align*}
  This leads to the generating function
  \begin{equation*}\label{PGF}
  P(z,y,\tilde y, \hat y)=\frac{1-zy+ z\hat{y}}{1-z(y+ M(z,\tilde y))},
  \end{equation*} 
  after all.
 Start with this assumption that in the Motzkin skeletons added directed edges not allowed to make multiple edges., see Figure~\ref{Gen2},  and then add the contribution of other all possible Motzkin tree skeletons with multiple edges as shown in Figure~\ref{PG2}.   Now we are ready to state the following result.
 \bp
 The exponential generating function for vertex-labeled general phylogenetic  networks with
 two reticulation nodes is 
 \begin{equation*} \label{fctN2}
 {G}_2(z)={G}_2^{\nshortparallel}(z)+{G}_2^{\shortparallel}(z)=z \cdot \frac{\tilde{a}_2(z^2)-\tilde{b}_2(z^2)\sqrt{1-2z^2}}{(1-2z^2)^{7/2}},
 \end{equation*}
 where
 \begin{equation*} \label{polyAB}
 \tilde{a}_2(z)=z^4-2z^3-\frac{1}{2}z^2+\frac{5}{2}z\quad\text{ and }\quad \tilde{b}_2(z)=-z^2+\frac{5}{2}z.
 \end{equation*}
 \ep  
 \bpf
  Consider the general phylogenetic networks arising from the Motzkin skeleton on the Figure~\ref{Gen2} $ (a) $ and complete the Motzkin skeleton by adding
  two egdes having start vertex $g_1$ and $g_2$, respectively.
  Due to this,  note that  pointings of  the green vertices do not violate the general phylogenetic network properties by  making  a  directed cyclic component.
  Also, to avoid  multiple edges, set up the generating function $ \tilde{M_{1}}(z,y_1+y_2) $ for the tree attached to the green vertex $ g_1 $. In general $ \tilde{M_{i}}(z,y_1+y_2)=(1-zy_{i})M(z,y_1+y_2) $ is the specification of unary root Motzkin trees such that pointer vertex which already marked by variable $ y_i $, is not allowed to point to the  root vertex. So this means  pointing to the root of this tree is forbidden for $ g_1 $ but not for $ g_2 $.
  For all the other trees there is no pointing restriction. The analysis of the vertices on the paths is done path by path.
  \begin{itemize}
  	\item Path $\ell$: No green vertex is allowed to point to the vertices of that path.
  	\item Path $k$: Except first node, pointing to all vertices is allowed for $g_2$, but $g_1$ may not point to that path at all. So we have 
  \end{itemize}
  \begin{align*}  	 
  G_{a}^{\nshortparallel}(z)&=
  \partial_{y_1}\partial_{y_2}z^2 \tilde M_{1} (z,y_1+y_2) P(z,y_2,y_1+y_2,0) P(z,0,y_1+y_2,0)\Big\vert_{y_1=0,y_2=0}. 
  \end{align*}
  Now, consider the Motzkin skeleton  $ (b) $. For the trees attached to the
  green vertices only pointing to the root is forbidden for parent vertices, for all the other trees there is no pointing
  restriction. The analysis of the vertices on the paths is done path by path.
  \begin{itemize}
  	\item Path $r$: No green vertex is allowed to point to the vertices of that path.
  	\item Path $k$: Pointing to all vertices is allowed for $g_2$, but $g_1$ may not point to that path at all. The situation for path $\ell$ is symmetric.
  \end{itemize}
  In this way, Motzkin skeletons which are not respecting the general condition are generated
  as well: Indeed, $ g_1 $ may point to the  vertex of $ \ell  $ and $ g_2 $ to the  vertex of $ k $, thus violating the generality condition by making a cycle. The factor $ \rcp{2} $ in the beginning of expression comes from the ``horizontal symmetry'' (This can be briefly shown by \textit{H-S}) of the Motzkin skeleton. 
  This yields the generating function
  \begin{align*}  	 
  G_{b}^{\nshortparallel}(z)&=
  \rcp2\partial_{y_1}\partial_{y_2} \left(\frac{z^3 \tilde M_{1}(z,y_1+y_2) \tilde M_{2}(z,y_1+y_2)}{1-zM(z,y_1+y_2)} P(z,y_2,y_1+y_2,y_2)
  P(z,y_1,y_1+y_2,y_1) \right.\\
  &-\left.\frac{z^3 M(z,0)^2 }{1-zM(z,0)} P(z,y_2,0,y_2)
  P(z,y_1,0,y_1)\Big\vert_{y_1=0,y_2=0}\right).
  \end{align*} 	
  The other case of general networks has the Motzkin skeleton as shown in Figure \ref{Gen2}, $ (c) $. The property of red-green leaf  entails, first one added directed edges connects $ g_1 $ to $ r_g $. After that, there is no restriction for pointing of $ r_g $ except the vertices on the paths. This gives
  \begin{align*}  	 
  G_{c}^{\nshortparallel}(z)=
  &\partial_{y_{r}}\frac{z^3M(z,y_{r})}{(1-zM(z,y_{r}))^3}\Big\vert_{y_{r_g}=0}.
  \end{align*}
  Now, consider the Motzkin skeleton $ (d)$  of Figure \ref{Gen2}. The double-green vertex $ g_g $ can point to all vertices (the pointing order does not matter, so we divide by $ 2 $) in the attached subtrees.
  \begin{align*}  	 
  G_{d}^{\nshortparallel}(z)&=\rcp2
  (\partial_{y_{g}})^2\frac{z}{1-zM(z,y_{g})}
  \Big\vert_{y_{g}=0}.
  \end{align*}
  For the final case, consider the Motzkin skeleton $ (e)$. Generality condition entails that $ r_g $ be only possible target vertex for pointing of $ g_2 $. For all the other trees there is no pointing restriction for $ r_g $. To avoid of multiple edges, the path $ k $ cannot be a simple edge.   To do that set  the generating function 
  \[ P^{\star}(z,y,\tilde y, \hat y)=P(z,y,\tilde y, \hat y)-1=\frac{z M(z,\tilde y)+ z\hat{y}}{1-z(y+ M(z,\tilde y))},\] 
  for a nonempty path.  Then we get
  \begin{align*}  	 
  G_{e}^{\nshortparallel}(z)&=
  \partial_{y_{r}}\frac{z}{1-zM(z,y_{r})}P^{\star}(z,0,y_r,0)
  \Big\vert_{y_{r}=0}=  \partial_{y_{r}}\frac{z^2M(z,y_{r}) }{(1-zM(z,y_{r}))^2}
  \Big\vert_{y_{r}=0}.
  \end{align*}
  The exponential generating function for vertex-labeled general networks (with no multiple edges) is obtained as
  $ G_{2}^{\nshortparallel}(z) = G_{a}^{\nshortparallel}(z) + G_{b}^{\nshortparallel}(z) +G_{c}^{\nshortparallel}(z) + G_{d}^{\nshortparallel}(z)+G_{e}^{\nshortparallel}(z)/4 $ after all. This gives the following result.
  \begin{align} \label{fctG2}
  {G}_2^{\nshortparallel}(z)=z \cdot \frac{a_2^{\nshortparallel}(z^2)-b_2^{\nshortparallel}(z^2)\sqrt{1-2z^2}}{(1-2z^2)^{7/2}},
  \end{align}
  where
  \begin{align}  \label{PolyG2}
  a_2^{\nshortparallel}(z)=z^4+\frac{1}{2}z^2+\frac{3}{2}z\quad\text{ and }\quad b_2^{\nshortparallel}(z)=z^2+\frac{3}{2}z.
  \end{align}
  
  Next we will set up the exponential  generating function for  general networks with at least one   multiple edges  on their structure (see Figure~\ref{PG2}). Altogether, we obtain
  \begin{align*} 
  G_{2}^{\shortparallel}(z)&= \rcp2\left(
  \partial_{y_2} z^3  M (z,y_2) \dfrac{P(z,y_2,y_2,0)}{1-zM(z,y_2)} \Big\vert_{y_2=0}+ \partial_{y_2}  \dfrac{z^4 M(z,y_2) \tilde M_{2} (z,y_2)P(z,y_2,y_2,y_2)  }{(1-zM(z,y_2))^2}\Big\vert_{y_2=0} \right.
  \\
  &\left. +\partial_{y_1} \frac{z^3 \tilde M_{1}(z,y_1)}{(1-zM(z,y_1))^2} \Big\vert_{y_1=0} +\partial_{y_r}\dfrac{z^2}{(1-zM(z,y_{r}))}\Big\vert_{y_r=0} \right)\\
  &+\dfrac{z^4 M(z,0)}{(1-zM(z,0))^2} +\rcp2 \dfrac{z^5 M^2(z,0)}{(1-zM(z,0))^3},   
  \end{align*}
  where the factor $ 2 $ appears for the expression of $ (1) $ to $ (4) $, because in these cases each general phylogenetic network  is generated two times. Note that, there is just a unique general network which arises from the  Case $ 5 $. Also, the factor $  2 $ appears in last term,  because of horizontal symmetry.  
 \epf
  So the exponential generating function for vertex-labeled general phylogenetic networks with
   two reticulation nodes is then $G_{2}(z)= G_{2}^{\nshortparallel}(z)+G_{2}^{\shortparallel}(z) $.
   As an easy consequence, we obtain the asymptotic number of networks.  
 \bcor
 Let $G_{2,n}$ denote the number of vertex-labeled  general phylogenetic networks with $n$ vertices
 and exactly two reticulation vertex. If $n$ is even then $G_{2,n}$ is zero, otherwise 
 \[
 G_{2,n}=n![z^n]G_{2}(z)=\(\frac{\sqrt2}{e}\)^n n^{n+3}\(\frac{\sqrt{2}}{16}-\frac{\sqrt{\pi}}{8}\cdot\frac{1}{\sqrt{n}}
 +{\mathcal O}\left(\frac{1}{n}\right)\),
 \]
 as $\nti$.
 \ecor
 
 \bpf
 This follows by singularity analysis as before.
 \epf
 \subsubsection{Explicit formula for vertex-labeled general networks with two reticulation vertices} 
  We can use generating functions $ G_{2}^{\nshortparallel}(z)$ and    $G_{2}^{\shortparallel}(z) $  to extract closed  formulas for vertex-labeled general networks. To see that, consider the contribution of each of them separately. Start with the exponential  generating function $ G_{2}^{\nshortparallel}(z) $ for general networks that do not have double edges in own structures.

  First, set $n=2m+1$. Then, from \ref{fctG2} we obtain 
  \[
  [z^n]{G}_2^{\nshortparallel}(z)=[z^m]\bar{G}_2^{\nshortparallel}(z)
  \]
  with 
  \[
  \bar{G}_2^{\nshortparallel}(z)= \frac{a_2^{\nshortparallel}(z)-b_2^{\nshortparallel}(z)\sqrt{1-2z}}{(1-2z)^{7/2}},
  \]
  where $a_2^{\nshortparallel}(z)$ and $b_2^{\nshortparallel}(z)$ are as in \ref{PolyG2}.
  So we have
  \[
  [z^m]\bar G_{2}^{\nshortparallel}(z)=[z^m]\dfrac{a_2^{\nshortparallel}(z)}{(1-2z)^\frac{7}{2}}-[z^m]\dfrac{b_2^{\nshortparallel}(z)}{(1-2z)^3}.
  \]
  After some  computation we have 
  \[
  [z^m]\bar{G}_2^{\nshortparallel}(z)=2^{m-2}\Big(P_1(m)\frac{\displaystyle 2m\binom{2m}{m}}{15(2m-1)4^m}-P_2(m)\Big),
  \]
  where 
  \begin{align} \label{polyP1P2}
  P_1(m)=30m^3+20m^2+15m-20\quad\text{ and }\quad P_2(m)=2m^2+m.
  \end{align}
   By replacing $ m=(n-1)/2 $ this implies
  \begin{align}\label{ExactG2}
  G_{2,n}^{\nshortparallel}= n!2^{(n-5)/2}\Big(P_1((n-1)/2)\frac{\displaystyle (n-1) \binom{n-1}{(n-1)/2}}{15(n-2)2^{n-1}}-P_2((n-1)/2)\Big).
  \end{align}
  Note that correspondent  generating function   for general networks with multiple edges is 
  \begin{align*} 
  { G}_2^{\shortparallel}(z)=z \cdot \frac{ a_2^{\shortparallel}(z^2)- b_2^{\shortparallel}(z^2)\sqrt{1-2z^2}}{(1-2z^2)^{1/2}},
  \end{align*}
  such that
  \begin{align*} 
  a_2^{\shortparallel}(z)=z^2+z\quad\text{ and }\quad  b_2^{\shortparallel}(z)=z.
  \end{align*} 
  In the same way, it   can be used to get  exact formula for vertex-labeled general networks that are belong to this subclass. We refrain from giving details and just list the obtained expressions. The reader is invited to derive them herself.
  \begin{align}\label{ExactGP2}
  G_{2,n}^{\shortparallel}= n!2^{(n-3)/2}(n-1)\Big(\frac{\displaystyle (n-1) \binom{n-1}{(n-1)/2}}{2^{n}}-\rcp{2})\Big).
  \end{align}
  
  After all by summing up \ref{ExactG2} and \ref{ExactGP2}  we have
  \begin{align} \label{ExactGen2}
  G_{2,n}=G_{2,n}^{\nshortparallel}+G_{2,n}^{\shortparallel}
  = n!2^{(n-3)/2}\Big(A((n-1)/2)\frac{\displaystyle (n-1) \binom{n-1}{(n-1)/2}}{15(n-2)2^{n-1}}-B((n-1)/2)\Big),
  \end{align}
  where
  \begin{align} 
  A(m)=30m^3+80m^2-15m-20\quad\text{ and }\quad B(m)=m^2+\frac{3}{2}m.
  \end{align}
  \subsubsection{Explicit formula for leaf-labeled general networks with two reticulation vertices} 
\begin{center}
  	\begin{figure}[h]
  		\begin{minipage}{1\textwidth}
  			\begin{center}
  				{\includegraphics[width=.6\textwidth]{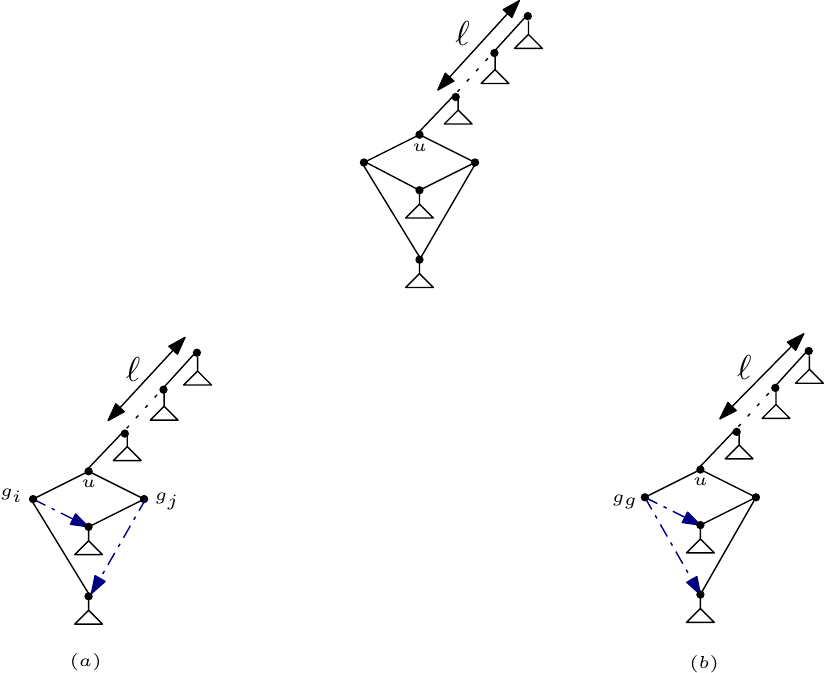}}
  			\end{center}
  		\end{minipage}
  		\caption{(Top) The structure of general network  with two reticulation vertices such that two vertices have the same set of descendant which can be generate from (a) by pointing green vertices to the root of each other attached trees or (b) a double-green vertex points to unary vertices with same parent.}
  		\label{Sym}
  	\end{figure}	
  \end{center}

  Note that, the Equation (\ref{GVT}) which comes from the described  procedure in Section \ref{LtoV}  for construction all vertex-labeled networks  from fixed leaf-labeled ones does not work anymore.  It is because by applying the method there are some leaf-labeled networks which generate some  vertex-labeled networks more than one (here twice).
  Thus for normalization, and deal with symmetry the correspondent  generating functions of such networks can be considered separately (see Figure \ref{Sym}).
  So we have 
  \[{G}_2^{\nshortparallel}(z)={\dot G}_2^{\nshortparallel}(z)+{\ddot  G}_2^{\nshortparallel}(z),  \]
  where $ {\dot G}_2^{\nshortparallel}(z) $ correspondent  generating function  for general networks such that  the procedure set out in Section \ref{LtoV} can be applied directly for them. 
  \begin{align*} 
  {\dot G}_2^{\nshortparallel}(z)=z \cdot \frac{\dot a_2^{\nshortparallel}(z^2)-\dot b_2^{\nshortparallel}(z^2)\sqrt{1-2z^2}}{(1-2z^2)^{7/2}},
  \end{align*}
  where,
  \begin{align*} 
  \dot a_2^{\nshortparallel}(z)=-4z^5+11z^4-9z^3+4z^2+z\quad\text{ and }\quad \dot b_2^{\nshortparallel}(z)=4z^4-6z^3+4z^2+z.
  \end{align*}
  For this set of  general networks  we can directly use  Equation (\ref{GVT}). Thus, same procedure as like before gives us
  \begin{align}
  \dot G_{2,\ell}^{\nshortparallel}= \ell!2^{\ell-1}\Big((6\ell^4+19\ell^3+18\ell^2-7\ell-3)\frac{\displaystyle (\ell+1) \binom{2\ell+2}{\ell+1}}{(6\ell-3)(2\ell+1)4^{\ell}}-(2\ell^2+5\ell+3)\Big).
  \end{align} 	 
  Now we set up  generating function for family of networks which are shown in the top of Figure \ref{Sym}. It is not hard to see that, by using the previous methods each related (fixed) leaf-labeled general network can construct a vertex-labeled general network exactly twice (because of symmetry). For this case Equation  (\ref{GVT}) can modify as    $ {\ddot G}^{\nshortparallel}_{2,\ell}=2\frac{\ell!}{n!}{\ddot G}^{\nshortparallel}_{2,n}. $ The generating function for this subfamily of general networks is 
  \begin{align*} 
  {\ddot G}_2^{\nshortparallel}(z)=z \cdot \frac{\ddot a_2^{\nshortparallel}(z^2)-\ddot b_2^{\nshortparallel}(z^2)\sqrt{1-2z^2}}{(1-2z^2)^{1/2}},
  \end{align*}
  where
  \begin{align*} 
  \ddot a_2^{\nshortparallel}(z)=-\rcp{2}z^2+\rcp{2}z\quad\text{ and }\quad \ddot b_2^{\nshortparallel}(z)=\rcp{2}z.
  \end{align*}
  After some manipulation we get 	
  \begin{align*}
  \ddot G_{2,\ell}^{\nshortparallel}= \ell!2^{\ell-1}\Big(\frac{\displaystyle (\ell+1) \binom{2\ell+2}{\ell+1}}{(2\ell-1)(2\ell+1)4^{\ell}}\Big).
  \end{align*}  	
  The explicit formula  for leaf-labeled general networks with no multiple edges and two reticulation vertices is	 
  \begin{align}\label{ExactGL2}
  \begin{aligned} 
  G_{2,\ell}^{\nshortparallel}&=\dot G_{2,\ell}^{\nshortparallel}+\ddot G_{2,\ell}^{\nshortparallel} \\
  &=\ell!2^{\ell-1}\Big((6\ell^4+19\ell^3+18\ell^2-7\ell)\frac{\displaystyle (\ell+1) \binom{2\ell+2}{\ell+1}}{(6\ell-3)(2\ell+1)4^{\ell}}-(2\ell^2+5\ell+3)\Big).
  \end{aligned} 
  \end{align}		  
  To complete the details, we can get explicit formula for the number of leaf-labeled networks that are generated by sparsened skeletons which as depicted in  Figure \ref{PG2}. Note that for this case, all generated networks belong to the first subclass of general networks which   the Equation \ref{GVT} can be used directly. So we have
  \begin{align*}
  G_{2,\ell}^{\shortparallel}= \ell!2^{\ell+1}(\ell+1)\Big(\frac{\displaystyle (\ell+1) \binom{2\ell+2}{\ell+1}}{4^{\ell+1}}-\rcp{2}\Big).
  \end{align*}	   
  
  After all, we get
  \begin{align}\label{G2L}
  \displaystyle \tilde{G}_{2,\ell}=G_{2,\ell}^{\nshortparallel}+G_{2,\ell}^{\shortparallel}=\ell!2^{\ell-1}\Big( A(\ell)\frac{\displaystyle (\ell+1) \binom{2\ell+2}{\ell+1}}{(6\ell-3)(2\ell+1)4^{\ell}}-B(\ell)\Big)\Big),
  \end{align}	   
  where
  \begin{align} 
  A(\ell)=6\ell^4+31\ell^3+30\ell^2-10\ell-3\quad\text{ and }\quad  B(\ell)=2\ell^2+\frac{41}{8}\ell+\frac{25}{8}.
  \end{align} 
 In the same way, the methods can be used to study of specifications for general phylogenetic networks with $ k\geq 3 $  reticulation nodes. Its obvious by increasing number of rediculation nodes we have to consider various number  of the  Motzkin skeletons to cover all possible cases. For more understanding, we invite the reader to look at Appendix section to see all the Motzkin skeletons and related specifications for $ k=3 $. As similar as $ k=2 $,   first we focus on  structures with no  multiple edges and then  for each Motzkin skeleton consider possible contributions of double edges on the structures and add them to the results.
After all this gives  $  G_3(z)=G_3^{\nshortparallel}(z)+G_3^{\shortparallel}(z)$ ( See Appendix A for more details) such that causes following results for vertex-labeled general networks with  fixed $ 3 $ reticulation vertices.
\bp
The exponential generating function for vertex-labeled  general phylogenetic networks with
three reticulation nodes is 
\begin{equation*} 
{G}_3(z)=G_3^{\nshortparallel}(z)+G_3^{\shortparallel}(z)= z \cdot \frac{a_3(z^2)-b_3(z^2)\sqrt{1-2z^2}}{(1-2z^2)^{11/2}},
\end{equation*}
where 
\begin{equation*} 
a_3(z)=z^6+5z^5-10z^4-\frac{23}{2}z^3+\frac{109}{4}z^2,
\end{equation*}
and,
\begin{equation*} 
b_3(z)=z^5-\frac{7}{2}z^4-5z^3+\frac{109}{4}z^2.\;\;\;\;\;
\end{equation*}
\ep
\bcor\label{fctN8}
Let $G_{3,n}$ denote the number of vertex-labelled  general phylogenetic networks with $n$ vertices
and exactly three reticulation vertex. If $n$ is even then $G_{3,n}$ is zero, otherwise 
\[
G_{3,n}=n![z^n]G_{3}(z)=\(\frac{\sqrt2}{e}\)^n n^{n+5}\(\frac{\sqrt{2}}{192}-\frac{\sqrt{\pi}}{64}\cdot\frac{1}{\sqrt{n}}
+{\mathcal O}\left(\frac{1}{n}\right)\),
\]
as $\nti$.\\
\ecor
Also,  consequently as similar as before we can take the explicit formulas for vertex and leaf-labeled general networks with $ 3 $ reticulation vertices. For vertex labelled case, as like before we set $n=2m+1$, so we have
\[
[z^n]{G}_3(z)=[z^m]\bar{G}_3(z),
\]
such that,
\[
[z^m]\bar G_{3}(z)=[z^m]\dfrac{a_{3}(z)}{(1-2z)^\frac{11}{2}}-[z^m]\dfrac{b_{3}(z)}{(1-2z)^5}.
\]
It gives 
\[
\mathcal{F}(m):=[z^m]\bar{G}_3(z)=\frac{2^{m-6}}{3}\Big(A_1(m)\frac{\displaystyle m (m-1)\binom{2m}{m}}{35(2m-1)4^{m-2}}-B_1(m)\Big),
\]
where 
\begin{equation} \label{polyPQG21}
A_1(m)=104m^4+836m^3+876m^2-454m-79 \quad \text{and,}\quad B_1(m)=48m^4+127m^3-60m^2-121m+6.
\end{equation}
By replacing $ m=(n-1)/2 $, we have,
\[
G_{3,n}= n!\cdot \mathcal{F}((n-1)/2).
\]
we need some more arguments to extract explicit formula for the leaf-labeled case. It is because of symmetry that we can see in some of the generated networks. 
For someone who is interest, complete details of steps can be found in the Appendix section.\\
 Now, the defined structure for paths of sparsened skeletons with  well defined generating function (\ref{fctM}) for attached trees, capable us to   prove the theorem \ref{thm-G1}.\\
\textit{Proof of Theorem \ref{thm-G1}}.
In particular note that function $ G(z,y) $ is  the form $ zM $ (\ref{fctM}), which $ z $ refers to vertices lie on the pathes of sparsened skeleton.
\begin{align}\label{func-G}
G(z,y)=a(z,y)-b(z,y)\sqrt{1+(y^2-2)z^2-2zy},
\end{align}
where $a(z,y),b(z,y)$ are polynomials in $z$ and $y$ with $ a(z, 0) = b(z, 0) = 1 $.
In  summary,  we have exponential generating function $ G_k $ for phylogenetic network in sum of terms of the form
\begin{align}\label{terms-N}
\partial_{y_1}\cdots\partial_{y_k}\frac{G_1(z,y)\cdots G_{s}(z,y)}{(1-G_{s+1}(z,y))\cdots(1-G_{s+t}(z,y))}\Big\vert_{y_1=0,\ldots,y_k=0},
\end{align}
Note that in this expression, numerator refers to generating function of subtrees which rooted at green vertices. The  denominator refers to sequences of subtrees which rooted the vertices on the paths of sparsened skeleton.  Also where the number of functions $G_{s+i}(z,y)$ is bounded by the number of edges of the sparsened skeleton increased by one (for the sequence of trees added above the root when constructing the Motzkin skeletons).  Now,  recall lemma $ 3.1 $ from  \cite{unpublishedkeyA} which  can be used for any similar structures as $ G(z,y)  $. 
We can apply this lemma after expanding (\ref{terms-N}) and obtain that
\begin{align}\label{gen-form}
G_k(z)=\frac{a_k(z)-b_k(z)\sqrt{1-2z^2}}{(1-2z^2)^{p}}.
\end{align}

We proceed to show that $p=2k-1/2$. For that,  observe  (\ref{terms-N}) without the
derivatives is of the general form given in (\ref{gen-form}) with the exponent of the denominator
equals $t/2$ which reaches its maximum for the sparsened skeleton with the maximal number of edges and is thus at most $k-1/2$. Also, from the above lemma, we see that each differentiation increases the exponent
by $1$. Thus, the exponent of (\ref{terms-N}) when written as (\ref{gen-form}) is at most $2k-1/2$. Adding up this terms gives   

\[
G_{k}(z)=\frac{{a}_k(z)-{b}_k(z)\sqrt{1-2z^2}}{(1-2z^2)^{2k-1/2}},
\]
where $a_k(z)$ and $b_k(z)$ are suitable polynomials.
Let $G_{k,n}$ denote the number of vertex-labeled general phylogenetic networks with $n$ vertices and $k$ reticulation vertices. If $n$ is even then $G_{k,n}$ is zero, otherwise there is a positive constant ${d}_{k}$ such that  as $ n\rightarrow\infty $,
\[
G_{k,n}=n![z^n]G_k(z)\sim {d}_k\left(\frac{\sqrt{2}}{e}\right)^n n^{n+2k-1}, 
\]
Where by singularity analysis and Stirling's formula we get
\[
{d}_k=\frac{2\sqrt{2\pi}{a}_k(1/\sqrt{2})}{4^k\Gamma(2k-1/2)}.
\]
\begin{flushright}
	$ \square $
\end{flushright}
\brem 
For the positivity claim, we already see in \cite{unpublishedkeyA} that 
corresponding constant $ \tilde{d}_k $  for normal and tree-child networks is positive  which is  lower bound of $ d_k $ for general networks. 
\erem
\bp
\it For the numbers of vertex-labeled  general phylogenetic networks $G_{k,n}$ and vertex-labeled tree-child networks $T_{k,n}$,
\begin{align}\label{G-T}
G_{k,n}=T_{k,n}\left({1+\mathcal O}(\frac{1}{n})\right),
\end{align}
as $\nti$.

\ep

\begin{proof}
First, observe that $ G_{k,n}-T_{k,n} $ is bounded by the number of networks which arise from all
types of Motzkin skeletons where for each green vertex, the considered all possibilities of adding an edge violates  the tree-child condition. Note that, the largest number will come from the sparsened skeletons where all pointer vertices are the leaves. 
Now, fix such a type of Motzkin skeletons and one of its green vertices. Then, for this vertex, we will have the following options.
\begin{itemize}
	\item The green vertex points to the  root of the subtree which is attached  to the one of green vertices in  the Motzkin skeletons. Note that if it points  to the root of its subtree, tree-child condition violates by making multiple edge.  For the exponential generating function this gives
	\begin{equation*}
	\partial_{y_2}\cdots\partial_{y_k}\frac{G_1(z,y)\cdots G_{s}(z,y)}{(1-G_{s+1}(z,y))\cdots(1-G_{s+2k-1}(z,y))}\Big\vert_{y_2=0,\ldots,y_k=0},
	\end{equation*}
	Here, and below $y$ is the sum of $y_i$'s with $2\leq i\leq k$ and not all of the $y_i$'s must be present; also which are present can differ from one occurrence to the next.
	\item There is a  red-green vertex on the Motzkin skeleton. Note that the red-green property entails that one another pointer vertex joints to this leaf by adding directed edge which reduces the number of the derivative by one. Then we get
	\begin{equation*}
	\partial_{y_2}\cdots\partial_{y_k}\frac{G_1(z,y)\cdots G_{s-1}(z,y)}{(1-G_{s+1}(z,y))\cdots(1-G_{s+2k-1}(z,y))}\Big\vert_{y_2=0,\ldots,y_k=0}.
	\end{equation*}
	\item There is double-green vertex in the Motzkin skeleton that points to the  branches of sparsened skeleton. Then, we have
	\begin{equation*}
	\partial_{y_3}\cdots\partial_{y_k}\frac{G_1(z,y)\cdots G_{s-2}(z,y)}{2\cdot(1-G_{s+1}(z,y))\cdots(1-G_{s+2k-1}(z,y))}\Big\vert_{y_2=0,\ldots,y_k=0}.
	\end{equation*}
	The existence of double green node in considered skeleton  is like that two green vertices are merged to each others. Consequently,  the number of edges reduce by two, which also leads to a contribution of smaller order.
\end{itemize}
\begin{center}
	\begin{figure}[h]
		\begin{minipage}{1\textwidth}
			\begin{center}
				{\includegraphics[width=.6\textwidth]{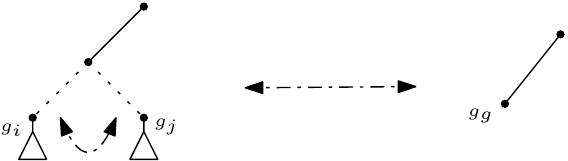}}
			\end{center}
		\end{minipage}
		\label{merge}
	\end{figure}	
\end{center} 

The exponential generating function of all networks arising from these Motzkin skeletons and the pointer vertices are a sum of generating functions of the above three types. Thus, we obtain that this generating function has the form
\[
\frac{c(z)-d(z)\sqrt{1-2z^2}}{(1-2z^2)^{p}},
\]
where $c(z)$ and $d(z)$ are suitable polynomials and the maximum of $p$ is as follows: note that without the derivatives in the above expressions, $p$ would be at most $k-1/2$. Also, because of above lemma, each derivative increases this bound by one. Thus, $p$ is at most $2k-\frac{3}{2}$.

Now, we obtain that the exponential generating function of the above number has the form
\[
\frac{{c}(z)-{d}(z)\sqrt{1-2z^2}}{(1-2z^2)^{2k-\frac{3}{2}}},
\]
where ${c}(z)$ and ${d}(z)$ are suitable polynomials. Singularity analysis gives then the bound
\[
{\mathcal O}\left(\left(\frac{\sqrt{2}}{e}\right)^nn^{n+2k-2}\right).
\]
Summing over all possible type of Motzkin skeletons and all green vertices, we obtain  the claimed result.
\end{proof}
\begin{table}[ht]
	\centering
	\begin{tabular}{|l||l|l|l||}
		\hline
		Vertex Labeled	& \;\;\;\;\;\;\;\; k=1 & \;\;\;\;\;\;\;\;\; k=2\;\;\;\;\;\; &\;\;\;\;\;\;\;\;\; k=3\;\;\;\;\;\;  \\ \cline{2-4} 
		Phylogenetic 	&  &  &  \\
		
		Networks	&$ c_1 $ \;\;\;\;\;\;\;\;\;\;$ c_1^\prime $  &$ c_2 $ \;\;\;\;\;\;\;\;\;\;$ c_2^\prime $  &$ c_3 $ \;\;\;\;\;\;\;\;\;\;$ c_3^\prime $    \\
		\hline
		&  &  &  \\
		$N_{k,n} $ & $ \frac{\sqrt{2}}{2} $ \;\;\;\;\; $ -\frac{3\sqrt{\pi}}{2} $ &$\frac{\sqrt{2}}{16} $ \;\;\;\;\; $ -\frac{3\sqrt{\pi}}{8} $ & $\frac{\sqrt{2}}{192} $ \;\;\;\;\; $ -\frac{3\sqrt{\pi}}{64} $   \\
		$  $ &  &  &        \\
		\hline
		&  &    &  \\
		$T_{k,n} $ & $ \frac{\sqrt{2}}{2} $ \;\;\;\;\; $ -\frac{\sqrt{\pi}}{2} $ & $ \frac{\sqrt{2}}{16} $ \;\;\;\;\; $ -\frac{\sqrt{\pi}}{8} $  & $\frac{\sqrt{2}}{192} $ \;\;\;\;\; $ -\frac{\sqrt{\pi}}{64} $     \\
		$  $ & &  &       \\
		\hline
		&  &  &   \\
		$G_{k,n}$ & $ \frac{\sqrt{2}}{2} $ \;\;\;\;\; $ -\frac{\sqrt{\pi}}{2} $ & $ \frac{\sqrt{2}}{16} $ \;\;\;\;\; $ -\frac{\sqrt{\pi}}{8} $  & $\frac{\sqrt{2}}{192} $ \;\;\;\;\; $ -\frac{\sqrt{\pi}}{64} $      \\
		& &  &       \\
		\hline
	\end{tabular}\caption{ The first two asymptotic orders of normal, tree-child and general phylogenetic networks with at most $ 3 $ reticulation vertices. For all of them  the first coefficient is same.   \label{TableNumbers}}
\end{table}

\subsection{Asymptotic counting of leaf-labeled general phylogenetic Networks }
In this part we want to prove Theorem \ref{thm-G2} and argue that for  the number of leaf-labeled general phylogenetic networks with $ k\geq 1 $ reticulation vertices (as  like leaf-labeled tree-child ($ \tilde{T}_{k,\ell} $) and normal networks, see \cite{unpublishedkeyA}) we  can use, 
\begin{align}\label{GA}
\tilde{G}_{k,\ell}\sim  2^{3k-1}d_k\left(\frac{2}{e}\right)^{\ell}\ell^{\ell+2k-1},\qquad(\ell\rightarrow\infty)
\end{align}
as a  relative precise estimate of leaf-labeled  general phylogenetic network, where $ d_k $ is as in Theorem \ref{thm-G1}.

\begin{center}
	\begin{figure}[h]
		\begin{minipage}{1\textwidth}
			\begin{center}
				{\includegraphics[width=.7\textwidth]{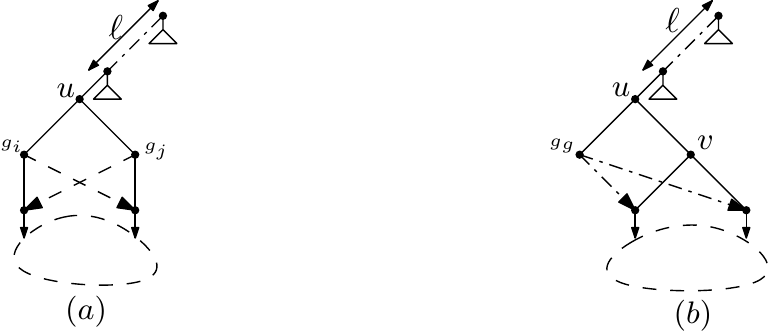}}
			\end{center}
		\end{minipage}
		\caption{The structures of general phylogenetic networks where pair of vertices have a same descendant set after adding the directed edges in Motzkin skeletons.}
		\label{LG}
	\end{figure}
\end{center}

It is enough to show that the number of a subfamily $ \mathcal{G} $ of  general networks such that some of their vertices   have the same set of descendant
are rare.
Indeed, $ \mathcal{G} $ consists of general networks that equation $(6) $  can not be used directly for them.  It's because of that in the described method (see \ref{LtoV} ) some of the fix leaf-labeled networks generate vertex-labeled networks more than one.
  In other words, having a pair of vertices with a set of the same descendant is necessary condition but not sufficient  to generate vertex-labeled networks twice or more. For instance, consider a leaf which is attached edge $ (u, g_i) $ in Figure \ref{LG} (a). Though, $ g_1 $ and $ g_2 $ have a set of the same descendant but applying the procedure (\ref{LtoV}), generates each vertex-labeled uniquely.
     

\textit{Proof of Theorem\ref{thm-G2}}.
Consider a subfamily $ \mathcal{G} $ as similar before. It is sufficient for our purposes to show that when $ \ell\rightarrow \infty$, the number of these networks are  asymptotically negligible.
Assume, without loss of generality, these networks are without multiple edges because each of them reduces the number of differentiations in
the expression for the exponential generating function by one, that causes the contribution of lower-order.\\
Note that, $ \mathcal{G} $ is bounded by the number of networks which
arise from two types of Motzkin skeletons that are depicted in Figure \ref{LG}.
First, when two green vertices point to the child vertices of each others (Figure \ref{LG}, (a)) and second, a double-green vertex points unary vertices with the same parent (b). In the former case, two green vertices and in the later case double-green vertex with vertex $v$ have set of the same descendant.  
Note that in each of described cases,  the number of derivatives and consequently, the power of denominator in exponential generating function will be reduce by two. So The first two asymptotic orders are as in theorem \ref{thm-G1}.
That implies
\begin{align}\label{GL}
G_{k,2\ell+2k-1}\sim\binom{2\ell+2k-1}{\ell}(\ell+2k-1)!\tilde{G}_{k,\ell}.
\end{align}
Now we have $ \tilde{G}_{k,\ell}\sim \frac{\ell!}{(2\ell+2k-1)!}G_{k,2\ell+2k-1} $, which an asymptotic result (\ref{GA})  follows by Theorem \ref{thm-G1} and Stirling's formula.
\begin{flushright}
	$ \square $
\end{flushright}
 \section*{Acknowledgment} I would like to thank the 
 Bernhard Gittenberger and Michael Fuchs for critically reading the manuscript and making suggestions that led to significant improvements to the content and clarity of the paper.

\bibliographystyle{plain}
\bibliography{Gen}

\begin{thebibliography}{10}

\bibitem{AlAl16}
Nikita Alexeev and Max~A. Alekseyev.
\newblock Combinatorial scoring of phylogenetic networks.
\newblock In {\em Computing and combinatorics}, volume 9797 of {\em Lecture
  Notes in Comput. Sci.}, pages 560--572. Springer, [Cham], 2016.

\bibitem{Bandelt1994}
photo Hans-Jurgen Bandelt.
\newblock Phylogenetic networks.
\newblock In {\em Verhandlungen des Naturwissenschaftlichen Vereins Hamburg},
  volume~34, pages 51--71, 1994.

\bibitem{Bo16}
Mikl\'os B\'ona.
\newblock On the number of vertices of each rank in {$k$}-phylogenetic trees.
\newblock {\em Discrete Math. Theor. Comput. Sci.}, 18(3):Paper No. 7, 7, 2016.

\bibitem{BoFl09}
Mikl\'os B\'ona and Philippe Flajolet.
\newblock Isomorphism and symmetries in random phylogenetic trees.
\newblock {\em J. Appl. Probab.}, 46(4):1005--1019, 2009.

\bibitem{BoSe16a}
Magnus Bordewich and Charles Semple.
\newblock Determining phylogenetic networks from inter-taxa distances.
\newblock {\em J. Math. Biol.}, 73(2):283--303, 2016.

\bibitem{MPM}
Mathilde Bouvel, Philippe Gambette, and Marefatollah Mansouri.
\newblock Counting phylogenetic networks of level 1 and 2. arxiv:1909.10460v2.
\newblock 2019.

\bibitem{Scornavacca}
Gabriel Cardona, Joan~Carles Pons, and Celine Scornavacca.
\newblock Correction: Generation of binary tree-child phylogenetic networks.
\newblock {\em PLOS Computational Biology}, 15(10):1--1, 10 2019.

\bibitem{CEJM13}
\'{E}va Czabarka, P\'{e}ter~L. Erd\H{o}s, Virginia Johnson, and Vincent
  Moulton.
\newblock Generating functions for multi-labeled trees.
\newblock {\em Discrete Appl. Math.}, 161(1-2):107--117, 2013.

\bibitem{DiRo17}
Filippo Disanto and Noah~A. Rosenberg.
\newblock Enumeration of ancestral configurations for matching gene trees and
  species trees.
\newblock {\em J. Comput. Biol.}, 24(9):831--850, 2017.

\bibitem{AnaCombi}
Philippe Flajolet and Robert Sedgewick.
\newblock {\em Analytic combinatorics}.
\newblock Cambridge University Press, Cambridge, 2009.

\bibitem{FoRo80}
Leslie~R. Foulds and Robert~W. Robinson.
\newblock Determining the asymptotic number of phylogenetic trees.
\newblock In {\em Combinatorial mathematics, {VII} ({P}roc. {S}eventh
  {A}ustralian {C}onf., {U}niv. {N}ewcastle, {N}ewcastle, 1979)}, volume 829 of
  {\em Lecture Notes in Math.}, pages 110--126. Springer, Berlin, 1980.

\bibitem{FoRo88}
Leslie~R. Foulds and Robert~W. Robinson.
\newblock Enumerating phylogenetic trees with multiple labels.
\newblock In {\em Proceedings of the {F}irst {J}apan {C}onference on {G}raph
  {T}heory and {A}pplications ({H}akone, 1986)}, volume~72, pages 129--139,
  1988.

\bibitem{unpublishedkeyA}
Michael Fuchs, Bernhard Gittenberger, and Marefatollah Mansouri.
\newblock Counting phylogenetic networks with few reticulation vertices:
  Tree-child and normal networks.
\newblock {\em Australasian Journal of Combinatorics}, (73)2:385--423, 2019.

\bibitem{GEL2003}
Dan Gusfield, Satish Eddhu, and Charles Langley.
\newblock Efficient reconstruction of phylogenetic networks with constrained
  recombination.
\newblock In {\em CSB03}, pages 363--374, 2003.
\newblock http://wwwcsif.cs.ucdavis.edu/~gusfield/ieeefinal.pdf.

\bibitem{huson2010phylogenetic}
D.H. Huson, R.~Rupp, and C.~Scornavacca.
\newblock {\em Phylogenetic Networks: Concepts, Algorithms and Applications}.
\newblock Cambridge University Press, 2010.

\bibitem{LinderRieseberg2004}
C.~Randal Linder and Loren~H. Rieseberg.
\newblock Reconstructing patterns of reticulate evolution in plants.
\newblock {\em American Journal of Botany}, 91(10):1700--1708, 2004.
\newblock http://www.amjbot.org/cgi/reprint/91/10/1700.pdf.

\bibitem{LSS13}
Simone Linz, Katherine St.~John, and Charles Semple.
\newblock Counting trees in a phylogenetic network is \#{P}-complete.
\newblock {\em SIAM J. Comput.}, 42(4):1768--1776, 2013.

\bibitem{MSW15}
Colin McDiarmid, Charles Semple, and Dominic Welsh.
\newblock Counting phylogenetic networks.
\newblock {\em Ann. Comb.}, 19(1):205--224, 2015.

\bibitem{Ro07}
Noah~A. Rosenberg.
\newblock Counting coalescent histories.
\newblock {\em J. Comput. Biol.}, 14(3):360--377, 2007.

\bibitem{Sc}
Ernst Schr\"oder.
\newblock Vier kombinatorische {P}robleme.
\newblock {\em Z. Math. Phys.}, 15:361--376, 1870.

\bibitem{Se16}
Charles Semple.
\newblock Phylogenetic networks with every embedded phylogenetic tree a base
  tree.
\newblock {\em Bull. Math. Biol.}, 78(1):132--137, 2016.

\bibitem{Se17}
Charles Semple.
\newblock Size of a phylogenetic network.
\newblock {\em Discrete Appl. Math.}, 217(part 2):362--367, 2017.

\bibitem{Semple2018}
Charles Semple and Jack Simpson.
\newblock When is a phylogenetic network simply an amalgamation of two trees?
\newblock {\em Bulletin of Mathematical Biology}, 80(9):2338--2348, Sep 2018.

\bibitem{SeSt06}
Charles Semple and Mike Steel.
\newblock {\em IEEE/ACM Trans. Comput. Biology Bioinform.}, 3:84--91, 2006.

\bibitem{10.2307/2412721}
P.~H.~A. Sneath.
\newblock Cladistic representation of reticulate evolution.
\newblock {\em Systematic Zoology}, 24(3):360--368, 1975.

\bibitem{tree-child}
Louxin Zhang.
\newblock Counting tree-child networks and their subclasses.
  arxiv:1908.01917v2.
\newblock 2019.

\bibitem{unknown}
Louxin Zhang.
\newblock Generating normal networks via leaf insertion and nearest neighbor
  interchange.
\newblock {\em BMC Bioinformatics}, 20(20):642, Dec 2019.

\end{thebibliography}
\newpage
\appendix
\section{General Phylogenetic Network With Three Reticulation Nodes}
\begin{center}
 	\begin{figure}[h]
 		\begin{minipage}{1\textwidth}
 			\begin{center}
 				{\includegraphics[width=1\textwidth]{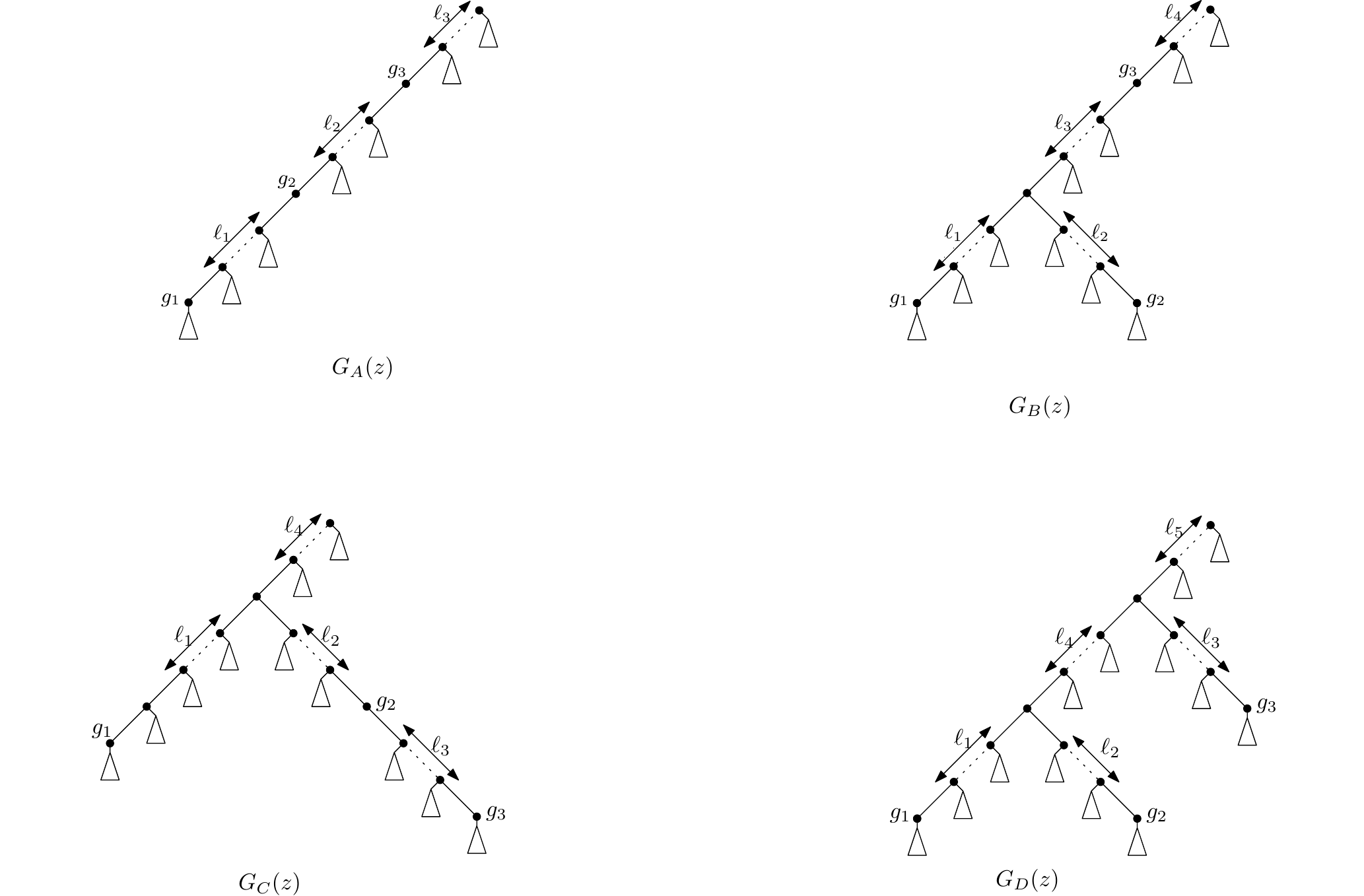}}
 			\end{center}
 		\end{minipage}
 		\caption{The structure of Motzkin skeletons of networks  with $ 3 $ reticulation vertices. All of them, originate from a sparsened skeleton which consists of only  green vertices.}
 		\label{Gen3}
 	\end{figure}
 \end{center}

This section presents a theoretical extension of the studied procedure for general phylogenetic network with three reticulation vertices. As like before,  we decompose the network
according to how the reticulation vertices are distributed in the networks. For more explicitly, first consider the Motzkin  skeletons with just green vertices ( Figure \ref{Gen3}). We can use them to figure out the rest of  Motzkin skeletons with red-green and double-green vertices as well. In the end, we add the contribution of the Motzkin skeletons with multiple edges.  
For $ i,\ j \in \{1,2,3,r,g\}  $,  use  $Y_{i,...,j}$ denotes the operator differentiating with respect to $y_i,...,\ y_j $ and setting $y_i=...=y_j=0$ afterwards, \emph{i.e.}, $\ Y_{i,...,j}
f(z,y_i,...,y_j)=\(\partial_{y_i}... \partial_{y_j}f\)(z,0,...,0)$. 
Now we investigate the details of extracting exponential generating function for cases in the Figure \ref{Gen3}. We follow the same procedure that used for general phylogenetic networks with two reticulation vertices. Start with simple case which that  the three green vertices lie on one path, \emph{i.e.}, one green
vertex is ancestor of another, which itself is ancestor of the third one. All possibilities for the pointings of the edges starting at $g_1, g_2$ and $g_3$  may target  any vertex in all the other trees. Concerning the vertices on the spine, we have some restrictions. The edge from $g_1$ may not end at any vertex from $\ell_1 $ , $ \ell_2 $ and the root of its attached subtree.  The edges from  $g_2$ may not point to first vertex of $ \ell_1 $ (to avoid of multiple edges)  and any vertex of $\ell_2$. Finally,
no green vertex may point to the vertex of $\ell_3$. Note that the contribution of  multiple edges will be considered in later cases. Overall, this yields the generating function
\begin{align*} 
G_{A}(z)&=Y_{1,2,3} \(\frac{z^3 \tilde M_1(z,y_1+y_2+y_3)}{1-zM(z,y_1+y_2+y_3)} P(z,y_3,y_1+y_2+y_3,0) P(z,y_2+y_3,y_1+y_2+y_3, y_3)\).
\end{align*} 
Next we will determine the generating function of all general networks belonging to case that  one green vertex is a common ancestor of the other two, but none of those two is ancestor of the other one.
As in the previous section we analyse the substructures. There are four vertices in the sparsened skeleton, yielding a factor $ z^4 $ .
Any non-root red vertices in the subtree attached to $g_1$ may be targets of the edge coming
from any green nodes and for root one, pointing is allowed for $ g_2 $ and $g_3$ (but not $ g_1 $,  to avoid multiple edges).  for the subtree attached to $ g_2 $ vice versa.

\begin{itemize}
	\item Paths $\ell_3$ and $\ell_4$: These paths are sequences of vertices, each with a subtree attached to it. For $\ell_4$ each green vertex is allowed to point at the red vertices in these  subtrees. Pointing to the vertices of the path is not allowed. Likewise, just the corresponding vertices on the path of
	$\ell_3$ are forbidden for $g_1$ and $g_2$ by the generality condition but $ g_3 $  may points non-first vertex of that as well.
	\item Paths $\ell_1$ and $\ell_2$: They are symmetric, so we discuss $\ell_1$. The  vertices of the subtrees are  allowed targets
	for the edge from all green vertices. The edge from $g_2$ and $g_3$ may end at each vertex of  the path. 
\end{itemize}
Note that the generality condition will be violated by making cyclic component, if $g_2$ points a red vetex on the path $ \ell_1 $ and $ g_1 $ vice versa. We subtract this cases from the result.
Overall, this gives, again using the operator $Y_{i,...,j}$ defined above, the generating function
\begin{align*} 
G_{B}(z)=&\rcp2Y_{1,2,3}\left(\frac{z^4 \tilde M_{1}(z,y_1+y_2+y_3)\tilde M_{2}(z,y_1+y_2+y_3)}{1-zM(z,y_1+y_2+y_3)}  P(z,y_1+y_3,y_1+y_2+y_3 ,y_1+y_3 ) \right.\quad\\
\times& P(z,y_3,y_1+y_2+y_3,0) P(z,y_3+y_2,y_1+y_2+y_3,y_3+y_2 ) \\
-&
\left.\frac{z^4 M(z,y_3)^2}{(1-zM(z,y_3))} P(z,y_1+y_3,y_3 ,y_1+y_3 ) P(z,y_2+y_3,y_3 ,y_2+y_3 )P(z,y_3,y_3 ,0 )\right)
\end{align*}	
Next we pay attention to the case that one green vertex is ancestor of another one, but not of both of them, and the third one is not ancestor of any other green vertex. The sparsened skeleton has $  4 $ vertices and the subtrees attached to $g_1$ and $g_3$. The red vertices of the subtree of $g_1$ and $g_3$   may be targeted by any edges starting from green vertices. Note that if $ g_1 $ and $ g_3 $ have the red root attached subtrees, they are not allowed to point at own attached red root vertex respectively to avoid multiple edges.
Next we inspect the paths:
\begin{itemize}
	\item Path $\ell_4$: All green vertices may point to the  vertices of the  subtrees. Pointing to the path itself is not allowed.
	\item Path $\ell_3$: The edge starting at $g_3$ may point to  vertices of the subtrees, but not to the vertices of the path itself.
	All but the first vertex for $g_2$  of the path as well as all tree vertices can be the end point of the edge starting at $g_1$ and $g_2$.
	\item Path $\ell_1$: Similar to $\ell_3$. The edges from $g_2$ and $g_3$ may point anywhere of the path. The vertices of the
	subtrees may be targeted by $g_1$ as well.
	\item Path $\ell_2$: All green vertices may point to the  vertices of the subtrees. To point at the vertices on the path  is only allowed for $g_1$.
\end{itemize}

Altogether, we obtain the generating function $G_{C}(z)$ with  the expression
{\fontsize{10}{13}\selectfont
	\begin{align*} 
	G_{C}(z)&=Y_{1,2,3}\left(\frac{z^4 \tilde M_{3}(z,y_1+y_2+y_3)\tilde M_{1}(z,y_1+y_2+y_3)}{1-zM(z,y_1+y_2+y_3)} P(z,y_2+y_3,y_1+y_2+y_3 ,y_2+y_3)\right. \quad\\
	&\quad P(z,y_1,y_1+y_2+y_3 ,y_1 )P(z,y_1+y_2,y_1+y_2+y_3 ,y_1 )\\
	& -\frac{z^4 M(z,0)^2}{1-zM(z,0)} P(z,y_1,0,y_1)^2
	P(z,y_2+y_3,0,y_2+y_3) \\
	&-\frac{z^4\tilde M_3(z,y_3)M_3(z,y_3)}{(1-zM(z,y_3))^2} P(z,y_1,y_3,y_1)
	P(z,y_2,y_3,y_2) \\
	&-\frac{z^4 M_2(z,y_2)^2}{1-zM(z,y_2)} P(z,y_1,y_2,y_1)
	P(z,y_1+y_2,y_2,y_1) P(z,y_3,y_2,y_3) \Bigg).
	\end{align*} }
In this way, Motzkin skeletons which are not respecting the generality condition are generated as well: Indeed, $g_1$ may point to the
vertex on the paths $\ell_2$  or $\ell_3$ when both or one of  $g_2$ and $ g_3 $ point to vertex of $\ell_1$, such that makes directed cyclic component.

\begin{center}
	\begin{figure}[h]
		\begin{minipage}{1\textwidth}
			\begin{center}
				{\includegraphics[width=1\textwidth]{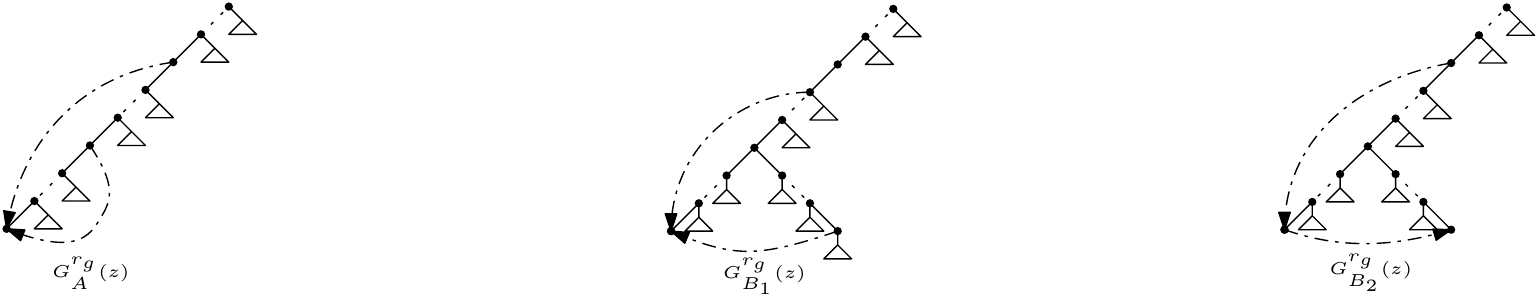}}
			\end{center}
		\end{minipage}
		\caption{The structure of  Motzkin skeletons with  red-green leaves that can be generated from the Motzkin skeletons $ G_A (z) $ and $ G_B (z) $. In the left and right figures,  $ g_1 $ is replaced with a red-green leaf, and in the right one two red-green leaves  are considered as a pointer vertex at the end of cherry.}	
		  		\label{redgreen1}
	\end{figure}
\end{center}

The last case of general networks has Motzkin skeletons as shown in Figure~\ref{Gen3}.
The restriction for the target vertex of the edges to be added at $g_1$, $g_2$ and $g_3$ follow the
analogous rules in order to meet the generality constraint. Setting up the generating function follows the same pattern as before.
We omit now the details and get from path analysis after all
{\fontsize{10}{13}\selectfont
	\begin{align*} 
	G_{D}(z)&=\rcp2
	Y_{1,2,3}\\
	&\left(\frac{z^5\tilde M_{3}(z,y_1+y_2+y_3) \tilde M_{2}(z,y_1+y_2+y_3) \tilde M_{1}(z,y_1+y_2+y_3)}{1-zM(z,y_1+y_2+y_3)}P(z,y_1+y_2,y_1+y_2+y_3 ,y_1+y_2) 
	\right. 
	\\
	&\times  P(z,y_1+y_3,y_1+y_2+y_3 ,y_1+y_3 ) P(z,y_2+y_3,y_1+y_2+y_3 ,y_2+y_3) P(z,y_3,y_1+y_2+y_3 ,y_3)
	\\ 
	&-\frac{z^5 \tilde M_{3}(z,y_3)M(z,y_3)^2}{(1-zM(z,y_3))^2} P(z,y_1+y_3,y_3 ,y_1+y_3 ) P(z,y_2+y_3,y_3 ,y_2+y_3 )P(z,y_3,y_3 ,y_3 )\Bigg)
	\\ 
	& -Y_{1,2,3} \left(
	\frac{z^5 \tilde M_{1}(z,y_3)M(z,y_1)^2}{(1-zM(z,y_1))^2} P(z,y_1+y_3,y_1 ,y_1+y_3 ) P(z,y_2+y_1,y_1 ,y_2+y_1 )P(z,y_3,y_1 ,y_3 )\right.
	\\ 
	& -
	\frac{z^5 M(z,0)^3}{(1-zM(z,0))^2} P(z,y_1,0,y_1 ) P(z,y_2,0 ,y_2 )P(z,y_3,0 ,y_3 ) \Bigg).
	\end{align*}}
So far we  just have considered the Motzkin skeletons in Figure \ref{Gen3}  with three reticulation vertices such that only green vertices are considered as  pointer set vertices.
Now we consider the structure of the Motzkin skeletons with red-green and double-green vertices and set up exponential generating functions for them separately. 
 Note that, the crucial point is that distribution of pointer nodes on the Motzkin skeleton must be such a way  that, after adding directed edges, we get  a general phylogenetic network with $ 3 $ reticulation vertices.
 Recall that, for any red-green leaf first we consider another pointer vertex such that connects to this nodes by adding a directed edge. 
Let's start with the Motzkin skeletons that contain at least one red-green vertex. Consider a case with three pointer vertices lie on a path (two green colored vertices with a red-green leaf), such a way that a red-green once lies on the bottom of the path (left of Figure \ref{redgreen1}).
Note that, we get two different expressions depends on our choice   that which green vertex ($ g_2 $ or $ g_3 $) is considered first to point red-green leaf.
\begin{align*} 
G_{A}^{r_g}(z)&=Y_{r,3}
\frac{z^3 P^{\star}(z, y_3, y_r+y_3,y_3) P(z, y_3,y_r+y_3,0) }{1-zM(z,y_r+y_3)}+Y_{r,2} \dfrac{z^3P(z, y_2,y_r+y_2,0)}{(1-zM(z,y_r+y_2))^2}.
\end{align*}

Note that to avoid of multiple edges, the path between $ g_2 $ and red-green vertex cannot be empty edge,  in the case of added directed edge connects $ g_2 $ to the red-green leaf. \\
As similar  before, there are two possible cases for  the general  networks arising from the Motzkin skeletons depicted on the middle
of Figure \ref{redgreen1}. 
In the first case, if we fix an  added  directed edge from $ g_2 $ to red-green leaf, the only restriction for pointing of $g_3$  will be the vertices on the path that connects it  to the root and its first child vertex (to avoid multiple edge). The red-green  vertex may point to any non-path vertex. The second term is regards the situation that a shortcut connects $ g_3  $ to the red-green vertex. After subtracting Motzkin skeletons which are not respecting general network condition, we obtain
\begin{align*} 
G_{B_1}^{r_g}(z)=&Y_{r,3}\left(\frac{z^4M(z,y_r+y_3)}{1-zM(z,y_r+y_3)}P(z,y_3,y_r+y_3,0) P(z,y_3,y_r+y_3,y_3)^2
\) \\
+&
Y_{r,2}\left(\frac{z^4\tilde M_2(z,y_2+y_r)}{(1-zM(z,y_r+y_2))^2}P(z,y_2,y_r+y_2,y_2) P(z,y_r,y_r+y_2,y_r)
\)
\\
-&
Y_{r,2}\left(\frac{z^4 M(z,0)}{(1-zM(z,0))^2}P(z,y_2,0,y_2) P(z,y_r,0,y_r)
\).
\end{align*}  
Another case such that one green vertex is a common ancestor of the other two red-green vertices, is  depicted in right of Figure \ref{redgreen1}. First, $ g_3 $ points to the one of red-green leaf then another directed edge  connects this leaf to second red-green leaf in the Motzkin skeleton. The edge starting at latter red-green leaf may point to any vertex except on the paths ones. This yields the generating function
\begin{align*} 
G_{B_2}^{r_g}(z)=&Y_r\left( \frac{z^4 }{(1-zM(z,y_r))^4}\right).
\end{align*}  
\begin{center}
	\begin{figure}[h]
		\begin{minipage}{1\textwidth}
			\begin{center}
				{\includegraphics[width=1\textwidth]{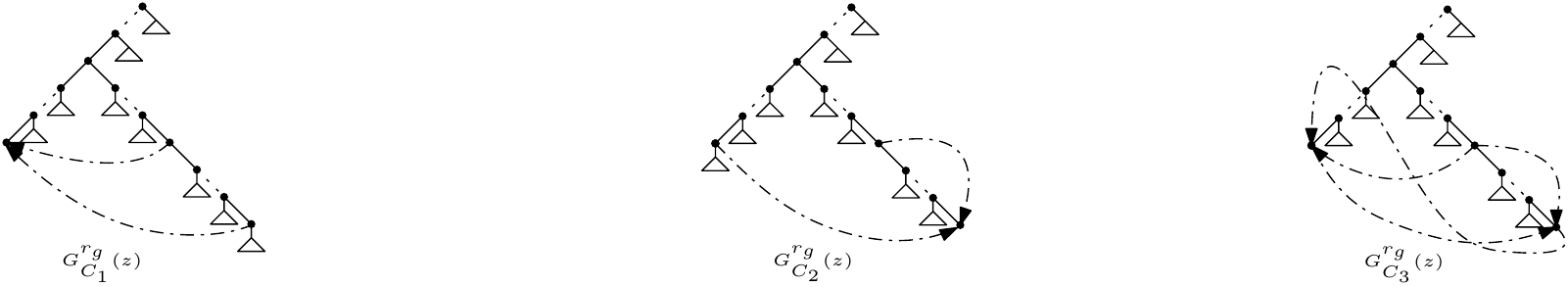}}
			\end{center}
		\end{minipage}
		\caption{The structures of Motzkin skeletons that are correspondent with the Motzkin skeleton $ G_C (z) $ by replacing respectively $ g_1 $,  $ g_3 $  or both of them  with red-green leaves.}
			\label{redgreen2}
	\end{figure}
\end{center}
Consider the Motzkin skeletons depicted in Figure \ref{redgreen2}. For the first one, the generating function is given by
\begin{align*} 
G_{C_1}^{r_g}(z)&=Y_{r,3}\left(\frac{z^4\tilde M_3(z,y_r+y_3)}{(1-zM(z,y_r+y_3))^2}P(z,y_r,y_r+y_3,y_r) P(z,y_3,y_r+y_3,y_3)\right.\\
&-\left.\frac{z^4 M(z,0)}{(1-zM(z,0))^2}P(z,y_r,0,y_r) P(z,y_3,0,y_3)
\right)
\\
&+
Y_{r,2}\left(\frac{z^4M(z,y_2+y_r)}{(1-zM(z,y_r+y_2))}P(z,y_2,y_r+y_2,y_2)  P(z,y_2,y_r+y_2,0)
\).
\end{align*}
For the Motzkin skeletons on the middle of Figure \ref{redgreen2}, we obtain 
\begin{align*} 
G_{C_2}^{r_g}(z)&=Y_{r,2}\left(\frac{z^4 M(z,y_r+y_2)}{(1-zM(z,y_r+y_2))^2}P(z,y_2,y_r+y_2,y_2) P(z,y_2,y_r+y_2,0)\)
\\
&+
Y_{r,1}\left(\frac{z^4 \tilde M_1(z,y_1+y_r)}{1-zM(z,y_r+y_1)}P(z,y_1,y_r+y_1,y_1) P^{\star}(z,y_1,y_r+y_1,y_1)  P(z,y_r,y_r+y_1,y_r)\right.\\
&
-\left.\frac{z^4 M(z,0)}{1-zM(z,0)}P(z,y_r,0,y_r) P(z,y_1,0,y_1)
P^{\star}(z,y_1,0,y_1)	\).
\end{align*}  
For the right one, we will take two terms for exponential generating function depending on which red-green leaf is pointed by $ g_2 $ first. After all, we get from path analysis 
\begin{align*} 
G_{C_3}^{r_g}(z)&=Y_r \left(	
\frac{z^4 }{(1-zM(z,y_r))^3}P^{\star}(z,0,y_r,0)+
\frac{z^4 }{(1-zM(z,y_r))^4}\right).
\end{align*}

\begin{center}
	\begin{figure}[h]
		\begin{minipage}{1\textwidth}
			\begin{center}
				{\includegraphics[width=.8\textwidth]{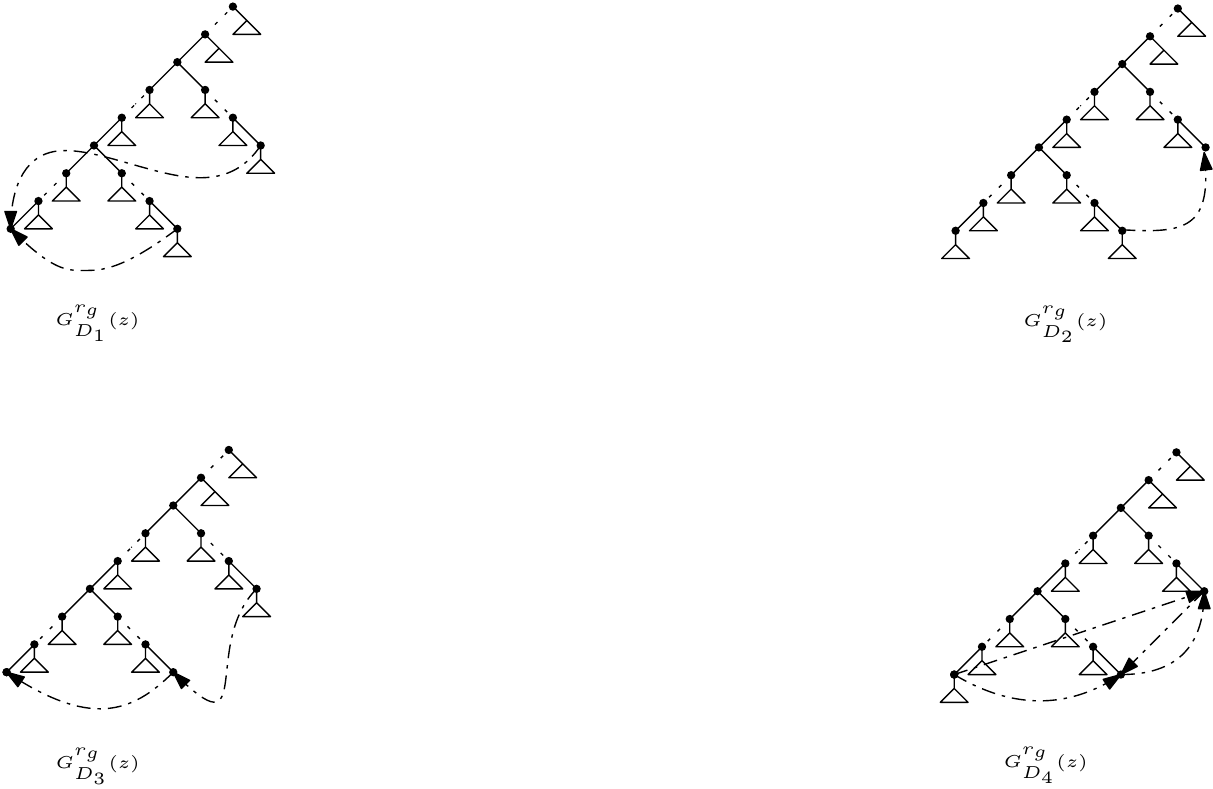}}
			\end{center}
		\end{minipage}
		\caption{ The Motzkin skeletons which arise from the $ G_D (z) $ by considering distribution of  all possible red-green leaves.}
		\label{redgreen3}
	\end{figure}
\end{center}
\begin{center}
	\begin{figure}[h]
		\begin{minipage}{1\textwidth}
			\begin{center}
				{\includegraphics[width=.8\textwidth]{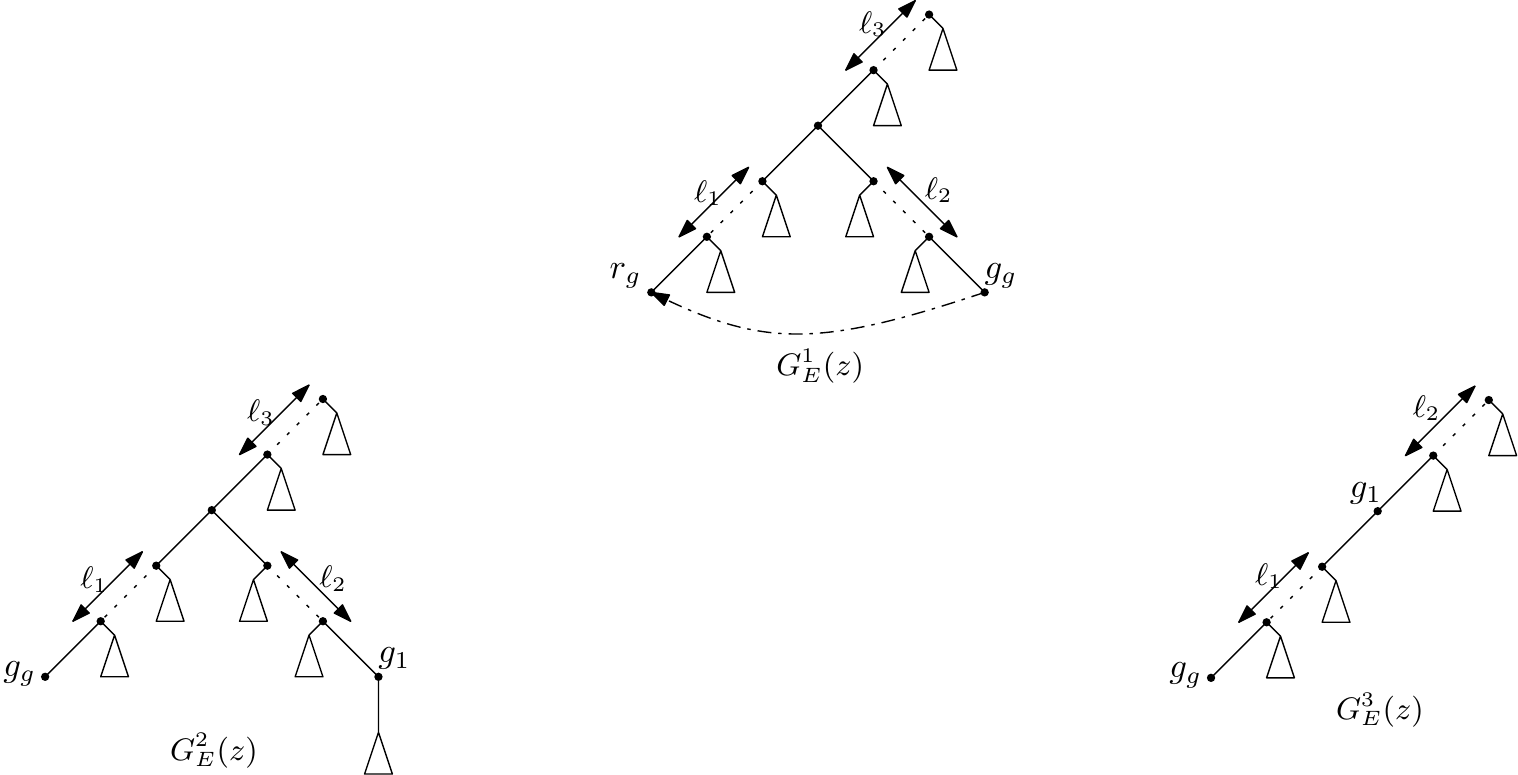}}
			\end{center}
		\end{minipage}
		\caption{Possible structures of Motzkin skeletons with three reticulation vertices and contribution of  double-green vertices.}
				\label{GE1M}
	\end{figure}
\end{center}

 
The last case of general networks with at least one red-green vertex have Motzkin skeletons as shown in Figure \ref{redgreen3}. The restriction
for the target vertex of the edges to be added at pointer set vertices follow the analogous rules in order
to meet the generality constraint. Setting up the generating function follows the same pattern as
before. We omit now the details and get from path analysis after all 
\begin{align*} 
G_{D_1}^{r_g}(z)&=Y_{r,3}\left(\frac{z^5\tilde M_3(z,y_r+y_3)M(z,y_r+y_3)}{1-zM(z,y_r+y_3)}P(z,y_r,y_r+y_3,y_r) P(z,y_3,y_r+y_3,y_3)^3
\right. \\
&-\left.\frac{z^5 M(z,0)^2}{(1-zM(z,0))}P(z,y_r,0,y_r) P(z,y_3,0,y_3)^3
\)
\\
&+
Y_{r,2}\left(\frac{z^5\tilde M_2(z,y_2+y_r)M(z,y_2+y_r)}{(1-zM(z,y_r+y_2))^2}P(z,y_2,y_r+y_2,y_2)^2 P(z,y_r,y_r+y_2,y_r)
\right.
\\
&-
\left.\frac{z^5 M(z,0)^2}{(1-zM(z,0))^2}P(z,y_2,0,y_2)^2 P(z,y_r,0,y_r)
\).
\end{align*}
\begin{align*} 
G_{D_2}^{r_g}(z)&=Y_{r,1}\left(\frac{z^5\tilde M_1(z,y_r+y_1)M(z,y_r+y_1)}{(1-zM(z,y_r+y_1))^2}P(z,y_r,y_r+y_1,y_r) P(z,y_1,y_r+y_1,y_1)^2
\right. \\
&-\left.\frac{z^5 M(z,0)^2}{(1-zM(z,0))^2}P(z,y_r,0,y_r) P(z,y_1,0,y_1)^2
\).\\ 
G_{D_3}^{r_g}(z)&=Y_r\left(\frac{z^5 M(z,y_r) }{(1-zM(z,y_r))^5}\right).
\\
G_{D_4}^{r_g}(z)&=Y_r\left(\frac{z^5 M(z,y_r) }{(1-zM(z,y_r))^5}+
\frac{z^5 M(z,y_r) }{(1-zM(z,y_r))^5}\right).
\end{align*}  
In the end, we consider the  Motzkin skeletons with the contribution double-green vertices as depicted in Figure \ref {GE1M}. Note that, The extra factor $ \rcp{2} $ appears in the expression of $ G_{E}^{2}(z) $ and $ G_{E}^{3}(z) $, because the order of pointing for double-green vertex  is not matter. After normalization we obtain 
\begin{align*} 
G_{E}^{1}(z)&= Y_{g,r}
\frac{z^3 P(z,y_g, y_r+y_g,y_g)}{(1-zM(z,y_g+y_r))^2}.
\\ 
G_{E}^{2}(z)&=\rcp{2} (Y_{g})^2 Y_{1} \frac{z^3\tilde M_1(z,y_1+y_g)}{1-zM(z,y_1+y_g)}P(z,y_g,y_1+y_g,y_g) P(z,y_1,y_g+y_r,y_1)\\
&-\rcp{2}(Y_{g})^2 Y_{1} \dfrac{z^3 M(z,0)}{1-zM(z,0)} P(z,y_1,0,y_1)P(z,y_g,0,y_g)\\
&-Y_{g} \dfrac{z^5 M(z,y_g)}{(1-zM(z,y_g))^5}.
\\ 
G_{E}^{3}(z)&=\rcp{2}(Y_{g})^2 Y_{1} \frac{z^2P(z,y_1,y_1+y_g,0) }{1-zM(z,y_1+y_g)}.
\end{align*}

Now, we  sum up  all obtained generating functions so far. For normalization, the result must be divided by $ 8 $, Since the procedure will
generate each general network eight times.
Overall, by collecting everything, 
the exponential generating function for vertex-labelled  general phylogenetic networks with
three reticulation nodes is 
\begin{equation*} 
{G}_3^{\nshortparallel}(z)= z \cdot \frac{a_3^{\nshortparallel}(z^2)-b_3^{\nshortparallel}(z^2)\sqrt{1-2z^2}}{(1-2z^2)^{11/2}},
\end{equation*}
where 
\begin{equation*} 
a^\nshortparallel_{3}(z)=3z^6+2z^5+4z^4+2z^3+\frac{69}{4}z^2,
\end{equation*}
and,
\begin{equation*} 
b^\nshortparallel_{3}(z)=z^5+\frac{9}{2}z^4+11z^3+\frac{69}{4}z^2.\;\;\;\;\;
\end{equation*}

\begin{center}
	\begin{figure}[h]
		\begin{minipage}{1\textwidth}
			\begin{center}
				{\includegraphics[width=.7\textwidth]{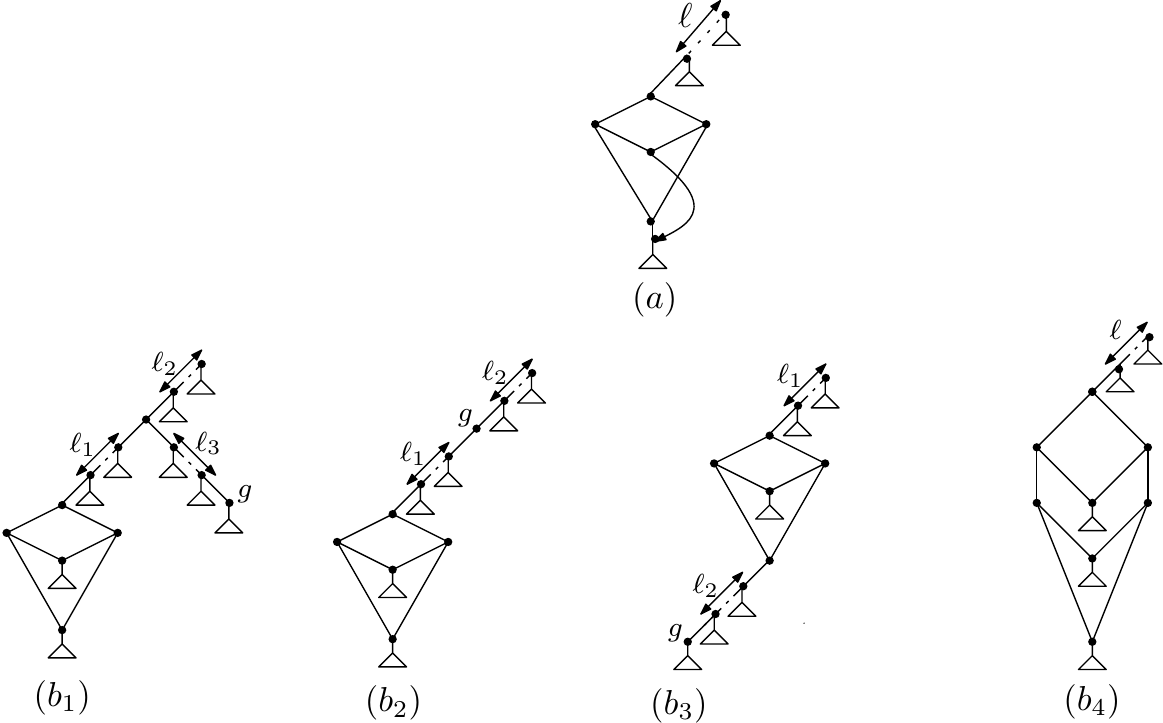}}
			\end{center}
		\end{minipage}
		\caption{The family  of general networks with generating function  $  {\ddot G^\nshortparallel}_3(z) $  such that each fixed leaf-labeled general network of (a) can construct corresponding  vertex-labeled general network  four times. For general network arose from second row shapes it will be  exactly two times. Note that, for $ b_1 $, $ b_2 $ and $ b_3 $ first we complete the structure by adding one more directed edge from green vertex to unary red vertex.   }
		\label{family1}
	\end{figure}
\end{center}

\begin{center}
	\begin{figure}
		\begin{minipage}{1\textwidth}
			\begin{center}
				{\includegraphics[width=.8\textwidth]{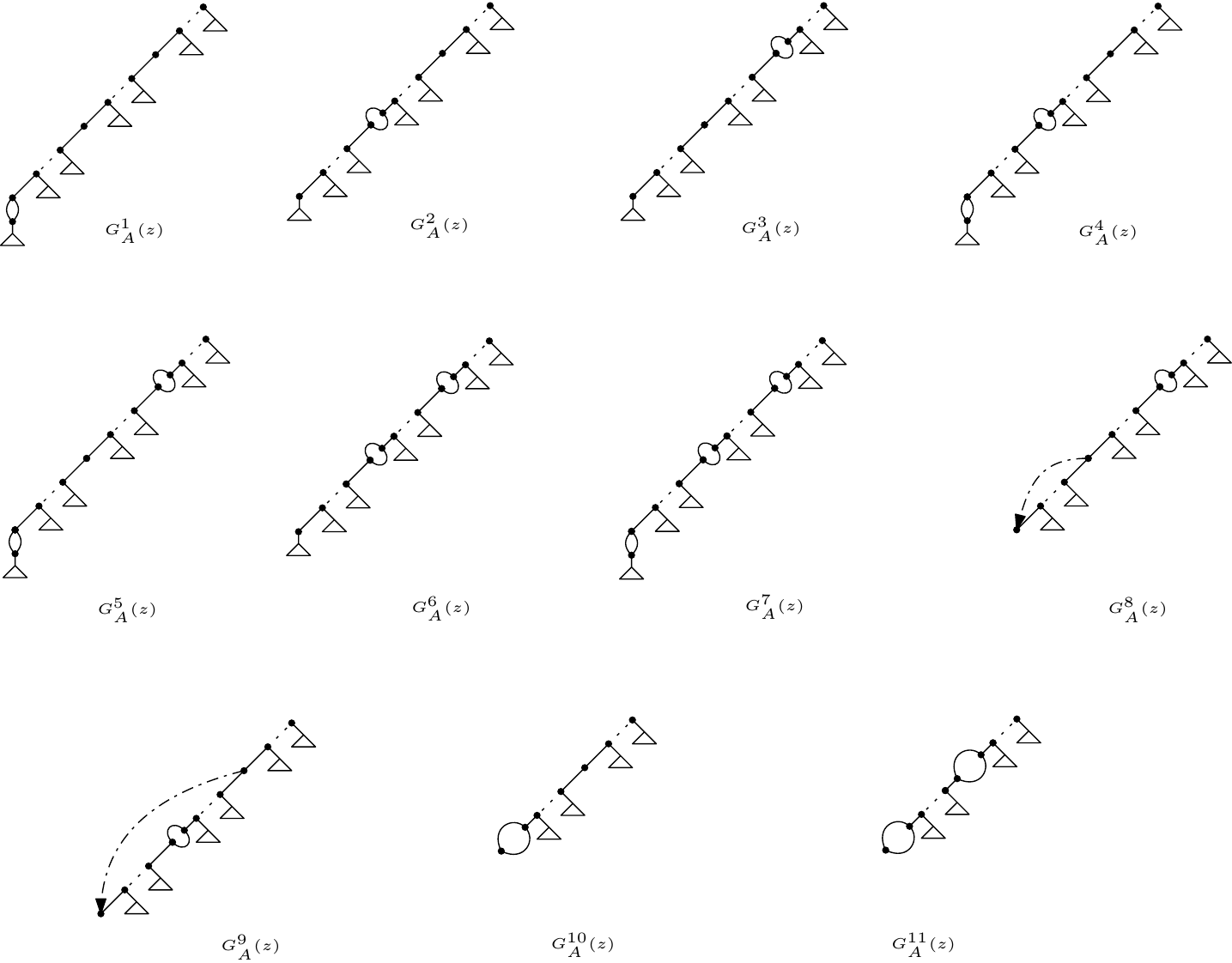}}
			\end{center}
		\end{minipage}
		\caption{The structures of the Motzkin skeletons of general phylogenetic networks with at least one multiple edges which are arised from $ G_A(z) $. }
		\label{G1}
	\end{figure}
\end{center}

\begin{center}
	\begin{figure}
		\begin{minipage}{1\textwidth}
			\begin{center}
				{\includegraphics[width=.7\textwidth]{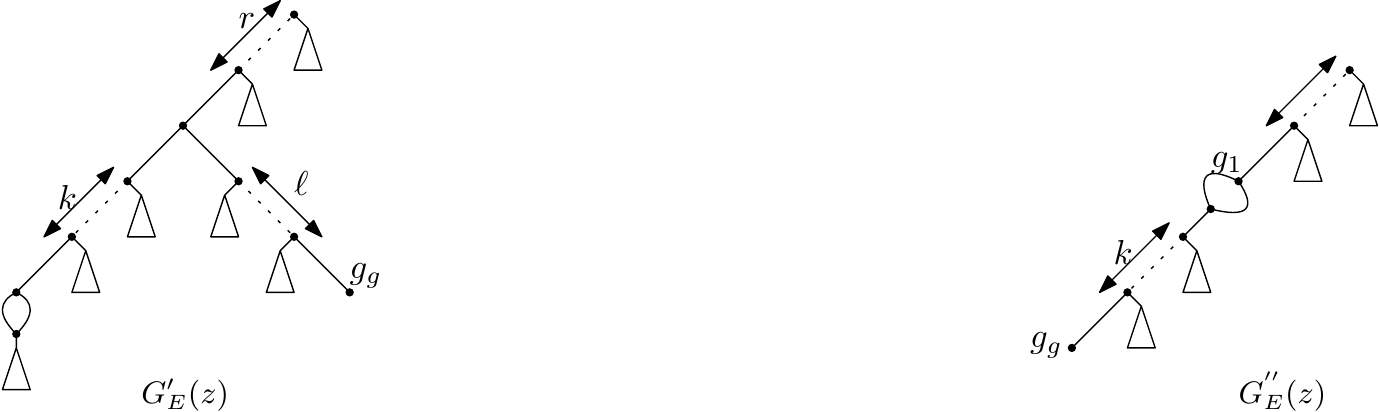}}
			\end{center}
		\end{minipage}
		\caption{	The structures of the Motzkin skeletons of general phylogenetic networks with at least one multiple edges which are arised from $ G_E^{2}(z) $ and $ G_E^{3}(z) $.}
				\label{GE}
	\end{figure}
\end{center}

Also,  consequently as similar as before we can take the explicit formulas for vertex and leaf-labeled general networks with $ 3 $ reticulation vertices. To see them set $n=2m+1$, so we have
\[
[z^n]{G}_3^\nshortparallel(z)=[z^m]\bar{G}_3^\nshortparallel(z),
\]
such that,
\[
[z^m]\bar G_{3}^\nshortparallel(z)=[z^m]\dfrac{a^\nshortparallel_{3}(z)}{(1-2z)^\frac{11}{2}}-[z^m]\dfrac{b^\nshortparallel_{3}(z)}{(1-2z)^5}.
\]
It gives 
\[
\mathcal{F^\nshortparallel}(m):=[z^m]\bar{G^\nshortparallel}_3(z)=\frac{2^{m-6}}{3}\Big(A^\nshortparallel_3(m)\frac{\displaystyle m (m-1)\binom{2m}{m}}{35(2m-1)4^{m-2}}-B^\nshortparallel_3(m)\Big).
\]
where 
\begin{align*} \label{polyPQG31}
A^\nshortparallel_3(m)&=104m^4+416m^3+596m^2-384m+61, \\
\quad B^\nshortparallel_3(m)&=48m^4+31m^3-12m^2-73m+6.
\end{align*}
By replacing $ m=(n-1)/2 $  we have,
$  G_{3,n}^\nshortparallel= n!\cdot \mathcal{F^\nshortparallel}((n-1)/2).
$\\
With some more steps but  similar as before we can present explicit formula for the number of leaf-labeled general network with three reticulation vertices. Let $ {\dot G^\nshortparallel}_3(z) $ denotes corresponding generating function 
for general networks that holds the situation of  equation (\ref{GVT}) and $ {\ddot G^\nshortparallel}_3(z) $ be generating function for general networks which arise from the Motzkin skeletons figure \ref{family1}.
 We have $ {G}_3^{\nshortparallel}(z)= {\dot G^\nshortparallel}_3(z) + {\ddot G^\nshortparallel}_3(z) $.
 So for the first subfamily (for $m>3$) we get
\begin{align}
\mathcal{\dot F^\nshortparallel}(m)=\frac{2^{m-5}}{3}\Big(\dot A^\nshortparallel_3(m)\frac{\displaystyle m (m-1)\binom{2m}{m}}{35(2n-5)(2n-3)(2m-1)4^{m-2}}-\dot B^\nshortparallel_3(m)\Big),     
\end{align} 
where,
\begin{align*} \label{polyPQG32}
\dot A^\nshortparallel_3(m)=280m^6-288m^5-1086m^4-2626m^3+9239m^2-7463m+4290, \quad 
\end{align*} 
  and
\begin{align*} 
\quad \dot B^\nshortparallel_3(m)=24m^4-\frac{31}{2}m^3+6m^2+\frac{85}{2}m-21.
\end{align*} 
Also for $ m=3 $ ($ \ell=1 $),  we have $\mathcal{\dot F^\nshortparallel}(3)=8 $.
Now we  consider the family  of general networks with $ 3 $ reticulation vertices such that there is a pair of vertices that have the set of same descendent and  applying the procedure (\ref{LtoV}) needs to cope with symmetry for them ; see Figure \ref{family1}. 
First, we set up generating function, let's show it $ G_{s_1}^{\nshortparallel}(z) $,  for case (a)
as shown at the top of Figure  \ref{family1}. Each fixed leaf-labeled general network which is arisen from this structure can generate corresponding vertex-labeled networks four times. 
So for this case we normalize  equation $ 6 $ by considering
$  {G^\nshortparallel}_{s_1,\ell}=4\frac{\ell!}{n!}G_{s_1,n}^{\nshortparallel}. $
%
Let's $ G_{s_2}^{\nshortparallel}(z) $ denotes the corresponding  generating function for second row structures of Figure  \ref{family1}.
 Note that each fixed leaf-labeled network belongs to these family can
 construct vertex-labeled network two times, so we get $  {G^\nshortparallel}_{s_2,\ell}=2\frac{\ell!}{n!}G_{s_2,n}^{\nshortparallel}. $ 
 Overall,  we obtain the $ {\ddot G^\nshortparallel}_3(z)= G_{s_1}^{\nshortparallel}(z)+G_{s_2}^{\nshortparallel}(z) $,
 where
\begin{align*} 
G_{s_1}^{\nshortparallel}(z)&=\rcp{4}
\frac{z^6 M(z,0)}{1-zM(z,0)},
\end{align*}
and then we get 
\begin{align}
\mathcal{\ddot F^\nshortparallel}_{s_1}(m):=[z^m]{\bar G^\nshortparallel}_{s_1}(z)=2^{m-2}\Big(\frac{\displaystyle m (m-1)\binom{2m}{m}}{(2n-3)(2m-1)4^{m}}\Big).
\end{align}
Also we have 
\begin{align*} 
G_{s_2}^{\nshortparallel}(z)&=\rcp{2} \partial y \frac{z^7\tilde M(z,y)M(z,y)^2}{4(1-zM(z,y))^2}P(z,y,y,y) \\
&+\rcp{2} \partial y \frac{z^6 M(z,y)^2}{4(1-zM(z,y))}P(z,y,y,0)\\
&+\rcp{2} \partial y \frac{z^6 \tilde M(z,y)M(z,y)}{2(1-zM(z,y))^2}
\\ 
&+\rcp{4} \frac{z^8M(z,0)^3 }{1-zM(z,0)}.
\end{align*}  

such that for or $ m>3 $ ( $ \ell>1 $) we have 
{\fontsize{10}{13}\selectfont
\begin{align}
	\mathcal{\ddot F^\nshortparallel}_{s_2}(m):=[z^m]{\bar G^\nshortparallel}_{s_2}(z)=2^{m-3}\Big((2m^3-15m^2+38m-34)\frac{\displaystyle m (m-1)\binom{2m}{m}}{(2m-5)(2n-3)(2m-1)4^{m-1}}-\rcp{2}(m-3)\Big),
\end{align}
	}
and  $ \mathcal{\ddot F^\nshortparallel}_{s_2}(3)=0 $. Obviously,  this means there is no any general network with one leaf which can be generated by second row structures of Figure \ref{family1}.  Overall, by collecting everything, we have $ G_{3,1}^\nshortparallel= 51$ and for $ \ell>1 $ we use 
\begin{align}
\begin{aligned}
  G_{3,\ell}^\nshortparallel&=
   \ell!\cdot \Big( \mathcal{\dot F^\nshortparallel}(\ell+2)+4\ddot F^{\nshortparallel}_{s_1}(\ell+2)+2\ddot F^{\nshortparallel}_{s_2}(\ell+2)\Big)\\
  &= \ell!\;  \Big(r_3(\ell)\displaystyle2^{-\ell}\binom{2\ell+4}{\ell+2}-2^{\ell}p_3(\ell)\Big),
  \end{aligned}   
\end{align}

where after manipulation  we get $ r_3(\ell) $ and $ p_3(\ell) $ as show in table (\ref{TablNumbers})
 for the number of leaf-labeled general networks with three reticulation vertices and no multiple edges. 
  In the following, we want to set up exponential generating functions for general networks with three reticulation vertices and at least one multiple edges.
It can be done by a case by case analysis of each sparsened skeletons which are depicted in Figures \ref{G1} to the \ref{GE}.
Note that, each factor of expression represents the number of generated times, so we use them to normalize final results. 
\newpage
\begin{align*}  
G_{A}^{1}(z)&=\rcp{4}Y_{2,3}
\frac{z^4 M(z,y_2+y_3)  }{1-zM(z,y_2+y_3)}P(z,y_2+y_3,y_2+y_3 ,y_3) 
P(z,y_3,y_2+y_3 ,0 ).
\\
G_{A}^{2}(z)&=\rcp{4}
Y_{1,3}
\frac{z^4 \tilde M_1(z,y_1+y_3)}{1-zM(z,y_1+y_3)}P(z,y_3,y_1+y_3 ,y_3) 
P(z,y_3,y_1+y_3 ,0 ).
\\
G_{A}^{3}(z)&=\rcp{4}Y_{1,2}
\frac{z^4 \tilde M_1(z,y_1+y_2)}{(1-zM(z,y_1+y_2))^2}P(z,y_2,y_1+y_2 ,0).\\
G_{A}^{4}(z)&=\rcp{2}Y_{3}
\frac{z^5  M(z,y_3) }{1-zM(z,y_3)}P(z,y_3,y_3,y_3)P(z,y_3,y_3,0).
\\
G_{A}^{5}(z)&=\rcp{2}Y_{2}
\frac{z^5  M(z,y_2) }{(1-zM(z,y_2))^2}P(z,y_2,y_2,0).
\\
G_{A}^{6}(z)&=
\rcp{2}Y_{1}
\frac{z^5  \tilde M(z,y_1) }{(1-zM(z,y_1))^3}.
\\
G_{A}^{7}(z)&=
\frac{z^6  M(z,0) }{(1-zM(z,0))^3}.
\\
G_{A}^{8}(z)&=\rcp{4} Y_{r}
\frac{z^4 P^{\star}(z,0, y_r,0)  }{(1-zM(z,y_r))^2}.\\
G_{A}^{9}(z)&=\rcp{4}Y_{r} \frac{z^4  }{(1-zM(z,y_r))^3}.\\
G_{A}^{10}(z)&=\rcp{4}Y_{r,3}
\frac{z^3 P(z,y_3, y_r+y_3,0)}{1-zM(z,y_3+y_r)}.\\
G_{A}^{11}(z)&=\rcp{2}Y_{r} \frac{z^4  }{(1-zM(z,y_r))^2}.
\end{align*} 

\begin{center}
	\begin{figure}
		\begin{minipage}{1\textwidth}
			\begin{center}
				{\includegraphics[width=1\textwidth]{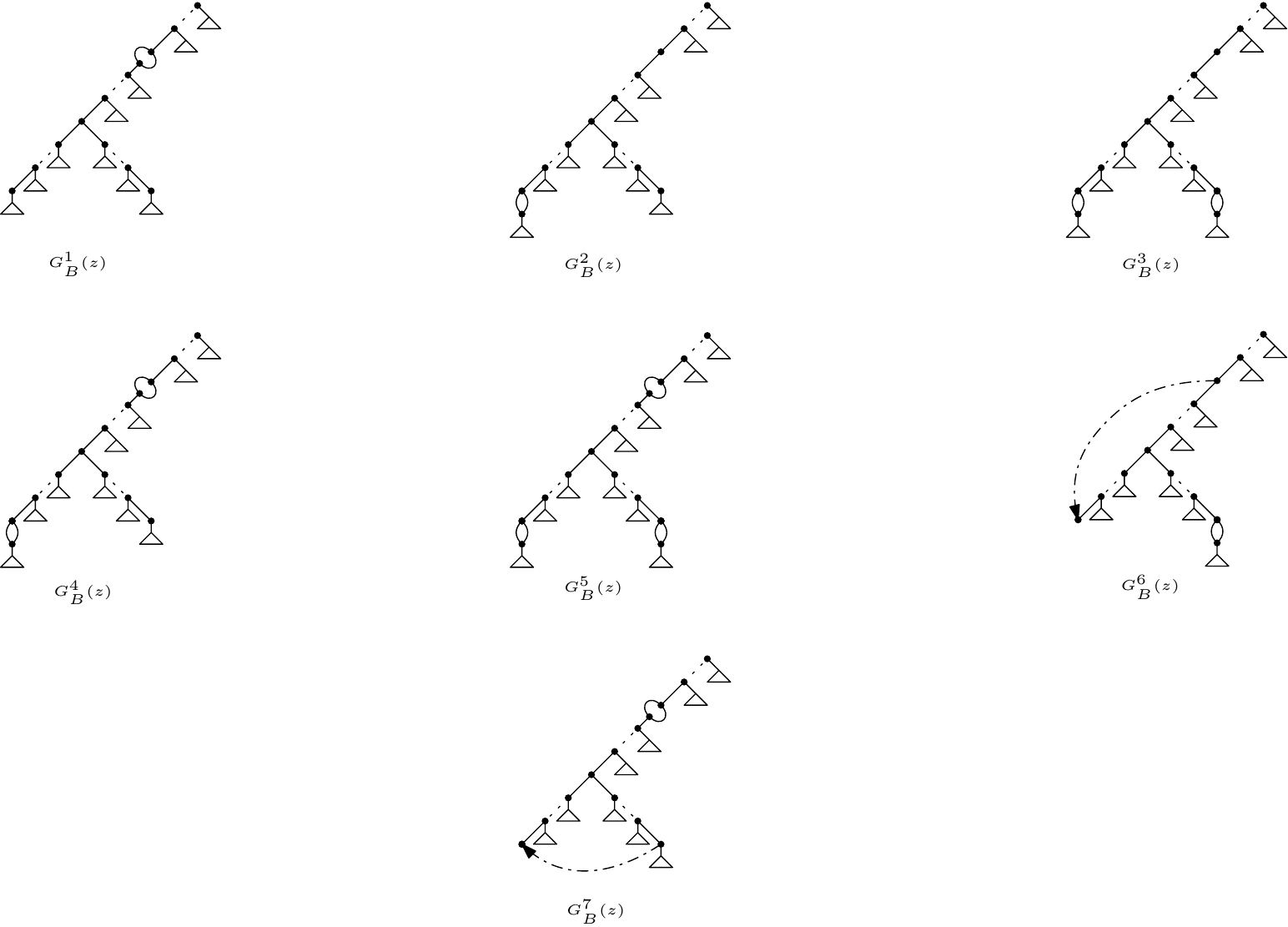}}
			\end{center}
		\end{minipage}
			\caption{	The structures of the Motzkin skeletons of general phylogenetic networks with at least one multiple edges which are arised from $ G_B(z) $.}
				\label{G2}
	\end{figure}
\end{center}
{\fontsize{10}{10}\selectfont
	\begin{align*}  
	G_{B}^{1}(z)=&\rcp8Y_{1,2}\Bigg(\frac{z^5 \tilde M_{1}(z,y_1+y_2)\tilde M_{2}(z,y_1+y_2)}{(1-zM(z,y_1+y_2))^2}  P(z,y_1,y_1+y_2 ,y_1 ) P(z,y_2,y_1+y_2,y_2) \\
	-&\frac{z^5 M(z,0)^2}{(1-zM(z,0))^2}P(z,y_1,0 ,y_1 ) 
	P(z,y_2,0,y_2)\Bigg).
	\\
	G_{B}^{2}(z)=&\rcp4Y_{2,3}\Bigg(
	\frac{z^5 M(z,y_2+y_3) \tilde M_{2}(z,y_2+y_3) }{1-zM(z,y_2+y_3)}P(z,y_2+y_3,y_2+y_3 ,y_2+y_3)\\ 
	\times	&  P(z,y_3,y_2+y_3 ,0 ) P(z,y_3,y_2+y_3 ,y_3)\Bigg). 
	\\
	G_{B}^{3}(z)=&\rcp4Y_{3}
	\frac{z^6 M(z,y_3)^2 }{1-zM(z,y_3)}P(z,y_3,y_3,y_3)^2P(z,y_3,y_3,0).\\
	G_{B}^{4}(z)=&\rcp2Y_{2}
	\frac{z^6 \tilde M_{2}(z,y_2) M(z,y_2) }{(1-zM(z,y_2))^3}P(z,y_2,y_2,y_2).\\	 	
	G_{B}^{5}(z)=&\rcp2
	\frac{z^7 M(z,0)^2 }{(1-zM(z,0))^4}.
	\\
	G_{B}^{6}(z)=&\rcp4Y_r
	\frac{z^5 M(z,y_r) }{(1-zM(z,y_r))^3}P(z,y_r,y_r,y_r).
	\\
	G_{B}^{7}(z)=&\rcp4Y_r\frac{z^5 M(z,y_r)}{(1-zM(z,y_r))^4}.	
	\end{align*} 
}
\begin{center}
	\begin{figure}
		\begin{minipage}{1\textwidth}
			\begin{center}
				{\includegraphics[width=.9\textwidth]{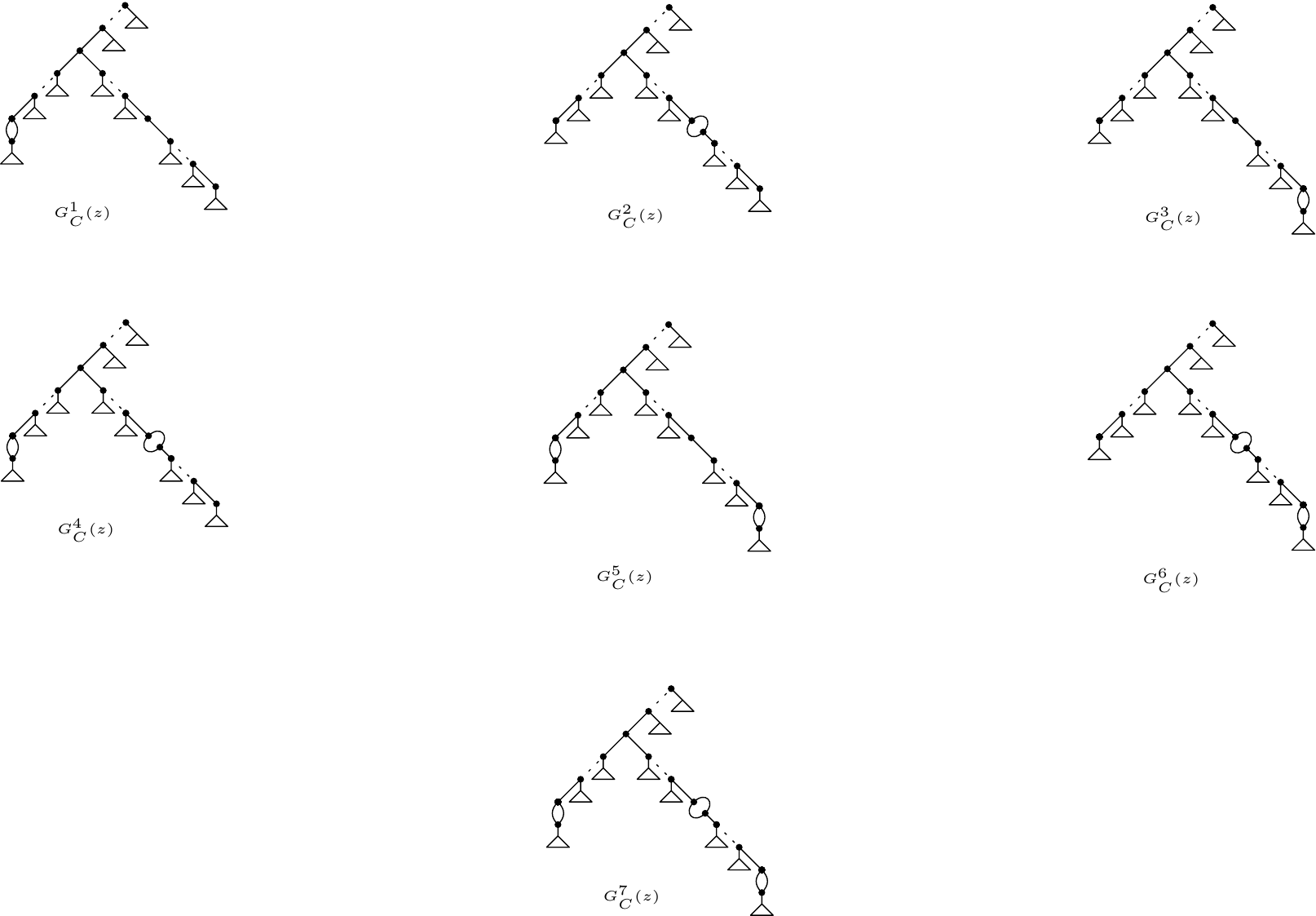}}
			\end{center}
		\end{minipage}
		\caption{	The structures of the Motzkin skeletons of general phylogenetic networks with at least one multiple edges which are arised from $ G_C(z) $.}
		\label{G3}
	\end{figure}
\end{center}  
\begin{align*}  
G_{C}^{1}(z)&=\rcp4Y_{2,3}
\frac{z^5 M(z,y_2+y_3) \tilde M_{3}(z,y_2+y_3) }{(1-zM(z,y_2+y_3))^2}P(z,y_2+y_3,y_2+y_3 ,y_2+y_3) 
P(z,y_2,y_2+y_3 ,0 ).
\\
G_{C}^2(z)&=\rcp4Y_{1,3}\left(\frac{z^5\tilde M_1(z,y_3+y_1)\tilde M_3(z,y_3+y_1)}{1-zM(z,y_3+y_1)}P(z,y_1,y_3+y_1,y_1)^2 P(z,y_3,y_3+y_1,y_3)\right.\\
&-\left.\frac{z^5 M(z,0)^2}{1-zM(z,0)}P(z,y_1,0,y_1)^2 P(z,y_3,0,y_3) \right). 
\\
G_{C}^3(z)&=\rcp4Y_{1,2}\left(\frac{z^5\tilde M_1(z,y_2+y_1) M(z,y_2+y_1)}{1-zM(z,y_2+y_1)}P(z,y_1,y_2+y_1,y_1) \right.\\ &\left.\times P(z,y_1+y_2,y_1+y_2,y_1)P(z,y_2,y_2+y_1,y_2)
\right. \\
&-\left.\frac{z^5 M(z,0)^2}{(1-zM(z,0))^2}P(z,y_1,0,y_1) P(z,y_2,0,y_2)
\).
\\
G_{C}^{4}(z)&=\rcp2Y_{3}
\frac{z^6 \tilde M_{3}(z,y_3) M(z,y_3) }{(1-zM(z,y_3))^3}P(z,y_3,y_3,y_3).
\\
G_{C}^{5}(z)&=\rcp2Y_{2}
\frac{z^6  M(z,y_2)^2 }{(1-zM(z,y_2))^2}P(z,y_2,y_2,y_2)P(z,y_2,y_2,0).
\\
G_{C}^{6}(z)&=\rcp2Y_{1}
\frac{z^6 \tilde M_{1}(z,y_1) M(z,y_1) }{(1-zM(z,y_1))^2}P(z,y_1,y_1,y_1)^2.
\\
G_{C}^{7}(z)&=\frac{z^7 M(z,0)^2 }{(1-zM(z,0))^4}.
\end{align*}  

\begin{center}
	\begin{figure}
		\begin{minipage}{1\textwidth}
			\begin{center}
				{\includegraphics[width=1\textwidth]{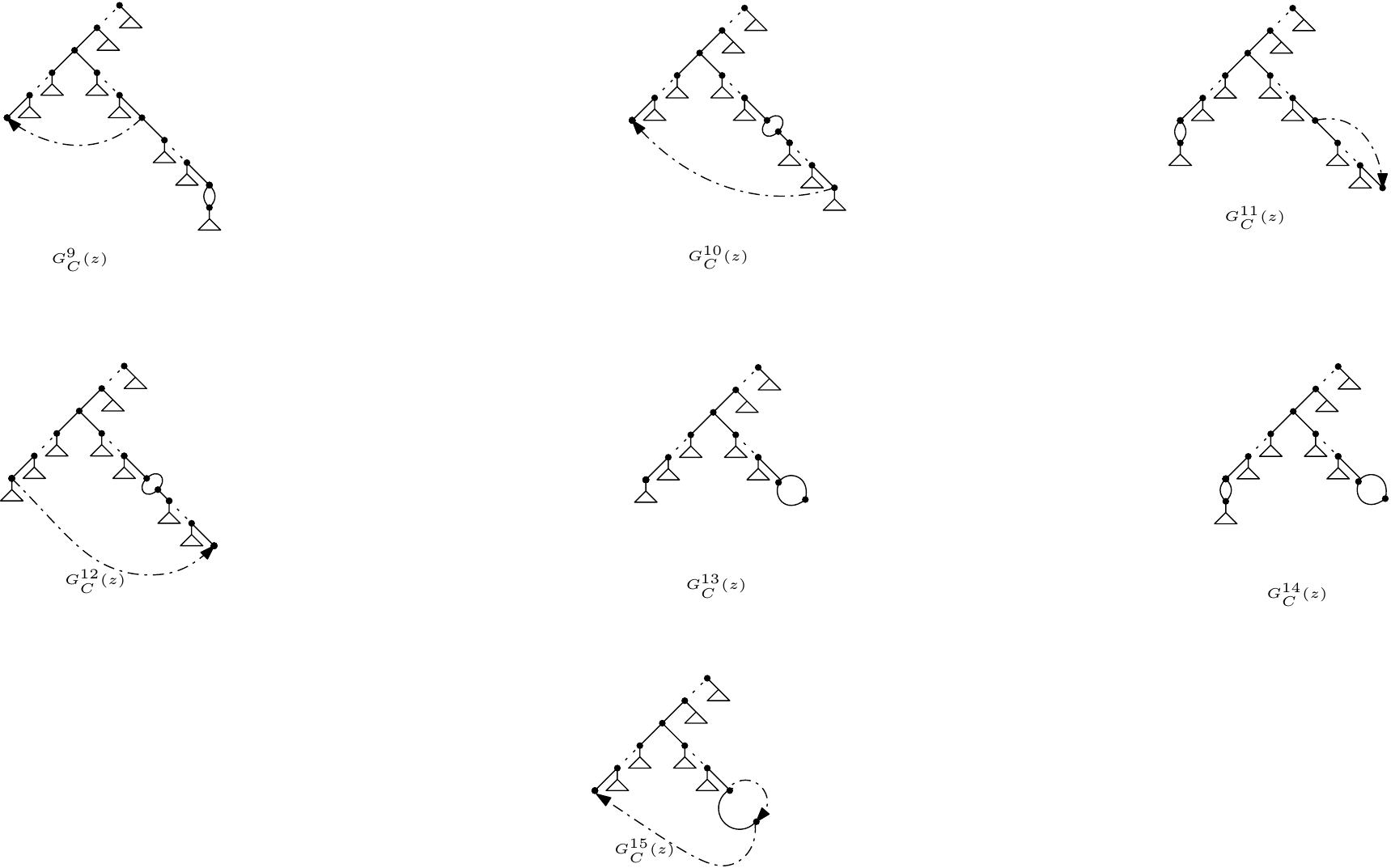}}
			\end{center}
		\end{minipage}
		\caption{	The structures of the Motzkin skeletons of general phylogenetic networks with at least one multiple edges which are arised from $ G_C(z) $.}
				\label{G3}
	\end{figure}
\end{center}
\begin{align*} 
G_{C}^{9}(z)&=\rcp4Y_r 
\frac{z^5 M(z,y_r)}{(1-zM(z,y_r))^3}P(z,y_r,y_r,y_r).
\\
G_{C}^{10}(z)&=\rcp4Y_r 
\frac{z^5 M(z,y_r)}{(1-zM(z,y_r))^4}.
\\
G_{C}^{11}(z)&=\rcp4Y_r
\frac{z^5 M(z,y_r)}{(1-zM(z,y_r))^2}P(z,y_r,y_r,y_r)P^{\star}(z,0,y_r,0).
\\
G_{C}^{12}(z)&=\rcp4Y_r\frac{z^5 M(z,y_r)}{(1-zM(z,y_r))^4}.
\\
G_{C}^{13}(z)&=\rcp4Y_{r,1}\left(
\frac{z^4 \tilde M_1(z,y_1+y_r) }{1-zM(z,y_1+y_r)}P(z,y_1,y_1+y_r,y_1)P(z,y_r,y_1+y_r,y_r)\right.\\
&-\left. \dfrac{M(z,0)}{1-zM(z,0)} P(z,y_1,0,y_1)P(z,y_r,0,y_r)\).
\\
G_{C}^{14}(z)&=\rcp2Y_r \left(
\frac{z^5 M(z,y_r)}{(1-zM(z,y_r))^2}P(z,y_r,y_r,y_r)\right).\\
G_{C}^{15}(z)&=\rcp4Y_r 
\frac{z^4 }{(1-zM(z,y_r))^3}.
\end{align*}    

\begin{center}
	\begin{figure}
		\begin{minipage}{1\textwidth}
			\begin{center}
				{\includegraphics[width=.9\textwidth]{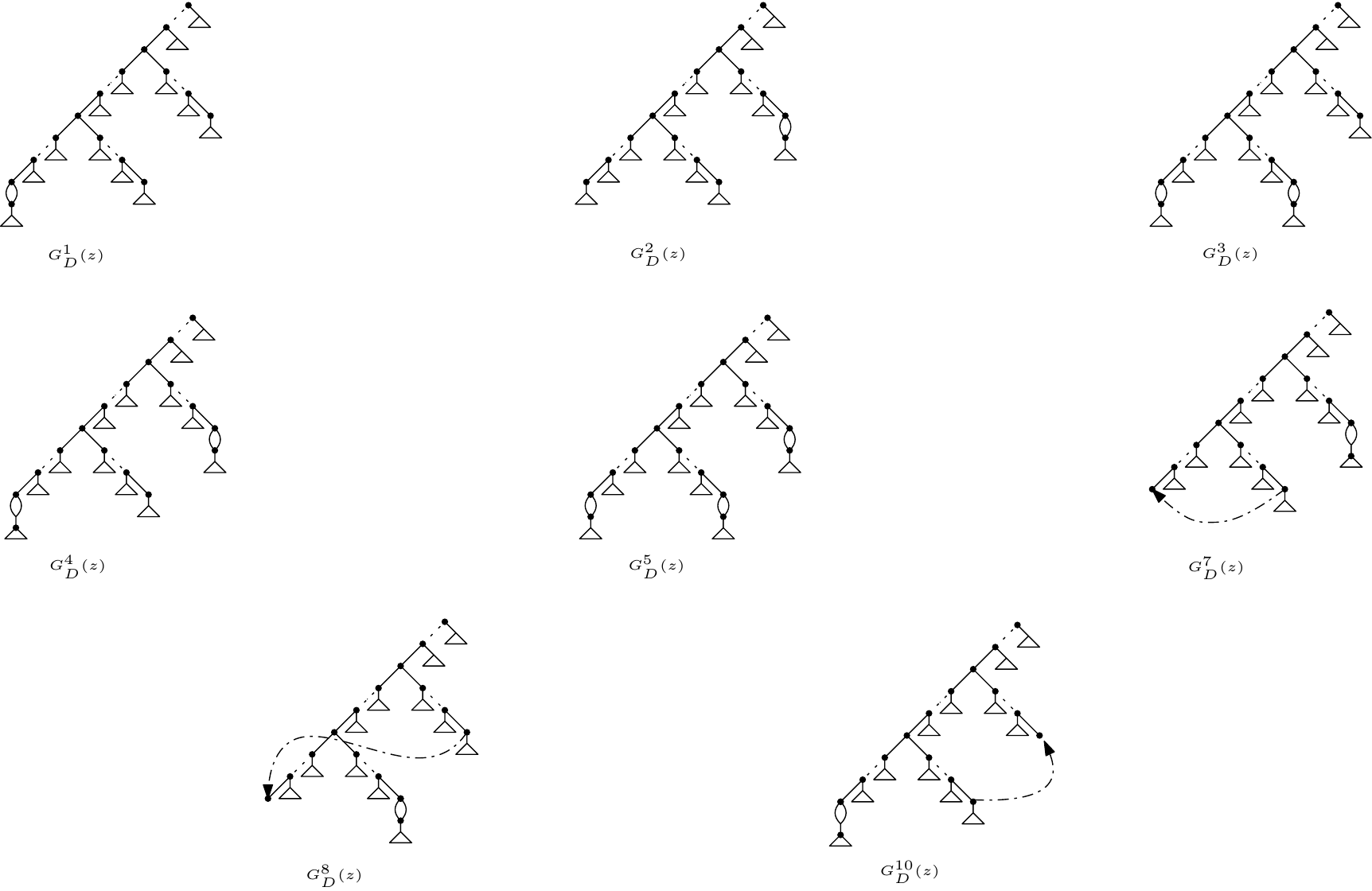}}
			\end{center}
		\end{minipage}
		\caption{	The structures of the Motzkin skeletons of general phylogenetic networks with at least one multiple edges which are arised from $ G_D(z) $.}
		\label{G4}
	\end{figure}
\end{center}  
{\fontsize{10}{10}\selectfont
	\begin{align*}  
	G_{D}^{1}(z)&=\rcp4Y_{2,3}
	\left(\frac{z^6 M(z,y_2+y_3) \tilde M_{2}(z,y_2+y_3) \tilde M_{3}(z,y_2+y_3)}{1-zM(z,y_2+y_3)}P(z,y_2+y_3,y_2+y_3 ,y_2+y_3) 
	\right. 
	\\
	&\times  P(z,y_2,y_2+y_3 ,y_2 ) P(z,y_3,y_2+y_3 ,y_3)^2 
	-
	\frac{z^6  M(z,0)^3}{(1-zM(z,0))^2} P(z,y_2,0 ,y_2 ) P(z,y_3,0 ,y_3 )^2\Bigg).\\
	G_{D}^{2}(z)&=\rcp8	Y_{1,2}
	\left(\frac{z^6 M(z,y_1+y_2) \tilde M_{1}(z,y_1+y_2) \tilde M_{2}(z,y_1+y_2)}{(1-zM(z,y_1+y_2))^2}P(z,y_1+y_2,y_1+y_2 ,y_1+y_2) 
	\right. 
	\\
	&\times  P(z,y_1,y_1+y_2 ,y_1 ) P(z,y_2,y_1+y_2 ,y_2) 
	-\frac{z^6  M(z,0)^3}{(1-zM(z,0))^3} P(z,y_1,0 ,y_1 ) P(z,y_2,0 ,y_2)\Bigg).\\
	G_{D}^{3}(z)&=\rcp4Y_{3}
	\frac{z^7 M(z,y_3)^2 \tilde M_{3}(z,y_3) }{(1-zM(z,y_3))^2}P(z,y_3,y_3,y_3)^3. \\
	G_{D}^4(z)&=\rcp2Y_{2}
	\frac{z^7 M(z,y_2)^2 \tilde M_{2}(z,y_2) }{(1-zM(z,y_2))^3}P(z,y_2,y_2,y_2)^2.
	\\
	G_{D}^5(z)&=\rcp2
	\frac{z^8 M(z,0)^3 }{(1-zM(z,0))^5}.
	\\
	G_{D}^6(z)&=\rcp4 \left(Y_r
	\frac{z^6 M(z,y_r)^2 }{(1-zM(z,y_r))^4}P(z,y_r,y_r,y_r)\right).
	\\
	G_{D}^{7}(z)&=\rcp4 \left(Y_r
	\frac{z^6 M(z,y_r)^2 }{(1-zM(z,y_r))^4}P(z,y_r,y_r,y_r)\right).
	\\
	G_{D}^{8}(z)&=\rcp{4}Y_r\left(
	\frac{z^6 M(z,y_r)^2 }{(1-zM(z,y_r))^4}P(z,y_r,y_r,y_r)\right). 	
	\end{align*}  
}

\begin{align*} 
G_{E}^{\prime}(z)&=\rcp{8} (Y_{g})^2  \frac{z^4 M(z,y_g)}{(1-zM(z,y_g))^2}P(z,y_g,y_g,y_g). \\ 
G_{E}^{''}(z)&=\rcp{8} (Y_{g})^2 \frac{z^3 }{(1-zM(z,y_g))^2}. \\
\end{align*}
\begin{center}
	\begin{figure}[h]
		\begin{minipage}{1\textwidth}
			\begin{center}
				{\includegraphics[width=.7\textwidth]{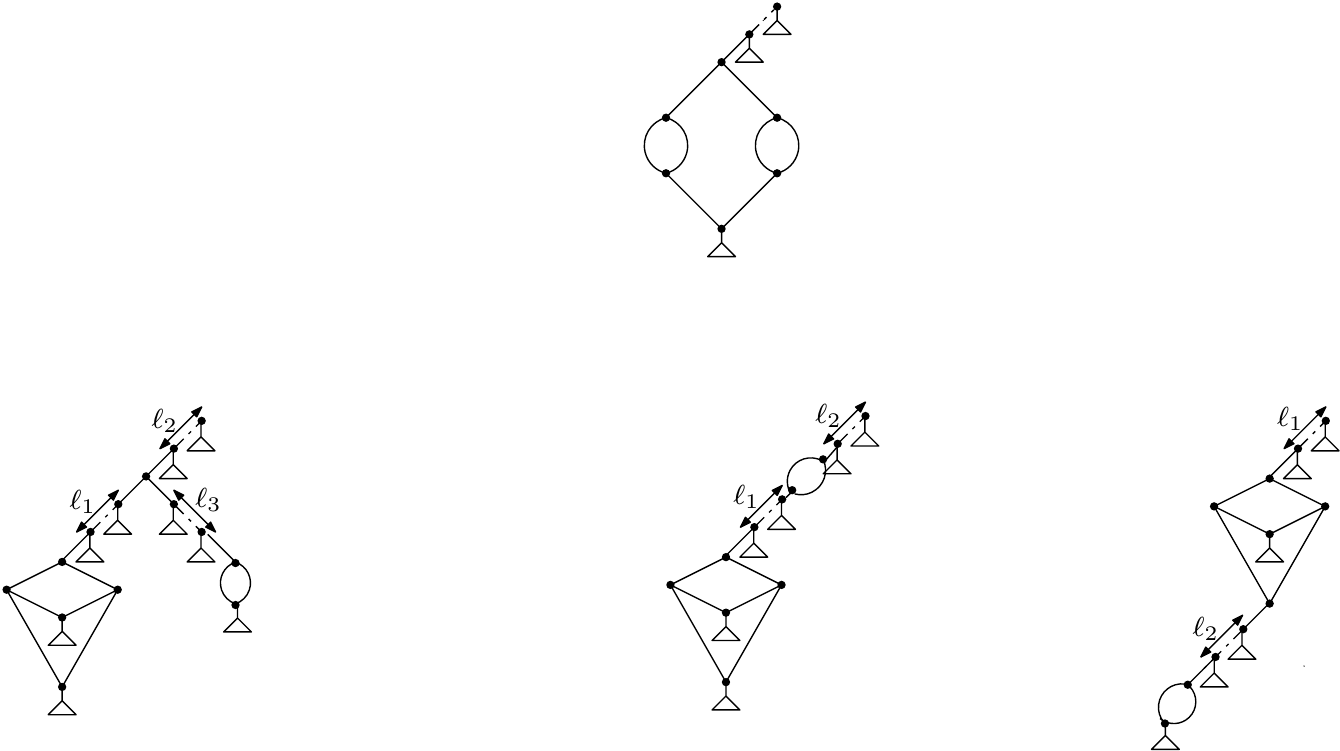}}
			\end{center}
		\end{minipage}
		\caption{General networks with multiple edges and corresponding generating function $ {\ddot G^\shortparallel}_3(z) $ such that any fixed leaf-labeled of them can generate all-vertex labeled exactly twice.}
				\label{family2}
	\end{figure}
\end{center} 
Overall, by collecting everything, we obtain the following result.  
\begin{equation*} 
{G}_3^{\shortparallel}(z)= z \cdot \frac{a_3^{\shortparallel}(z^2)-b_3^{\shortparallel}(z^2)\sqrt{1-2z^2}}{(1-2z^2)^{11/2}},
\end{equation*}
where 
\begin{equation*} 
a^\shortparallel_{3}(z)=z^5-z^4+\frac{13}{2}z^3+10z^2 \quad
\text{and},\quad
b^\shortparallel_{3}(z)=4z^3+10z^2.\;\;\;\;\;
\end{equation*}   
After some computation, it gives 
\begin{align}
\mathcal{F^\shortparallel}(m):=[z^m]\bar{G^\shortparallel}_3(z)=2^{m-1}\Big(A^\shortparallel_3(m)\frac{\displaystyle m (m-1)\binom{2m}{m}}{3(2m-1)4^{m}}-B^\shortparallel_3(m)\Big).
\end{align}
where 
\[ \label{polyPQG3p}
A^\shortparallel_3(m)=6m^3+4m^2-m-2 \quad \text{and,}\quad B^\shortparallel_3(m)=m^3-\rcp{2} m^2-\rcp{2}.
\]
By replacing $ m=(n-1)/2 $  we have,
$  G_{3,n}^\shortparallel= n!\cdot \mathcal{F^\shortparallel}((n-1)/2)
$ for the number of vertex-labeled general phylogenetics with $ 3 $ reticulation vertices and at least one multiple edge in their structures.
 
Now we set up generating function for leaf-labeled. we Consider
$ {G^\shortparallel}_3(z)= {\dot G^\shortparallel}_3(z)+{\ddot G^\shortparallel}_3(z)$ which respectively right side of equation  denote generating functions for two subfamilies of this class (general networks with multiple edges) that we can use the equation directly or not (needs to cope with symmetry); see Figure \ref{family2}. For the first subfamily we get
\[
\mathcal{\dot F^\shortparallel}(m):=[z^m]\dot{ \bar G}^\nshortparallel_3(z)=2^{m-2}\Big(\dot A^\nshortparallel_3(m)\frac{\displaystyle m (m-1)\binom{2m}{m}}{3(2n-3)(2m-1)4^{m-1}}-\dot B^\nshortparallel_3(m)\Big).
\]

where,
\begin{align} \label{polyPQG3pp}
\dot A^\shortparallel_3(m)=6m^4-5m^3-7m^2-2m+6 \quad \text{and,}\quad \dot B^\shortparallel_3(m)=2m^3-m^2-m.
\end{align}    
Also the generating function corresponding to the general networks in Figure \ref{family2} is   
\begin{align*}  
\ddot  G^\shortparallel_3(z)&=\rcp{2}  \frac{z^6 M(z,0)}{1-zM(z,0)}+\rcp{4}  \frac{z^8 M(z,0)^3}{(1-zM(z,0))^3}+\rcp{4}  \frac{z^7 M(z,0)^2}{(1-zM(z,0))^2}+\rcp{2}  \frac{z^7 M(z,0)^2}{(1-zM(z,0))^2} \\ 
&=\rcp{2}  \frac{z^3 }{(1-2z^2)^\frac{3}{2}}, \\
\end{align*}   
such that
\begin{align}
\mathcal{\ddot F^\shortparallel}(m):=[z^m]\ddot{\bar G}^\nshortparallel_3(z)= 2^{m-1}m (m-1)(m-2)\displaystyle\Big(\frac{\displaystyle \binom{2m}{m}}{(2n-3)(2m-1)4^{m}}\Big).
\end{align}    
 Note that, every member of leaf-labeled general networks arising from Figure \ref{family2} construct corresponding  vertex-labeled networks twice.
Overall, by replacing $ m=\ell+2 $  we have
\begin{align}
\begin{aligned} 
 G_{3,\ell}^\shortparallel&=  \displaystyle \ell!\cdot \Big(\mathcal{\dot F^\shortparallel}(\ell+2)+2\mathcal{\ddot F^\shortparallel}(\ell+2)\Big)\\
 &= \ell!\cdot 2^{\ell}\cdot  \displaystyle \Big(\dfrac{(\ell+1)(\ell+2)^2 (6\ell^3+31\ell^2+45\ell+15)\displaystyle\binom{2\ell+4}{\ell+2}}{3(2\ell+1)(2\ell+3)4^{\ell+1}}- (2\ell^3+11\ell^2+19\ell+10) \Big),
 \end{aligned}
\end{align} 
for the number of leaf-labeled general networks with three reticulation vertices and at least one multiple edge. Finally, we have $ \tilde{G}_{3,\ell}= G_{3,\ell}^\nshortparallel+G_{3,\ell}^\shortparallel $, for the number of all general phylogenetic networks with three reticulation vertices.


\end{document}